\documentclass{amsart}
\usepackage[utf8]{inputenc}
\usepackage{amssymb,amsfonts,amsthm,amssymb,amscd,amsmath}
\usepackage[cm]{fullpage}
\usepackage{wrapfig}
\usepackage{mathtools}
\usepackage{graphicx}
\usepackage{caption}
\usepackage{subcaption}
\usepackage{float}
\usepackage{yfonts}
\usepackage{mathrsfs}
\usepackage{color}
\usepackage[table]{xcolor}
\usepackage{enumitem}
\usepackage{hyperref}
\usepackage{setspace}
\usepackage{xfrac}
\usepackage[T1]{fontenc}
\usepackage{amsrefs}

\newtheorem{Lemma}{Lemma}[section]
\newtheorem{Theorem}{Theorem}
\newtheorem{Proposition}[Lemma]{Proposition}
\newtheorem{Corollary}[Lemma]{Corollary}

\DeclareMathAlphabet{\mathpzc}{OT1}{pzc}{m}{it}

\DeclarePairedDelimiter\floor{\lfloor}{\rfloor}

\newcommand{\R}{\mathbb{R}}

\newcommand{\tl}[1]{\tilde{#1}}

\newcommand{\beq}{\begin{equation}}
\newcommand{\eeq}{\end{equation}}
\newcommand{\ba}{\begin{align}}
\newcommand{\ea}{\end{align}}

\newcommand{\re}{\mathrm{e}}

\begin{document}

\title{\bf Asymptotic approximation of a modified compressible Navier-Stokes system }

\author[F1.Goh]{Ryan Goh}
\author[F2.Wayne]{C. Eugene Wayne}
\author[F3.Welter]{Roland Welter}
\address{Boston University \\
Department of Mathematics and Statistics \\
111 Cummington Mall \\
Boston, MA 02215}
\email{rgoh@bu.edu, cew@bu.edu, rwelter@bu.edu}
\thanks{\textbf{Acknowledgments:} The authors would like to acknowledge partial support by the National Science Foundation through grants DMS-1603416, DMS-2006887 (RG), DMS-1813384 (CEW) and (RW).  CEW thanks Th. Gallay for numerous discussions of the use of the Hermite expansion in the analysis of the asymptotic behavior of fluids.}

\begin{abstract}
	We study the long time asymptotics of a modified compressible Navier-Stokes system (mcNS) inspired by the previous work of Hoff and Zumbrun \cite{hoffzumbrun1}.  We introduce a new decomposition of the momentum field into its irrotational and incompressible parts, and a new method for approximating solutions of the heat equation in terms of Hermite functions in which $n^{th}$ order approximations can be computed for solutions with $n^{th}$ order moments.  We then obtain existence of solutions to the mcNS system and show that the approximation in terms of Hermite functions gives the leading order terms in the long-time asymptotics, and under certain assumptions can be evaluated explicitly.
\end{abstract}

\maketitle

\large

\section{Introduction}

The compressible Navier-Stokes equations are given by
\begin{equation} \label{CompNavStokes} \begin{split}
\partial_t \rho + \nabla \cdot \vec{m} & = 0 \\
\partial_t \vec{m} + \Big [ \nabla \cdot \left ( \frac{\vec{m}\otimes \vec{m}}{\rho} \right ) \Big ]^{T} + \nabla P & = \epsilon \Delta \left (\frac{\vec{m}}{\rho} \right ) + \eta \nabla \left ( \nabla \cdot \left ( \frac{\vec{m}}{\rho} \right ) \right ) 
\end{split} \end{equation}
These equations model the flow of a fluid with density $\rho$, momentum $\vec{m}$ and pressure $P$.  We assume that the fluid is barotropic, hence $P = P(\rho)$ is a function only of the density.  In the present paper, we are motivated by the question of stability of the constant density, constant momentum solution $(\rho^*,\vec{m}^*)^T$ to the compressible Navier-Stokes system in three dimensions, which without loss of generality we can take $(\rho^*,\vec{m}^*)^T=(1,0)^T$.  

Kawashima appears to have been the first to partially answer this question in dimension $d \geq 1$.  In \cite{kawashima}, he proves existence of global solutions for a general class of hyperbolic-parabolic systems which include (\ref{CompNavStokes}) and proves these solutions decay in $L^p$ at a given rate for $p \geq 2$.  In the course of his analysis, he studies a system with an artificial viscosity term, then shows that the results hold in the limit that the viscosity term is taken to zero. 

Hoff and Zumbrun (\cite{hoffzumbrun1}, \cite{hoffzumbrun2}) study the asymptotic behavior of small perturbations from the constant state for the compressible Navier Stokes equations.  Following Kawashima, they prove existence of solutions $u(t) = (\rho(t),m(t))^T$ for initial data $u_0 \in L^1 \cap H^s$ chosen sufficiently small for some $s \geq [\frac{d}{2}] +1$, and find that the solutions decay as $\| u(\cdot,t) \|_{L^p} \leq C t^{-\frac{d}{2}(1-\frac{1}{p})} $ for $p \geq 2$.   They show that there exists a unique linear, artificial-viscosity system associated with $(\ref{CompNavStokes})$ given by
\begin{equation} \label{LinArtViscSys} \begin{split}
\partial_t \rho + \nabla \cdot \vec{m} & = \frac{1}{2} \big ( \epsilon + \eta \big ) \Delta \rho \\
\partial_t \vec{m} + c^2 \nabla \rho & = \epsilon \Delta \vec{m} + \frac{1}{2} \big ( \eta - \epsilon \big ) \nabla \left ( \nabla \cdot \vec{m} \right )
\end{split} \end{equation}
which can be used to approximate solutions to $(\ref{CompNavStokes})$, in the sense that
\[ \| u(\cdot,t) - \tilde{G}(t)\ast u_0  \|_{L^p} \leq C t^{-\frac{d}{2}(1-\frac{1}{p})-\frac{1}{2}} \]
where $\tilde{G}$ is the Green's matrix of $(\ref{LinArtViscSys})$.  They go further by obtaining decay rates in $L^p$ for $1\leq p <2$, and show that the momentum field can be decomposed into an irrotational and incompressible piece, and that the solutions are asymptotically irrotational as measured in $L^p$ for $1\leq p <2$ and asymptotically incompressible for $p >2$.  Furthermore if one additionally has $(1+|x|)u_0 \in L^1$, then the solutions can be approximated by the explicit function
\[ \| u(\cdot,t) - \tilde{G}(t)U  \|_{L^p} \leq C t^{-\frac{d}{2}(1-\frac{1}{p})-\frac{1}{2}} \]
where $U = \int u_0dy$ is the total mass vector.  The Green's matrix is shown to consist of diffusing Gaussians which are convected by the fundamental solution of a wave equation, which they call diffusion waves.

Recently, Kagei and Okita \cite{kageiokita} extended the result of Hoff and Zumbrun by computing a second order approximation to the solutions of $(\ref{CompNavStokes})$ in dimension $d \geq 3$.  Among their findings, they prove that
\[ \big \| u(t) - G(t) \ast u_0 - \sum_{i=1}^d \partial_{x_i} G_1 (t,\cdot) \int \nolimits_0^{\infty} \int \nolimits_{\mathbb{R}^d} \mathcal{F}_i^0dyds \big \|_{L^p} \leq C\log(1+t)(1+t)^{-\frac{d}{2}(1-\frac{1}{p}) -\frac{3}{4}}  \]
for $p\geq 2$, where $G(t)$ is the Green's matrix for the linearization of $(\ref{CompNavStokes})$, $G_1(t)$ is a low frequency cutoff of $G(t)$, and the $\mathcal{F}_i^0$ are quantities which can be computed with knowledge of the solution $\rho(t), m(t)$, as well as knowledge of the pressure $P$ and its derivatives.  Their results also show that the solutions can be explicitly approximated by Gaussian functions 
\[ \Big \| u(t) - G_1(t) \int u_0dy + \sum_{|\alpha|=1} \partial_{x}^{\alpha} G_1 (t,\cdot) \int y^{\alpha}u_0dy - \sum_{i=1}^d \partial_{x_i} G_1 (t,\cdot) \int_0^{\infty} \int_{\mathbb{R}^d} \mathcal{F}_i^0dyds \Big \|_{L^p} \leq C\log(1+t)(1+t)^{-\frac{d}{2}(1-\frac{1}{p}) -\frac{3}{4}}  \]
if one includes the additional correction factor given by the $\mathcal{F}_i^0$ terms.

On the other hand, Gallay and Wayne (\cite{gallaywayne}, \cite{gallaywayne2}) study the asymptotic behavior of solutions of the incompressible equations in two and three dimensions.  They show that one can obtain an approximation of any desired order if one is willing to assume more spatial locality.  Specifically, if one chooses $\frac{3}{2} < \mu  \leq 2$ and $n \in \mathbb{Z}_{\geq 0}$ such that $n > 2\mu + \frac{1}{2}$ the for initial data $(1+|x|)^{n} u_0 \in L^2$ one can obtain the existence of solution $u(t)$ and approximations $u_{app,k}(t)$ such that
\[ \| u(t) - \sum_{k=1}^{n} u_{app,k}(t) \|_{L^p} \leq Ct^{-\frac{d}{2}(1-\frac{1}{p})+1-\mu } \] 
where the approximation terms $u_{app,k}(t)$ are also given in terms of diffusing Gaussians and their derivatives.  The key insight is that for each moment that the solution has one can compute a higher order, and hence more accurate, approximation.

We aim to use the tools developed in \cite{gallaywayne}, \cite{gallaywayne2} to extend the approximation of solutions to the compressible Navier-Stokes in \cite{hoffzumbrun1}, \cite{kageiokita} to a higher order.  To do so we begin with the modified compressible Navier-Stokes system
\begin{equation} \label{ModCompNavSto} \begin{split}
\partial_t \rho + \nabla \cdot \vec{m} & = \frac{1}{2} \big ( \epsilon + \eta \big ) \Delta \rho \\
\partial_t \vec{m} + \big [ \nabla \cdot \left ( \vec{m}\otimes \vec{m} \right ) \big ]^{T} + c^2 \nabla \rho & = \epsilon \Delta \vec{m} + \frac{1}{2} \big ( \eta - \epsilon \big ) \nabla \left ( \nabla \cdot \vec{m} \right )
\end{split} \end{equation}
obtained from (\ref{CompNavStokes}) by adding an artificial viscosity term and assuming that the density is bounded away from zero.  From the work of Hoff and Zumbrun, we know that the leading order long-time asymptotics of \eqref{ModCompNavSto} are the same as those of the compressible Navier-Stokes equations, but this model is somewhat simpler from a technical point of view, so we defer the consideration of \eqref{CompNavStokes} to future work.  For incompressible fluids, the work of Brandolese \cite{brandolese} showed that solutions of the Navier-Stokes equation which initially have finite moments may fail to have finite moments at later times.  This
problem does not occur if one works with the vorticity.  While it is not known if the momentum field of the compressible Navier-Stokes equation has this same property, we avoid its possible appearance by working with the curl and divergence of $\vec{m}$.  If one lets $a = \nabla \cdot \vec{m}$, $\vec{\omega} = \nabla \times \vec{m}$, and $u(t) = (\rho(t),a(t),\vec{\omega}(t))^T$ and computes the divergence and curl of $(\ref{ModCompNavSto})$, one arrives at the curl-divergence form of the modified compressible Navier-Stokes system:
\begin{equation} \label{RhoAOmegaSys} \partial_t u = \mathcal{L}u - \mathcal{Q}(u,u) \end{equation}
where we let $\nu = \frac{1}{2} (\epsilon + \eta)$, $I_3$ be the $3\times3$ identity matrix and where
\[ \mathcal{L} = \begin{pmatrix}
\nu \Delta & -1 & \\ -c^2 \Delta & \nu \Delta & \\  & & \epsilon \Delta I_3
\end{pmatrix} \hspace{.5 cm} \text{ , } \hspace{.5 cm} \mathcal{Q}(u,u) = \begin{pmatrix}  0 \\ \nabla \cdot \big [ \sum_{j=1}^3 \partial_{x_j} \big ( m_j \vec{m} \big ) \big ] \\ \nabla \times \big [ \sum_{j=1}^3 \partial_{x_j} \big ( m_j \vec{m} \big ) \big ] \end{pmatrix} \]
We take $(\ref{RhoAOmegaSys})$ as our starting point, and address the question of equivalence to the original system $(\ref{ModCompNavSto})$ in the course of our analysis.  Our main results can then be summarized in the following theorem:

\begin{Theorem}\label{t:Summary}
	Let $\vec{u}_0 = (\rho_0,a_0,\vec{\omega}_0)^T$ where $a_0,\vec{\omega}_0$ have zero total mass (ie $\int_{\mathbb{R}^3} a_0(x) dx = 0$), and suppose $(1+|x|)^n \vec{u}_0 \in W^{1,p} \times L^{p}\times \big (L^{p} \big )^3$ for some $0\leq n \leq 2$ and for all $1 \leq p \leq \frac{3}{2}$.  If $k \geq 1$ is fixed and if
	\[ E_n = \sup_{1\leq p \leq 3/2} \big ( \big \| (1+|x|)^n \rho_0 \big \|_{W^{1,p}} + \big \| (1+|x|)^n a_0 \big \|_{L^p} + \big \| (1+|x|)^n \vec{\omega}_0 \big \|_{(L^p)^3} \big ) \]
	is chosen sufficiently small, then there exists a unique mild solution $(\rho(t),a(t),\vec{\omega}(t))$ of $(\ref{RhoAOmegaSys})$ such that for a blowup rate $r_{\alpha,p}$ and decay rates $\ell_{n,p,\mu}, \tilde{\ell}_{n,p,\mu}$ defined below we have
	\[ \begin{split}\| \partial_{x}^{\alpha} \rho(t)\|_{L^p} & \leq  CE_n t^{-r_{\alpha,p}}(1+t)^{-\ell_{n,p,0}+\frac{1}{2}} \\ \| \partial_{x}^{\alpha} a(t)\|_{L^p} & \leq C E_n t^{-r_{\alpha,p}}(1+t)^{-\ell_{n,p,0}} \\ \| \partial_{x}^{\alpha} \vec{\omega}(t)\|_{\mathbb{L}^p} & \leq CE_n t^{-r_{\alpha,p}}(1+t)^{-\tilde{\ell}_{n,p,0}}
	\end{split} \]
	for $|\alpha| \leq k - 1$, where $C$ depends only on $n,k,\nu,\epsilon$.  Furthermore, for $n \geq 1$ there exist efficiently computable functions $(\rho_{app},a_{app},\vec{\omega}_{app})^T$ such that
	\[ \begin{split}\| \partial_{x}^{\alpha} \big ( \rho(t) - \rho_{app}(t) \big ) \|_{L^p} & \leq  CE_n t^{-r_{\alpha,p}}(1+t)^{-\ell_{n,p,0}+\frac{1}{2}-\frac{1}{2}} \\ \| \partial_{x}^{\alpha} \big ( a(t) - a_{app}(t) \big ) \|_{L^p} & \leq C E_n t^{-r_{\alpha,p}}(1+t)^{-\ell_{n,p,0}-\frac{1}{2}} \\ \| \partial_{x}^{\alpha}\big ( \vec{\omega}(t) - \vec{\omega}_{app}(t) \big )\|_{\mathbb{L}^p} & \leq CE_n t^{-r_{\alpha,p}}(1+t)^{-\tilde{\ell}_{n,p,0}-\frac{1}{2}}
	\end{split} \]
	and for $n=2$ one can take these functions to be explicitly computable with knowledge only of the moments of order $\lfloor n \rfloor$ of the initial data.
\end{Theorem}

For $n,\mu \in \mathbb{R}_{\geq 0}$, let $\lfloor n \rfloor_1 = \min(n,1)$ and $\lfloor \mu \rfloor_1 = \min (\mu,1)$, and we define the rates via
\begin{equation} \label{LinDecayRates1} \begin{split}
r_{\alpha,p} = \left \{ \begin{array}{cl}  \frac{|\alpha|}{2} & \text{ for } 1 \leq p \leq \frac{3}{2} \\ \frac{3}{2}(\frac{2}{3}-\frac{1}{p}) + \frac{|\alpha|}{2} & \text{ for } p \geq \frac{3}{2} \end{array} \right . \hspace{.2 cm} \text{ , } \hspace{.2 cm} \tilde{\ell}_{n,p,\mu} = \left \{ \begin{array}{cl}  \frac{3}{2}(1-\frac{1}{p}) + \frac{\lfloor n \rfloor_1 +\lfloor \mu \rfloor_1 }{2} -\mu & \text{ for } 1 \leq p \leq \frac{3}{2} \\ \frac{1}{2} + \frac{\lfloor n \rfloor_1 + \lfloor \mu \rfloor_1 }{2} -\mu & \text{ for } p \geq \frac{3}{2} \end{array} \right .
\end{split} \end{equation}
\begin{equation} \label{LinDecayRates2} \ell_{n,p,\mu} = \left \{ \begin{array}{cl}  \frac{5}{2}(1-\frac{1}{p}) -\frac{1}{2} + \frac{\lfloor n \rfloor_1 }{2} - \mu & \text{ for } 1 \leq p \leq \frac{3}{2} \\ (1-\frac{1}{p}) + \frac{\lfloor n \rfloor_1}{2} - \mu  & \text{ for } p \geq \frac{3}{2} \end{array} \right .  \end{equation}

In section 2, we prove a number of inequalities for later use in our existence and asymptotic analysis.  We also introduce an expansion for solutions of the heat equation which we call the Hermite expansion, and demonstrate how it works for related systems.  In section 3, we prove that \eqref{ModCompNavSto} has unique solutions, and that these solutions remain in the same weighted Lebesgue spaces as the initial data.  We obtain asymptotic rates for these solutions in weighted spaces, and find our solutions are asymptotically irrotational as measured in $L^p$ for $1\leq p < 2$ and asymptotically incompressible for $p >2$.  In section 4, we prove results about the accuracy of the linear approximation, and then show how this approximation can be improved if the initial data is appropriately localized.

\subsection{Mild formulation}

The nonlinear term in $(\ref{RhoAOmegaSys})$ still depends on $\vec{m}$, and hence we introduce the operators
\begin{equation} \label{PiAndBiotSavart} \begin{split}
\Pi a & = \nabla \big ( \Delta^{-1} a \big ) \\
B \vec{\omega} & = -\nabla \times \big ( \Delta^{-1} \vec{\omega} \big )
\end{split} \end{equation}
and 
\[ N(a,\vec{\omega}) = \sum_{j=1}^3 \partial_{x_j} \Big ( \big (\Pi a + B \vec{\omega} )_j \big ( \Pi a + B \vec{\omega} \big ) \Big ) \]
Note that the inverse Laplacian is well-defined only when we make a suitable choice of function spaces for $a$ and $\vec{\omega}$.  We will do so below in subsection $\ref{sect:PiAndBiotSavartOperators}$, and then obtain estimates for the action of $\Pi$ and $B$ over these spaces. 

We can now apply Duhamel's formula to obtain an integral formulation of $(\ref{RhoAOmegaSys})$:
\begin{equation} \label{IntegralForm} \begin{split}
\rho(t) & = \partial_t w (t) \ast K_\nu (t) \ast \rho_0  -  w (t) \ast K_\nu (t) \ast a_0 + \int_0^t w (t-s) \ast K_\nu (t-s) \ast \big [ \nabla \cdot N \big (a(s),\vec{\omega}(s) \big ) \big ] ds \\
a(t) & = -\partial_t^2 w (t) \ast K_\nu (t) \ast \rho_0 + \partial_t w (t) \ast K_\nu (t) \ast a_0 - \int_0^t \partial_t w (t-s) \ast K_\nu (t-s) \ast \big [ \nabla \cdot N \big (a(s),\vec{\omega}(s) \big ) \big ] ds \\
\vec{\omega}(t) & =  \mathbb{K}_\epsilon (t) \ast \vec{\omega}_0 - \int \nolimits_0^t \mathbb{K}_\epsilon (t-s) \ast \big [ \nabla \times N \big (a(s),\vec{\omega}(s) \big ) \big ]ds
 \end{split} \end{equation}
Here we use the fact that the Green's matrix $G$ for the linear part of the hyperbolic-parabolic system for $\rho,a$ above can be decomposed as the composition of the wave evolution with the heat evolution 
\begin{equation} \label{GreensMatrix}  G(t) \ast \begin{pmatrix}
\rho_0 \\ a_0 
\end{pmatrix} = G_W(t) \ast \big [ K_{\nu}(t)I_2 \ast \begin{pmatrix}
\rho_0 \\ a_0 
\end{pmatrix} \big ]  \end{equation}
in which 
\[ G_W(t) = \begin{pmatrix}
\partial_t w(t) & -w(t)  \\ 
-\partial_t^2 w(t) & \partial_t w(t) 
\end{pmatrix} \]
is the Green's matrix for the wave evolution, $K_\nu(t)= \frac{1}{(4\pi \nu t)^{3/2}} \exp \Big [ - \frac{|x|^2}{4\nu t} \Big ] $ is the scalar heat kernel and $I_2$ is the $2\times 2$ identity matrix.  The wave operator $w(t)$ is the Fourier multiplier defined by
\[ \hat{w}(\xi,t) = \frac{\sin (ct|\xi|)}{c|\xi|} \]
which together with its temporal derivatives determine the components of the wave evolution for various initial data.  For sufficiently smooth functions this can be expressed via Kirchhoff's formula, which in odd dimension $d\geq 3$ is as follows:
\begin{equation}
\label{Kirchhoff}
\begin{split}
(w\ast h)(x,t) & = \sum_{0\leq |\alpha| \leq \frac{d-3}{2}} b_{\alpha,0} (ct)^{|\alpha| + 1} \int \nolimits_{|z| = 1} D^{\alpha} h(x+ctz)z^{\alpha} dS(z)\\
(\partial_t w*h)(x,t) &= \sum_{0\leq |\alpha|\leq \frac{d-1}{2}} b_{\alpha,1} (ct)^{|\alpha|} \int \nolimits_{|z|=1} D^\alpha h(x+ ctz) z^\alpha dS(z)\\
(\partial_{t}^2 w*h)(x,t) &= \sum_{1\leq |\alpha|\leq \frac{d+1}{2}} b_{\alpha,2} (ct)^{|\alpha|-1}\int \nolimits_{|z|=1} D^\alpha h(x+ ctz) z^\alpha dS(z)
\end{split}
\end{equation}
with $S_z$ the surface element on the unit sphere, and some constants $b_{\alpha,i}$.  Finally, $\mathbb{K}_{\epsilon}(t)$ is the diagonal matrix having the heat kernel $K_{\epsilon}(t)$ for each entry on the diagonal.

We want to prove existence of mild solutions to $(\ref{RhoAOmegaSys})$ in some function space and determine the asymptotic behavior of these solutions.  We'll see that the natural setting for our analysis is found in the homogeneous, algebraically weighted Lebesgue spaces
\[ \mathring{L}^p(n) = \{ f\, : \,  \|f(x) \|_{\mathring{L}^p(n)} = \Big ( \int\nolimits_{\mathbb{R}^3} |x|^{np}|f(x)|^pdx \Big )^{1/p} < \infty \} \]
and their inhomogeneous counterparts
\[ L^p(n) = \{ f\, : \, \|f(x) \|_{L^p(n)} = \Big ( \int\nolimits_{\mathbb{R}^3} (1+|x|)^{np}|f(x)|^pdx \Big )^{1/p} < \infty  \} \]
We let $W^{k,p}(n)$ be the subspace of the Sobolev space $W^{k,p}$ consisting of algebraically weighted, weakly differentiable functions:
\[ W^{k,p}(n) = \{ f \in W^{k,p} : \|f\|_{W^{k,p}(n)}^p =  \sum_{|\alpha|\leq k} \|\partial_{x}^{\alpha} f\|_{L^p(n)}^p < \infty \} \]
We also introduce the vector-valued function space $\mathbb{L}^p = \big ( L^p \big )^d$ with norm
\[ \|\vec{\omega} \|_{\mathbb{L}^p} = \max_{i=1,2,3} \|\omega_i \|_{L^p} \]
as well as the function spaces $\mathring{\mathbb{L}}^p(n) = \big ( \mathring{L}^p(n) \big )^d$, $\mathbb{L}^p(n) = \big ( L^p(n) \big )^d$ and $\mathbb{W}^{k,p}(n) = (W^{k,p}(n))^d$ with analogous norms.  Furthermore let $\mathbb{L}^p_{\sigma}$ be the closure of the space of divergence free vector fields in the space $\mathbb{L}^p$, and let $\mathring{\mathbb{L}}^p_{\sigma}(n)$, $\mathbb{L}^p_{\sigma}(n)$ and $\mathbb{W}^{k,p}_{\sigma}(n)$ be the closures in the analogous spaces.  Finally, we will make use of Schwartz class functions as tools in our analysis, and hence we will write $\mathcal{S}$ for the space of Schwartz class functions and $\mathbb{S}_{\sigma}$ for the space of Schwartz class divergence free vector fields.

\section{Preliminary analysis}

\subsection{The $\Pi$ and B operators}

\label{sect:PiAndBiotSavartOperators}

The results here apply to the $d=3$ case.  We first define the operators $\Pi$ and $B$ for $(a,\vec{\omega}) \in \mathcal{S}\times \mathbb{S}_{\sigma}$ via $(\ref{PiAndBiotSavart})$.  Note that the inverse Laplacian is well defined on the space of Schwartz class functions, and for such functions we have
\[ \Pi a = -\frac{1}{4\pi} \int \nolimits_{\mathbb{R}^3} \frac{(x-y)}{|x-y|^3} a(y)dy \hspace{1 cm} \text{ and } \hspace{1 cm} B\vec{\omega} = \frac{1}{4\pi} \int \nolimits_{\mathbb{R}^3} \frac{(x-y) \times \vec{\omega}(y) }{|x-y|^3} dy \]

In the following proposition we obtain estimates on the action of $\Pi$ and $B$, which then allow us to extend these operators to be defined on all of $L^p(n)\times \mathbb{L}_{\sigma}^p(n)$, for suitable choices of $p$ and $n$. 

\begin{Proposition} \label{p:BiotSavart}
	Let $a \in \mathcal{S}$ and $\vec{\omega} \in \mathbb{S}_{\sigma}$
	
	\textbf{(a)} Suppose that $1 < p_1 < \infty$.  Then there exists a constant $C_1$ depending only on $p_1$ such that
	\begin{equation} \label{CalderonZygmundEst} \begin{split}  \|\partial_{x_i} \Pi a \|_{L^{p_1}} \leq C_1 \|a\|_{L^{p_1}} \hspace{.5 cm} \text{ , } \hspace{.5 cm} \|\partial_{x_i} B \vec{\omega} \|_{\mathbb{L}^{p_1}} \leq C_1 \| \vec{\omega} \|_{\mathbb{L}^{p_1}} \end{split} \end{equation}
	
	\textbf{(b)} Suppose that $n \in [0,2)$ and $1 < p_3 < p_2 < \infty$ are such that
	\begin{equation} \label{HardyLittlewoodSobolevIndexCondition} \frac{1}{p_2} = \frac{1}{p_3} - \frac{1}{3} \end{equation}
	and $p_3$ satisfies the constraint
	\[ \frac{1-n}{3} < \frac{1}{p_3} < \frac{3-n}{3} \]
	Then there exists a constant $C_2$ depending only on $n,p_3$ such that
	\begin{equation} \label{BiotSavartEst} \begin{split}  \|\Pi a \|_{L^{p_2}(n)} \leq C_2 \|a\|_{L^{p_3}(n)} \hspace{.5 cm} \text{ , } \hspace{.5 cm} \|B \vec{\omega} \|_{\mathbb{L}^{p_2}(n)} \leq C_2 \|\vec{\omega} \|_{\mathbb{L}^{p_3}(n)} \end{split} \end{equation}
	
	\textbf{(c)} Suppose $n \in [1,3)$, $1 < p_3 < p_2 < \infty$ solve (\ref{HardyLittlewoodSobolevIndexCondition}) and $p_3$ satisfies the new constraint
	\[ \frac{3-n}{3} < \frac{1}{p_3} < \frac{4-n}{3} \]
	If, in addition, $a$ and $\vec{\omega}$ are such that
	\begin{equation} \label{ZeroMass} \int \nolimits_{\mathbb{R}^3} a(x) dx = 0 \hspace{.5 cm} \text{ , } \hspace{.5 cm} \int \nolimits_{\mathbb{R}^3}\vec{\omega}(x)dx = 0 \end{equation}
	then there exists a (possibly different) constant $C_2$ depending only on $n,p_3$ such that $(\ref{BiotSavartEst})$ holds.
	
\end{Proposition}

The proof of these estimates follows closely the strategy used to the study the $B$ operator in Proposition B.1 of \cite{gallaywayne2}, but we extend
the results to general values of $p$ and $n$, rather than focusing on the $L^2$ based spaces in that
reference, as well as studying the operator $\Pi$.   We defer the proof to Appendix \ref{app:BiotSavart}.  The following Corollary is immediate from the definition of the $\Pi$,$B$ operators for $a,\vec{\omega} \in L^{p_3}(n) \times \mathbb{L}^{p_3}_{\sigma}(n)$:

\begin{Corollary}\label{c:PiAndBiotSavartOperators}
	\textbf{(a)} Suppose $p_1, C_1$ are as in Prop. $\ref{p:BiotSavart}$ part \textbf{(a)}.  Then for $a \in L^{p_1}$, $\vec{\omega} \in \mathbb{L}_{\sigma}^{p_1}$ $(\ref{CalderonZygmundEst})$ holds.
	
	\textbf{(b)} Suppose that $n, p_2, p_3, C_2$ are as in Prop. $\ref{p:BiotSavart}$ part \textbf{(b)}.  Then for $a \in L^{p_3}(n)$ and $\vec{\omega} \in \mathbb{L}_{\sigma}^{p_3}(n)$ $(\ref{BiotSavartEst})$ holds.
	
	\textbf{(c)} Suppose that $n, p_2, p_3,C_2$ are as in Prop. $\ref{p:BiotSavart}$ part \textbf{(c)}.  If $a \in L^{p_3}(n)$ and $\vec{\omega} \in \mathbb{L}_{\sigma}^{p_3}(n)$ satisfy (\ref{ZeroMass}), then $(\ref{BiotSavartEst})$ holds.
	
\end{Corollary}

\subsection{Heat evolution estimate}

In any dimension $d\geq1$, the heat evolution tends to dissipate the $L^p$ norms of a function.  We have
\[ \| \partial_x^{\alpha} K_{\nu}(t) \ast f \|_{L^p} \leq C(\nu t)^{-\frac{|\alpha|}{2}-\frac{d}{2}(\frac{1}{q}-\frac{1}{p})} \| f\|_{L^q} \]
using Young's inequality for $1 \leq q \leq p \leq \infty$ and $f \in L^q$.  In weighted spaces, one can obtained faster decay under certain conditions described in the following proposition, which is an extension of Proposition A.3 found in \cite{gallaywayne2}.  We defer the proof to Appendix \ref{app:HeatEst}.

\begin{Proposition}\label{p:HeatEst}
	For dimension $d\geq 1$, let $1 \leq q \leq p \leq \infty$ be Lebesgue indices, let $n, \mu \in \mathbb{R}_{\geq 0}$ be weight indices such that $n \geq \mu$ and that $\exists$ $\tilde{n} \in \mathbb{Z}_{\geq 0}$ such that $d(1-\frac{1}{q}) + \tilde{n} < n < d(1-\frac{1}{q}) + \tilde{n} +1$, and let $f \in L^{q}(n)$ be such that its moments up to order $\tilde{n}$ are zero, ie for all multi-indices $\beta \in \mathbb{N}^d$, $|\beta| \leq \tilde{n}$ we have
	\[ \int \nolimits_{\mathbb{R}^3} x^{\beta} f(x) dx = 0 \]
	Then there exists a $C > 0$ depending only on $d,p,q,n,\mu,\alpha$ such that
	\begin{equation} \label{HeatEst} \| \partial_x^{\alpha} K_{\nu}(t) \ast f \|_{\mathring{L}^p(\mu)} \leq C (\nu t)^{-\frac{|\alpha|}{2}-\frac{d}{2}(\frac{1}{q}-\frac{1}{p})} (1+\nu t)^{- \frac{n-\mu}{2}} \|f\|_{L^q(n)}  \end{equation}
\end{Proposition}

\subsection{Heat-wave evolution estimate}

We obtain the following bounds on the heat-wave operators of the linear evolution of the $\rho,a$ system in homogeneous weighted spaces in a general odd dimension $d \geq 3$.

\begin{Proposition}\label{p:HypParaEst} For general odd dimension $d\geq 3$, Lebesgue index $q\geq1$ and weight $n\geq 0$ there exists a $C>0$ depending only on $d,c,\nu,n$ such that the following estimates hold:
	\begin{equation}
	\label{HypParaEst}
	\begin{split}
	\| w(t) \ast K_{\nu}(t) \|_{\mathring{L}^q(n)} & \leq C t^{1+\frac{n}{2} - \frac{d}{2}(1-\frac{1}{q})}(1+t)^{\frac{n}{2}+\frac{d-3}{4} -\frac{d-1}{2}(1-\frac{1}{q})}\\
	\| \partial_t w (t)\ast K_{\nu}(t) \|_{\mathring{L}^q(n)} & \leq C  t^{\frac{n}{2} - \frac{d}{2}(1-\frac{1}{q})}(1+t)^{\frac{n}{2}+\frac{d-1}{4} -\frac{d-1}{2}(1-\frac{1}{q})} \\
	\| \partial_{t}^2 w(t) \ast K_{\nu}(t) \|_{\mathring{L}^q(n)} & \leq C t^{\frac{n}{2} -\frac{1}{2} - \frac{d}{2}(1-\frac{1}{q})}(1+t)^{\frac{n}{2}+\frac{d-1}{4} -\frac{d-1}{2}(1-\frac{1}{q})} 
	\end{split}
	\end{equation}
	
\end{Proposition}

We defer the proof to Appendix \ref{app:HeatWaveEst}.  Note that the term $\partial_t^2 w(t)\ast K_{\nu}(t)$ blows up as $t \to 0$ as a result of the fact that $K_{\nu}(t)$ tends to a delta function, and hence the $L^p$ norms of derivatives of $K_{\nu}(t)$ in the formula (\ref{Kirchhoff}) become arbitrarily large.  However, when the heat-wave operator  $\partial_t^2 w(t)\ast K_{\nu}(t)$ acts on a function with a little bit of smoothness we can obtain the following improved estimate with milder blow up, the proof of which we defer to Appendix \ref{app:HeatWaveEst2}:

\begin{Proposition}
	\label{p:HypParaEst2}
	
	Suppose $\rho_0 \in W^{1,q}(n)$ for some $q \geq 1$.  There exists a $C>0$ such that for $p\geq q$ and $\mu \leq n$ we have
	\[ \|\partial_t^2 w(t) \ast K_{\nu}(t)\ast \rho_0 \|_{\mathring{L}^p(\mu)} \leq C t^{-\frac{d}{2}(\frac{1}{q}-\frac{1}{p})}(1+t)^{\mu-\frac{1}{2}+\frac{d-1}{4}-\frac{d-1}{2}(\frac{1}{q}-\frac{1}{p})} \|\rho_0 \|_{W^{1,q}(n)} \]
\end{Proposition}

\subsection{Hermite expansion}

We aim to study the asymptotic behavior of solutions to $(\ref{RhoAOmegaSys})$ by computing an expansion of the solution using Hermite functions.  This is the point where we begin to diverge strongly 
from the approach of \cite{hoffzumbrun1} or \cite{kageiokita}.  We illustrate this process first for the heat equation.  To do so, we define
\[ \phi_0(x) = (4\pi)^{-\frac{d}{2}} \exp \big [ -\frac{|x|^2}{4} \big ] \]
and let $H_{\alpha}$ be the $\alpha$th Hermite polynomial given by
\[ H_{\alpha}(x) = \frac{2^{|\alpha|}}{\alpha !} e^{\frac{|x|^2}{4}} \partial_{x}^\alpha (e^{-\frac{|x|^2}{4}}) \]
Note that these satisfy the orthonormality property:
\begin{equation}
	\label{HermitePolyOrthoProp}
	\langle H_{\alpha}(\cdot), \partial_{x}^{\beta} \phi_0(\cdot) \rangle = \delta_{\alpha\beta}
\end{equation}

\begin{Proposition}
	\label{p:HermiteExp}
	For arbitrary dimension $d\geq 1$, suppose that $u_0 \in L^1(n)$ for $n \geq 0$.  If $u(t) = K_{\nu}(t) \ast u_0$ is the solution of the heat equation in $C^0 \big [ [0,\infty), L^1(n) \big ]$, then we can write
	\[ u(x,t) = \sum_{|\alpha|\leq \floor{n}} \big \langle H_{\alpha}, u_0 \big \rangle \partial_{x}^{\alpha} K_{\nu}(t)\ast  \phi_0(x) + R(x,t) \]
	where for any $\mu\leq n$
	\[ \| R(\cdot,t) \|_{\mathring{L}^{p}(\mu)} \leq C\|u_0\|_{L^1(n)}(\nu t)^{-\frac{d}{2}(1-\frac{1}{p})-\frac{n-\mu}{2}} \]
\end{Proposition}  

\begin{proof}
	If we write
	\[ u(x,t) = \sum_{|\alpha|\leq \floor{n}} \big \langle H_{\alpha}, u_0 \big \rangle \partial_{x}^{\alpha} K_{\nu}(t)\ast  \phi_0(x) + R(x,t) \]
	then we note that the remainder term $R(x,t)$ is itself a solution of the heat equation.  Furthermore, we note that at time $t=0$ we have
	\[ \big \langle H_{\beta},R(\cdot,0) \big \rangle  = 0 \]
	for all $|\beta|\leq \floor{n}$.  Therefore $R_j$ satisfies the moment zero condition required in Proposition $\ref{p:HeatEst}$, which then gives us our result.
\end{proof}

The Hermite expansion illustrates a few of the features of the heat evolution.  We note that orders of this expansion decay sequentially faster, and the remainder at least matches the fastest decay rate.  The Hermite functions are self similar under the heat evolution, in the sense that the heat evolution acts on these functions by dilation and scaling.  See \cite{gallaywayne} for details.  Importantly, the Hermite expansion illustrates how the heat evolution dissipates the moments of a function.  The $\alpha$th moment evolves according to the $\alpha$th term in the Hermite expansion.  For instance, the zeroth order Hermite function gives an explicit example of an initial condition for which the heat evolution preserves the $L^1$ norm, yet has any degree of algebraic decay one could ask for, and hence the estimate in $(\ref{HeatEst})$ is sharp with respect to the zero mass condition.  However, the $L^\infty$ norm decays, so here the heat evolution is spreading mass around, but it conserves the total signed mass.  The first order Hermite function provides an example where the total signed mass is zero, and we see that its $L^1$ norm does decay.  The Hermite expansion can be used to show that this holds in general, and similar statements can be made about higher order moments.

\subsubsection{Hermite expansion for the hyperbolic-parabolic system}

We need a Hermite expansion for the hyperbolic-parabolic system
\begin{equation} \label{HypParaSys} \begin{split}
\partial_t \rho_L & = \nu \Delta \rho_L - a_L \\
\partial_t a_L & = -c^2 \Delta \rho_L + \nu \Delta a_L
\end{split} \end{equation}
As in (\ref{RhoAOmegaSys}) we can write the solution of the linear equation in terms of the heat-wave operators via
\begin{equation} \label{RhoALinearEvolution} \begin{split}
\rho_L(t) & = \partial_t w (t) \ast K_\nu (t) \ast \rho_0 -  w (t) \ast K_\nu (t) \ast a_0 \\
a_L(t) & = -\partial_t^2 w (t) \ast K_\nu (t) \ast \rho_0 + \partial_t w (t) \ast K_\nu (t) \ast a_0
\end{split} \end{equation}
Since the heat and wave operators commute, we can apply them sequentially, and since $K_{\nu}(t)\ast \rho_0$ and $K_{\nu}(t) \ast a_0$ are solutions of the heat equation, we can use the scalar Hermite expansion.  We define
\begin{equation}
\label{RhoAHermiteProfiles}
\begin{split}
\begin{pmatrix}
\rho_1(t) \\ a_1(t) 
\end{pmatrix} = \begin{pmatrix} \partial_t w (t) \ast K_\nu (t) \ast \phi_0 \\ - \partial_t^2 w (t) \ast K_\nu (t) \ast \phi_0
\end{pmatrix} \hspace{0.5 cm} \text{ and } \hspace{0.5 cm}  \begin{pmatrix}
\rho_{2}(t) \\ a_{2}(t)
\end{pmatrix} = \begin{pmatrix} -w (t) \ast K_\nu (t) \ast \phi_0 \\ \partial_t w (t) \ast K_\nu (t) \ast \phi_0
\end{pmatrix}
\end{split}
\end{equation}
where $\hat{e}_j$ are the standard unit two vectors.  We determine these asymptotic profiles explicitly in Appendix \ref{app:AsympProfs} below.  We then have the following analogue of the Hermite expansion, where for convenience we assume that $\rho$ has at least one weak derivative:

\begin{Proposition}
	\label{p:HypParaHermiteExp}
	
	For general odd dimension $d\geq 3$, suppose that $\rho_0 \in W^{1,1}(n)$, $a_0 \in L^1(n)$ for $n \geq 0$.  If $(\rho_L(t), a_L(t))^T$ is the solution of $(\ref{HypParaSys})$ in $C^0 \big [ [0,\infty) , L^1(n) \times L^1(n) \big ] $, then we can write
	\[ \begin{split}
	\begin{pmatrix}
	\rho_L(x,t) \\ a_L(x,t)
	\end{pmatrix} = \sum_{\substack{i\leq 2 \\ |\alpha|\leq \lfloor n \rfloor }} \big \langle H_{\alpha} \vec{e}_i, \begin{pmatrix}
	\rho_0 \\ a_0 
	\end{pmatrix} \big \rangle \text{ } \partial_{x}^{\alpha} \begin{pmatrix}
	\rho_i(x,t) \\ a_i(x,t)
	\end{pmatrix}  + \begin{pmatrix}
	\rho_{LR}(x,t) \\ a_{LR}(x,t)
	\end{pmatrix}
	\end{split} \]
	where for any $\mu\leq n$
	\begin{equation}
	\label{HypParaRemEst}
	\begin{split}
	\| \rho_{LR}(\cdot,t) \|_{\mathring{L}^{p}(\mu)} \leq C \big ( \|\rho_0\|_{L^1(n)} + \|a_0\|_{L^1(n)} \big ) t^{-\frac{d}{2}(1-\frac{1}{p})} (1+t)^{\frac{d+1}{4}-\frac{d-1}{2}(1-\frac{1}{p})-\frac{n}{2}+\mu} \\ \| a_{LR}(\cdot,t) \|_{\mathring{L}^{p}(\mu)} \leq C\big ( \|\rho_0\|_{L^1(n)} + \|a_0\|_{L^1(n)} \big )t^{-\frac{d}{2}(1-\frac{1}{p})} (1+t)^{\frac{d-1}{4}-\frac{d-1}{2}(1-\frac{1}{p})-\frac{n}{2}+\mu}
	\end{split} \end{equation}
\end{Proposition}   

\noindent The proof makes use of Props \ref{p:HypParaEst}, \ref{p:HypParaEst2} and is similar to that of Prop \ref{p:HermiteExp}, and we leave it to the reader.

\subsubsection{Hermite expansion for divergence free vector fields}

We will assume the dimension $d = 3$ for the remainder of the paper.  When considering the asymptotics of the vorticity equation, we will need a Hermite expansion for divergence free vector fields.  If we write 
\begin{equation}
\label{OmegaLinearEvolution} \vec{\omega}_L(t) = \mathbb{K}_{\epsilon}(t)\ast \vec{\omega}_0
\end{equation}
and naively expand each component of $\vec{\omega}(t)$ using the scalar Hermite expansion, the terms we obtain are not, in general, divergence free.  For the purposes of this paper, we will only consider Hermite expansions out to moments of order 2, hence we define these asymptotic profiles in the following table.  Higher order Hermite expansion can be defined and we let $\vec{p}_{\tilde{\alpha},j} = \vec{f}_{\tilde{\alpha},j} = 0$ for all $|\tilde{\alpha}| \leq 3$ not listed below.  We determine these profiles explicitly in Appendix \ref{app:AsympProfs} below.

\begin{table}[H]
	\setlength\extrarowheight{5pt}
	\begin{center}
		\begin{tabular}{|c|c|c|c|}
			\rowcolor{gray!50}
			\hline 
			$\tilde{\alpha}$ & j &$\vec{p}_{\tilde{\alpha},j}$ & $\vec{f}_{\tilde{\alpha},j}$ \\ \hline \hline
			(1,1,0) & 1 & $(-\frac{1}{2}x_2,\frac{1}{2}x_1,0)^T$ &  $\nabla \times (\phi_0 \vec{e}_3)$ \\
			
			(1,0,1) & 1 & $(\frac{1}{2}x_3,0,-\frac{1}{2}x_1)^T$ &  $\nabla \times (\phi_0 \vec{e}_2)$ \\
			
			(0,1,1) & 1 & $(0,-\frac{1}{2}x_3,\frac{1}{2}x_2)^T$ &  $\nabla \times (\phi_0 \vec{e}_1)$ \\
			\hline \hline 
			
			(2,1,0) & 1 & $(\frac{1}{2}x_1x_2,-\frac{1}{4}x_1^2,0)^T$ & $\nabla \times ( \partial_{x_1}\phi_0 \vec{e}_3)$\\
			
			(1,2,0) & 1 & $(\frac{1}{4}x_2^2,-\frac{1}{2}x_1x_2,0)^T$ & $\nabla \times ( \partial_{x_2}\phi_0 \vec{e}_3)$ \\
			
			(2,0,1) & 1 & $(-\frac{1}{2}x_1x_3,0,\frac{1}{4}x_1^2)^T$ & $\nabla \times ( \partial_{x_1}\phi_0 \vec{e}_2)$ \\
			
			(1,0,2) & 1 & $(-\frac{1}{4}x_3^2,0,\frac{1}{2}x_1x_3)^T$ & $\nabla \times ( \partial_{x_3}\phi_0 \vec{e}_2)$ \\
			
			(0,2,1) & 1 & $(0,\frac{1}{2}x_2x_3,-\frac{1}{4}x_2^2)^T$ & $\nabla \times ( \partial_{x_2}\phi_0 \vec{e}_1)$\\
			
			(0,1,2) & 1 & $(0,\frac{1}{4}x_3^2,-\frac{1}{2}x_2x_3)^T$ &  $\nabla \times ( \partial_{x_3}\phi_0 \vec{e}_1)$ \\
			\hline 
			(1,1,1) & 1 & $(x_2x_3,0,0)^T$ &  $\nabla \times ( \partial_{x_3}\phi_0 \vec{e}_3)$ \\
			
			(1,1,1) & 2 & $(0,0,-x_1x_2)^T$ & $\nabla \times ( \partial_{x_1}\phi_0 \vec{e}_1)$ \\
			\hline 
		\end{tabular}
	\end{center}
	
	\setlength\extrarowheight{0pt}
	
	\caption{Asymptotic profiles for the divergence-free vector field Hermite expansion.  See also \cite{gallaywayne2}. }
	\label{tab:DivFreeHermiteExp}
\end{table}
All of the profiles $\vec{f}_{\tilde{\alpha},j}$ are clearly divergence-free, and straightforward computations show that for $\vec{p}_{\tilde{\alpha},j}, \vec{f}_{\tilde{\alpha},j}$ defined above we have the orthonormality condition
\[ \big \langle \vec{p}_{\tilde{\alpha},j}, \vec{f}_{\tilde{\beta},k} \big \rangle = \delta_{jk} \delta_{\tilde{\alpha}\tilde{\beta}} \]
We then have the an analogue of the Hermite expansion, and we again leave the proof to the reader:


\begin{Proposition}
	\label{p:VecHermiteExp}
	
	Suppose that $\vec{\omega}_0 \in \mathbb{L}_{\sigma}^1(n)$ for $0\leq n \leq 2$.  If $\vec{\omega}_L(t)$ is the solution of the heat equation in $C^0 \big [ [0,\infty), \mathbb{L}^1_{\sigma}(n) \big ]$ given by (\ref{OmegaLinearEvolution}), then we can write
	\[  \vec{\omega}_L(x,t) = \sum_{\substack{j\leq 2 \\ |\tilde{\alpha}|\leq \lfloor n \rfloor +1}} \big \langle \vec{p}_{\tilde{\alpha},j},\vec{\omega}_0  \big \rangle \mathbb{K}_{\epsilon}(t) \ast \vec{f}_{\tilde{\alpha},j}(x) + \vec{\omega}_{LR}(x,t) \]
	where for any $\mu\leq n$
	\begin{equation} \label{HeatRemEst} \begin{split}
	\| \vec{\omega}_{LR}(\cdot,t) \|_{\mathring{\mathbb{L}}^{p}(\mu)} \leq C \|\vec{\omega}_0\|_{\mathbb{L}^1(n)} (\nu t)^{-\frac{3}{2}(1-\frac{1}{p})-\frac{n-\mu}{2}}
	\end{split} \end{equation}
\end{Proposition}

\section{Existence and uniqueness of solutions for the $(\rho,a, \vec{\omega})$ system}

Note that from the form of $(\ref{IntegralForm})$, if we can prove the existence of $a$ and $\vec{\omega}$, we can get the solution for $\rho$ by integration.  Hence we need to choose a function space for $(a,\vec{\omega})$.  In the Hermite expansions above, we saw that we could obtain higher order approximations by increasing the spatial localization of the initial conditions.  Hence for a given $n \in \mathbb{R}_{\geq 0}$ we might choose $(\rho_0,a_0,\vec{\omega}_0) \in L^1(n)\times L^1(n)\times \mathbb{L}^1_{\sigma}(n)$ as a sufficiently general space to start with, and expect to obtain solutions with $\lfloor n \rfloor$ orders of asymptotic profiles.  Note however that we expect that $a$ and $\vec{\omega}$ come from a velocity vector field via $a = \nabla \cdot \vec{m}$ and $\vec{\omega} = \nabla \times \vec{m}$, hence we can assume they have zero total mass as in (\ref{ZeroMass}).  Since $\vec{m}$ is assumed to have at least one derivative, we assume that $\rho$ has at least one as well, hence we assume $(\rho_0,a_0,\vec{\omega}_0) \in W^{1,1}(n)\times L^1(n)\times \mathbb{L}^1_{\sigma}(n)$

It will be desirable that the moments be continuous functions of time.  To obtain this we will see that we need a slightly stronger assumption: we require that $(\rho_0,a_0,\vec{\omega}_0)$ belong to $W^{1,\tilde{p}}(n)\times L^{\tilde{p}}(n)\times \mathbb{L}^{\tilde{p}}_{\sigma}(n)$ for all $1 \leq \tilde{p} \leq 3/2$.  We therefore define the function space
\begin{equation}
\label{IntersectionSpace1}
Z_{n}^{0} = \bigcap \limits_{1\leq p < \frac{3}{2}} C^{0} \big [ [0,\infty) , L^p(n) \times \mathbb{L}^p_{\sigma}(n) \big ] 
\end{equation}
Due to the smoothing properties of the heat evolution the solutions have more regularity for $t > 0$, so if we fix a degree of smoothness $k\geq 1$ we define
\begin{equation}
\label{IntersectionSpace2}
Z_{n,k}^{+} = \bigcap \limits_{1\leq p \leq \infty} C^{0} \big [ (0,\infty) , W^{k,p}(n) \times \mathbb{W}^{k,p}_{\sigma}(n) \big ] 
\end{equation}

Our existence analysis begins by studying the linear part of the evolution in \ref{RhoAOmegaSys}.  To this end we let $(\rho_L(t), a_L(t), \vec{\omega}_L(t))^T$ be defined by (\ref{RhoALinearEvolution}) and (\ref{OmegaLinearEvolution}) for $t > 0$ and $(\rho_L(t), a_L(t), \vec{\omega}_L(t))^T = (\rho_0,a_0,\vec{\omega}_0)^T$ for $t = 0$.  
In Appendix $\ref{app:LinEvol}$, we determine the smoothness properties and decay rates of these functions.  Based on our findings we look for solutions of $(\ref{IntegralForm})$ in the function space
\begin{equation} \label{FunctionSpace} \begin{split}
X_{n,k} = \Big \{  (a,\vec{\omega}) \in  Z_{n}^{0} \cap Z_{n,k}^{+} : \int \nolimits_{\mathbb{R}^3} a(x,t)dx = 0 \text{ and } \int \nolimits_{\mathbb{R}^3} \vec{\omega}(x,t)dx = 0  \Big \}
\end{split}  \end{equation}
with norm
\[ \begin{split}
\|(a,\vec{\omega})\|_{X_{n,k}} = \sup_{|\alpha|\leq k} \sup_{1\leq p \leq \infty} \sup_{0\leq \mu \leq n} \sup_{0 < t < \infty} \Big [ t^{r_{\alpha,p}}(1+t)^{\ell_{n,p,\mu}+\hat{\ell}_{k,p,\alpha}} \| \partial_x^{\alpha} a(t)\|_{\mathring{L}^p(\mu)} + t^{r_{\alpha,p}}(1+t)^{\tilde{\ell}_{n,p,\mu}} \|\partial_x^{\alpha} \vec{\omega}(t)\|_{\mathring{\mathbb{L}}^p(\mu)} \Big ] \end{split} \]
where $r_{\alpha,p}$, $\tilde{\ell}_{n,p,\mu}$ and $\ell_{n,p,\mu}$ are as in (\ref{LinDecayRates1}), (\ref{LinDecayRates2}), and $\hat{\ell}_{k,p,\alpha}$ is defined by
\begin{equation}
\label{NonlinearDecayRate}
\hat{\ell}_{k,p,\alpha} = \left \{ \begin{array}{cl} 0 & \text{ for } |\alpha| < k \text{ and for } |\alpha| = k \text{ , } 1 \leq p \leq 2  \\ -\frac{2}{3}(1-\frac{1}{p}) + \frac{1}{3}  & \text{ for } |\alpha| = k \text{ , } p \geq 2 \end{array} \right .
\end{equation}
The factor $\hat{\ell}_{k,p,\alpha}$ accounts for a slightly slower admissible decay rate for the highest order derivative in $L^p$, $p > 2$ as compared to the linear evolution.  Note that $X_{n,k}$ is a Banach space with this norm.  We will also need to define
\[ L_{n,\tilde{n}}(t) = \log(1+t) \text{ when } n=\tilde{n} \text{ and } L_{n,\tilde{n}}(t) = 1 \text{ otherwise } \]

\begin{Theorem}\label{t:Existence}
	Fix $n \in [ 0,2 ]$, $k \geq 1$ and let $(\rho_0,a_0,\vec{\omega}_0)$ belong to $W^{1,p}(n)\times L^{p}(n)\times \mathbb{L}_{\sigma}^{p}(n)$ for all $1 \leq p \leq \frac{3}{2}$, and suppose that $a_0$ and $\vec{\omega}_0$ have zero total mass.  If 
	\begin{equation} \label{InitialNorm} E_n = \sup_{1\leq p \leq 3/2} \big ( \|\rho_0 \|_{W^{1,p}(n)} + \|a_0 \|_{L^p(n)} + \|\vec{\omega}_0 \|_{\mathbb{L}^p(n)} \big ) \end{equation}
	is chosen sufficiently small, then there exists a unique solution $(a(t),\vec{\omega}(t))$ of $(\ref{IntegralForm})$ belonging to $X_{n,k}$ such that $(a(0),\vec{\omega}(0)) = (a_0,\vec{\omega}_0)$.
\end{Theorem}

\begin{proof}
Having chosen an initial condition satisfying the above, define the map $F_{(\rho_0,a_0,\vec{\omega}_0)}$ on $X_{n,k}$ sending $(a(s),\vec{\omega}(s))^T$ to a new function of space and time by letting $F_{(\rho_0,a_0,\vec{\omega}_0)}[(a,\vec{\omega})](0) = (a_0,\vec{\omega}_0)^T$ and 
\[ F_{(\rho_0,a_0,\vec{\omega}_0)} \big [ a , \vec{\omega} \big ](t) = \begin{pmatrix}
-\partial_t^2 w \ast K_\nu \ast \rho_0 + \partial_t w \ast K_\nu \ast a_0 - \int_0^t \big [\partial_t w \ast K_\nu \big ] (t-s) \ast \Big [ \nabla \cdot N \big (a(s),\vec{\omega}(s) \big ) \Big ]ds \\
\mathbb{K}_\epsilon \ast \vec{\omega}_0 - \int_0^t \mathbb{K}_\epsilon (t-s) \ast \Big [ \nabla \times N \big (a(s),\vec{\omega}(s) \big ) \Big ]ds
\end{pmatrix} \]
for $t > 0$.  For convenience, we'll drop the subscript.  We claim that $F$ maps $X_{n,k}$ into itself and has Lipschitz constant equal to $1/2$ on a ball of radius $R$ centered at the origin, which we prove below.  Given these two claims, we can conclude our proof as follows.  If $(a_L,\vec{\omega}_L)$ are as above, we note that each of the bounds determined in Appendix \ref{app:LinEvol} depend on the magnitude of the initial condition, hence 
\[ \big \| (a_L,\vec{\omega}_L ) \big \|_{X_{n,k}} \leq CE_n  \]

Therefore if we choose the initial condition sufficiently small, (ie $E_n \leq \frac{R}{2C}$) we then have
\[ \begin{split} \| F \big (
a , \vec{\omega} \big ) - (a_L,\vec{\omega}_L) \|_{X_{n,k}} & = \| F \big (
a , \vec{\omega}  \big ) - F \big ( 
0 , 0 \big ) \|_{X_{n,k}} \\ & \leq \frac{1}{2} \| \big (
a , \vec{\omega}  \big ) - \big (a_L,\vec{\omega}_L \big ) \|_{X_{n,k}} +\frac{1}{2} \| \big (a_L,\vec{\omega}_L \big ) \|_{X_{n,k}} \leq \frac{R}{2}
\end{split} \]
for $(a,\vec{\omega}) \in B\big ( (a_L,\vec{\omega}_L),\frac{R}{2} \big )$, the closed ball of radius $\frac{R}{2}$ centered at $(a_L,\vec{\omega}_L)^T$.  Therefore $F$ maps $B\big ( (a_L,\vec{\omega}_L)^T,\frac{R}{2} \big )$ into itself, and since $F$ is a contraction here, the unique solution of $(\ref{IntegralForm})$ is given by the fixed point of $F$.

\paragraph{Claim One: $F:X_{n,k} \mapsto X_{n,k}$ }

We begin by proving that for $(a,\vec{\omega}) \in X_{n,k}$ the $X_{n,k}$ norm of $F(a,\vec{\omega})$ is finite and that $F(a,\vec{\omega})\in Z_{n}^0\cap Z_{n,k}^+$.  We note again that the decay rates and smoothness requirements to belong to $X_{n,k}$ were found to be more than satisfied by those of the linear terms in Appendix \ref{app:LinEvol}, so we need only analyze the evolution of the Duhamel terms.  Furthermore we note that is sufficient to bound the $\mathring{L}^p(\mu)$ norms for $\mu = 0$ and $\mu= n$ since we can interpolate via 
\[ \| a \|_{\mathring{L}^p(\mu)} \leq \big ( \| a \|_{\mathring{L}^p(n)} \big )^{\frac{\mu}{n}} \big ( \| a \|_{L^p} \big )^{1-\frac{\mu}{n}}  \]
For $\mu$ fixed either as $\mu=0$ or $\mu=n$, we need only bound the $\mathring{L}^p(\mu)$ norms $p = 1,2,\infty$ for times $t > 1$ and  $L^p$ norms for $p = 1, 3/2,\infty$ for times $t < 1$, and the result then follows from interpolation via
\[ \| a \|_{\mathring{L}^r(\mu)} \leq  \| a \|_{\mathring{L}^p(\mu)} \| a \|_{\mathring{L}^q(\mu)}   \]
for $r^{-1} = p^{-1} + q^{-1}$.

We begin by bounding the unweighted $L^p$ norms of the Duhamel term corresponding to $a(t)$ using our estimates above.  First we use Young's inequality, then split the integral into two parts:
\[ \begin{split}
\int \nolimits_0^{t} \Big \| \partial_t w(t-s) \ast & \partial_x^{\alpha} K_{\nu}(t-s) \ast \big [ \nabla \cdot N (a(s),\vec{\omega}(s) ) \big ] \Big \|_{L^{p}}ds \\ & \leq \int \nolimits_0^{t} \Big \| \partial_t w(t-s) \ast K_{\nu}(\frac{t-s}{2}) \Big \|_{L^q} \Big \| \partial_x^{\alpha} K_{\nu}(\frac{t-s}{2}) \ast \big [ \nabla \cdot N (a(s),\vec{\omega}(s) ) \big ] \Big \|_{L^{q_1}}ds \\  & \leq \Bigg ( \int_0^{t/2} + \int_{t/2}^{t} \Bigg ) (t-s)^{-\frac{3}{2}(\frac{1}{q_1}-\frac{1}{p}) }(1+t-s)^{\frac{1}{2}-(\frac{1}{q_1}-\frac{1}{p})} \Big \| \partial_x^{\alpha} K_{\nu}(\frac{t-s}{2}) \ast \big [ \nabla \cdot N (a(s),\vec{\omega}(s) ) \big ] \Big \|_{L^{q_1}}ds \\  & =: I_1 + I_2 
\end{split}
\] 
Here $1+ \frac{1}{p} = \frac{1}{q}+\frac{1}{q_1}$.  We can then bound the integrals for $s \in (0,t/2)$ and $s\in (t/2,t)$ separately.  

First we handle the $I_1$ term.  We use the heat estimate to pull the divergence and the $\partial_x^{\alpha}$ derivative off of the nonlinearity:
\[ \begin{split}
I_1 & \leq \int \nolimits_0^{t/2} (t-s)^{-\frac{3}{2}(\frac{1}{q_1}-\frac{1}{p})-\frac{1+|\alpha|}{2}}(1+t-s)^{\frac{1}{2}-(\frac{1}{q_1}-\frac{1}{p})} \big \|N (a(s),\vec{\omega}(s) ) \big \|_{L^{q_1}}ds \\ & \leq \max_{ijl} \int \nolimits_0^{t/2} (t-s)^{-\frac{3}{2}(\frac{1}{q_1}-\frac{1}{p})-\frac{1+|\alpha|}{2}}(1+t-s)^{\frac{1}{2}-(\frac{1}{q_1}-\frac{1}{p})} \big \| \partial_{x_i} ( m_j ) m_l \big \|_{L^{q_1}}ds
\end{split} \]
We can then use our above estimates on $\Pi$, $B$ in Cor. \ref{c:PiAndBiotSavartOperators} parts \textbf{(a)}, \textbf{(b)} to bound the nonlinear term:
\begin{equation}
\label{UnweightedNonlinearTermBound}
\begin{split}
\big \| \partial_{x_i}(m_j ) m_l \big \|_{L^{q_1}} & \leq \big \| \partial_{x_i} m_j \|_{L^{p_1}} \big \| m_l \big \|_{L^{p_2}} \leq  C \big ( \| a \|_{L^{p_1}} + \| \vec{\omega} \|_{\mathbb{L}^{p_1}} \big ) \big ( \| a \|_{L^{p_3}} + \| \vec{\omega} \|_{\mathbb{L}^{p_3}} \big ) \\ & \leq  C s^{-r_{0,p_1}-r_{0,p_3}}(1+s)^{-\min (\ell_{n,p_1,0},\tilde{\ell}_{n,p_1,0}) - \min (\ell_{n,p_3,0},\tilde{\ell}_{n,p_3,0})}\|(a,\vec{\omega})\|_{X_{n,k}}^2
\end{split}
\end{equation}
Note that the use of Young's inequality, H\"older's inequality, $(\ref{CalderonZygmundEst})$ and $(\ref{BiotSavartEst})$ puts the following restrictions on the set of admissible values for $p_1,p_3$:
\begin{equation} \label{Constraints} \begin{split}
1 < p_1 < \infty \hspace{0.75 cm} \text{ , } \hspace{0.75 cm} 1 < p_3 < 3 \hspace{0.75 cm} \text{ , } \hspace{0.75 cm}
\frac{1}{p} \leq \frac{1}{p_1}+ \frac{1}{p_3} -\frac{1}{3} \leq 1
\end{split} \end{equation}
We choose $q_1 = 1$ hence we require $\frac{1}{p_1}+ \frac{1}{p_3} -\frac{1}{3} = 1$.  Letting $p_1 = p_3 = \frac{3}{2}$ $(\ref{UnweightedNonlinearTermBound})$ becomes
\[ \begin{split}
\big \| \partial_{x_i}(m_j) m_l \big \|_{L^{1}} & \leq C(1+s)^{-\frac{2}{3} - \lfloor n \rfloor_1 } \|(a,\vec{\omega})\|_{X_{n,k}}^2
\end{split} \]
hence putting this together we have
\begin{equation} \label{I1Est} \begin{split}
I_1 & \leq \int \nolimits_0^{t/2} (t-s)^{-\frac{3}{2}(1-\frac{1}{p})-\frac{1+|\alpha|}{2} }(1+t-s)^{\frac{1}{2}-(1-\frac{1}{p})}(1+s)^{-\frac{2}{3}- \lfloor n \rfloor_1}  \|(a,\vec{\omega})\|_{X_{n,k}}^2 ds \\ & \leq C\|(a,\vec{\omega})\|_{X_{n,k}}^2 t^{-\frac{3}{2}(1-\frac{1}{p})-\frac{1+|\alpha|}{2}}(1+t)^{\frac{1}{2}-(1-\frac{1}{p})+\max (\frac{1}{3}- \lfloor n \rfloor_1,0) }L_{n,1/3}(t)
\end{split}
\end{equation}
for $t \geq 1$.  Thus the $L^p$ norms of $I_1$ have sufficiently fast decay as $t\to \infty$ for all $1\leq p\leq \infty$ such that the $X_{n,k}$ norm remains bounded.  For $t < 1$ we have
\[ \begin{split}
I_1 \leq \int_0^{t/2} (t-s)^{-\frac{3}{2}(1-\frac{1}{p})-\frac{1+|\alpha|}{2} }(1+t-s)^{\frac{1}{2}-(1-\frac{1}{p})}(1+s)^{-\frac{2}{3}- \lfloor n \rfloor_1}  \|(a,\vec{\omega})\|_{X_{n,k}}^2 ds \leq C\|(a,\vec{\omega})\|_{X_{n,k}}^2 t^{-\frac{3}{2}(1-\frac{1}{p})+\frac{1-|\alpha|}{2}}
\end{split}
\]
and hence we see the $L^p$ norms have the right behavior for $1\leq p \leq \infty$ such that the $X_{n,k}$ norms remain bounded.  Furthermore we note that for $1\leq p < 3/2$ and $|\alpha| = 0$ the $L^p$ norms tend to zero, which is consistent with the continuity of $F(a,\vec{\omega})$ at $t = 0$.

For $I_2$ we use the heat estimate to pull the divergence off of the nonlinearity:
\[ \begin{split}
I_2 & = \int \nolimits_{t/2}^t (t-s)^{-\frac{3}{2}(\frac{1}{q_1}-\frac{1}{p})}(1+t-s)^{\frac{1}{2}-(\frac{1}{q_1}-\frac{1}{p})} \big \| \partial_x^{\alpha} K_{\nu}(\frac{t-s}{2}) \ast \big [ \nabla \cdot N (a(s),\vec{\omega}(s) ) \big ] \big \|_{L^{q_1}}ds \\ & \leq \max_{ijl} \int \nolimits_{t/2}^t (t-s)^{-\frac{3}{2}(\frac{1}{q_1}-\frac{1}{p})-\frac{1}{2}}(1+t-s)^{\frac{1}{2}-(\frac{1}{q_1}-\frac{1}{p})} \big \|\partial_x^{\alpha} \partial_{x_i} ( m_j m_l ) \big \|_{L^{q_1}}ds
\end{split} \]
For an arbitrary multi-index $\beta$, we can use the estimates in Cor. \ref{c:PiAndBiotSavartOperators} parts \textbf{(a)}, \textbf{(b)} to obtain
\begin{equation}
\label{UnweightedDerivNonlinearTermBound}
\begin{split}
\big \| \partial_x^{\beta} \partial_{x_i}(m_j m_l) \big \|_{L^{q_1}} & \leq  \sum_{\gamma_1 + \gamma_2 = \beta} \big \| \partial_x^{\gamma_1} \partial_{x_i} m_j \|_{L^{p_1}} \big \|\partial_x^{\gamma_2} m_l \big \|_{L^{p_2}} \\ & \leq  C\sum_{\gamma_1 + \gamma_2 = \beta}  \big ( \|\partial_x^{\gamma_1} a \|_{L^{p_1}} + \| \partial_x^{\gamma_1}\vec{\omega} \|_{\mathbb{L}^{p_1}} \big ) \big ( \| \partial_x^{\gamma_2} a \|_{L^{p_3}} + \| \partial_x^{\gamma_2}\vec{\omega} \|_{\mathbb{L}^{p_3}} \big ) \\ & \leq  C s^{-r_{0,p_1}-r_{0,p_3}-\frac{|\beta|}{2}}(1+s)^{-\min (\ell_{n,p_1,0},\tilde{\ell}_{n,p_1,0}) - \min (\ell_{n,p_3,0},\tilde{\ell}_{n,p_3,0})}\|(a,\vec{\omega})\|_{X_{n,k}}^2
\end{split}
\end{equation}
provided that the constraints in $(\ref{Constraints})$ are met.  Here we take $\beta = \alpha$.  We must also ensure that the singularity at $s = t$ is integrable.  For $1\leq p < 3/2$ we can choose $p_1= p_3 = 3/2$ as before, and we obtain
\begin{equation} \label{I2Bound} \begin{split}
I_2 & \leq \int \nolimits_{t/2}^t (t-s)^{-\frac{3}{2}(1-\frac{1}{p})-\frac{1}{2}}(1+t-s)^{\frac{1}{2}-(1-\frac{1}{p})} s^{-\frac{|\alpha|}{2}} (1+s)^{-\frac{2}{3}-\lfloor n \rfloor_1} \|(a,\vec{\omega})\|_{X_{n,k}}^2 ds \\ & \leq  Ct^{-\frac{3}{2}(1-\frac{1}{p})+\frac{1-|\alpha|}{2} }(1+t)^{\frac{1}{2}-(1-\frac{1}{p})-\frac{2}{3}-\lfloor n \rfloor_1 } \|(a,\vec{\omega})\|_{X_{n,k}}^2
\end{split}
\end{equation}
for $0 < t < \infty$, hence these $L^p$ norms have the right behavior as $t \to 0$ and as $t\to \infty$, and tend to zero for $|\alpha|=0$ which is consistent with continuity at $t =0$.  Similarly, for $3/2 \leq p \leq 2$ we can choose $p_1 = p_3 = 2$ in $(\ref{UnweightedDerivNonlinearTermBound})$ and obtain the pointwise bound
\[ \big \| \partial_x^{\alpha}\partial_{x_i}(m_j m_l )\big \|_{L^{\frac{3}{2}}} \leq Cs^{-\frac{1}{2}-\frac{|\alpha|}{2}}(1+s)^{-1-\lfloor n \rfloor_1 } \|(a,\vec{\omega})\|_{X_{n,k}}^2 \]
from which it follows that
\[ \begin{split}
I_2 & \leq \int \nolimits_{t/2}^t (t-s)^{-\frac{3}{2}(\frac{2}{3}-\frac{1}{p})-\frac{1}{2}}(1+t-s)^{\frac{1}{2}-(\frac{2}{3}-\frac{1}{p})}s^{-\frac{1}{2}-\frac{|\alpha|}{2} }(1+s)^{-1-\lfloor n \rfloor_1} \|(a,\vec{\omega})\|_{X_{n,k}}^2 ds \\ & \leq  Ct^{-\frac{3}{2}(\frac{2}{3}-\frac{1}{p})-\frac{|\alpha|}{2}}(1+t)^{\frac{1}{2}-(1-\frac{1}{p})-\frac{2}{3}-\lfloor n \rfloor_1} \|(a,\vec{\omega})\|_{X_{n,k}}^2
\end{split}
\]
for $0 < t < \infty$, hence these $L^{p}$ norms also have the right behavior as $t \to 0$ and as $t\to \infty$.  Finally, we can obtain bounds on the $L^{\infty}$ norm by choosing $p_1 = 8$, $p_3 = 8/3$ in $(\ref{UnweightedDerivNonlinearTermBound})$ to obtain the pointwise bound
\[ \big \| \partial_x^{\alpha}\partial_{x_i}(m_j m_l )\big \|_{L^{6}} \leq Cs^{-\frac{5}{4}-\frac{|\alpha|}{2}}(1+s)^{-1-\lfloor n \rfloor_1 } \|(a,\vec{\omega})\|_{X_{n,k}}^2 \]
from which we then obtain the following bound on the integral for $0 < t < \infty$:
\begin{equation} \label{SlowInfinityBound} \begin{split}
I_2 \leq  Ct^{-1-\frac{|\alpha|}{2}}(1+t)^{-1+\frac{1}{3}-\lfloor n \rfloor_1} \|(a,\vec{\omega})\|_{X_{n,k}}^2
\end{split}
\end{equation}
Note this is slower than the linear evolution rate.  For $|\alpha| < k$ we can make an improved estimate to match the linear rate as follows.  With $p=\infty$, we keep all derivatives on the nonlinearity when using the heat estimate, and we obtain
\begin{equation} \label{LargeTimeAltEst}
I_2 \leq \max_{ij} \int \nolimits_{t/2}^t (t-s)^{-\frac{3}{2q_1}}(1+t-s)^{\frac{1}{2}-\frac{1}{q_1}} \big \|\partial_x^{\alpha} \partial_{x_i} \partial_{x_j} ( m_i m_j ) \big \|_{L^{q_1}}ds 
\end{equation}
We can then use the estimate in \ref{UnweightedDerivNonlinearTermBound} by taking $\beta = \alpha + e_j$, and we choose $p_1 = p_3 = 12/5$ to obtain
\[I_2 \leq C\int \nolimits_{t/2}^{t} (t-s)^{-\frac{3}{4}}s^{-\frac{3}{4}-\frac{|\alpha|+1}{2}}(1+s)^{-1-\lfloor n \rfloor_1} ds \|(a,\vec{\omega})\|_{X_{n,k}}^2 \leq C\|(a,\vec{\omega})\|_{X_{n,k}}^2t^{-1-\frac{|\alpha|}{2}}(1+t)^{-1-\lfloor n\rfloor_1} \]

For $n = 0$ we are done.  For $n > 0$ we bound the weighted norms when $\mu = n$ of the Duhamel term corresponding to $a(t)$, and the results then follow by interpolation.  We first bound the weighted norm of the convolution in terms of the weighted norms of each of its components using Young's inequality:
\begin{equation} \label{WeightedEstimateDuhamel} \begin{split}
\int \nolimits_0^t \Big \| \partial_t & w(t-s) \ast \partial_x^{\alpha} K_{\nu}(t-s) \ast \big [ \nabla \cdot N (a(s),\vec{\omega}(s) ) \big ] \Big \|_{\mathring{L}^p(n)}ds \\  & \leq \int \nolimits_0^t \Big \| \partial_t w(t-s) \ast K_{\nu}(\frac{t-s}{2}) \Big \|_{\mathring{L}^{\tilde{q}}(n)} \Big \|\partial_x^{\alpha} K_{\nu}(\frac{t-s}{2})\ast  \big [ \nabla \cdot N (a(s),\vec{\omega}(s) ) \big ] \Big \|_{L^{\tilde{q}_1}}ds \\ & \hspace{1 cm} + \int \nolimits_0^t \Big \| \partial_t w(t-s) \ast K_{\nu}(\frac{t-s}{2}) \Big \|_{L^q} \Big \| \partial_x^{\alpha} K_{\nu}(\frac{t-s}{2})\ast \big [ \nabla \cdot N (a(s),\vec{\omega}(s) ) \big ] \Big \|_{\mathring{L}^{q_1}(n)}ds 
\end{split}
\end{equation}
For the first term, we can use the weighted estimate of the heat-wave operator in Prop \ref{p:HypParaEst} and then repeat the analysis used above for the unweighted norm of the nonlinearity line by line to obtain the appropriate bounds for this term.  So we need only bound the second term.

For the second term we use the unweighted estimate in Prop \ref{p:HypParaEst} and split the integral as before:
\[ \begin{split}
\int \nolimits_0^t \Big \| \partial_t w(& t-s)  \ast K_{\nu}(\frac{t-s}{2}) \Big \|_{L^q} \Big \| \partial_x^{\alpha} K_{\nu}(\frac{t-s}{2})\ast \big [ \nabla \cdot N (a(s),\vec{\omega}(s) ) \big ] \Big \|_{\mathring{L}^{q_1}(n)}ds \\ & \leq
\int \nolimits_0^t (t-s)^{-\frac{3}{2}(\frac{1}{q_1}-\frac{1}{p})} (1+t-s)^{\frac{1}{2}-(\frac{1}{q_1}-\frac{1}{p})} \Big \| \partial_x^{\alpha} K_{\nu}(\frac{t-s}{2})\ast  \big [ \nabla \cdot N (a(s),\vec{\omega}(s) ) \big ] \Big \|_{\mathring{L}^{q_1}(n)}ds \\ & \leq \Bigg ( \int_0^{t/2} + \int_{t/2}^{t} \Bigg ) (t-s)^{-\frac{3}{2}(\frac{1}{q_1}-\frac{1}{p})}(1+t-s)^{\frac{1}{2}-(\frac{1}{q_1}-\frac{1}{p})} \Big \| \partial_x^{\alpha} K_{\nu}(\frac{t-s}{2})\ast \big [ \nabla \cdot N (a(s),\vec{\omega}(s) ) \big ] \Big \|_{\mathring{L}^{q_1}(n)}ds \\ & = I_1 + I_2
\end{split} \]
The next step is to use our heat estimate, and then we will need bounds for the weighted norm of the nonlinear term analogous to (\ref{UnweightedNonlinearTermBound}), (\ref{UnweightedDerivNonlinearTermBound}).  Note however that these bounds are essentially the same, so here we will derive both at once.  The derivation is similar to (\ref{UnweightedDerivNonlinearTermBound}), but one must always place the weight on the term with fewer derivatives in order to use Cor. \ref{c:PiAndBiotSavartOperators} part \textbf{(a)}.  For $0 < n < 2$ we make the estimate
\begin{equation} 
\label{WeightedNonlinearTermBound} \begin{split}
\|\partial_x^{\beta} \partial_{x_i} (m_j m_l) \|_{\mathring{L}^{q_1}(n)} & \leq \| \partial_x^{\beta} \partial_{x_i} (m_j m_l) \|_{L^{q_1}(n)} \\ & \leq C\sum_{\gamma_1+\gamma_2 = \beta} \big ( \|\partial_x^{\gamma_1}a\|_{L^{p_1}} + \|\partial_x^{\gamma_1}\vec{\omega}\|_{\mathbb{L}^{p_1}} \big ) \big ( \|\partial_x^{\gamma_2} a \|_{L^{p_3}(n)} + \| \partial_x^{\gamma_2} \vec{\omega} \|_{\mathbb{L}^{p_3}(n)} \big ) \\ & \leq C s^{-r_{0,p_1}-r_{0,p_3}-\frac{|\beta|}{2}}(1+s)^{-\min (\ell_{n,p_1,0},\tilde{\ell}_{n,p_1,0}) - \min (\ell_{n,p_3,n},\tilde{\ell}_{n,p_3,n})}\|(a,\vec{\omega})\|_{X_{n,k}}^2 \end{split} \end{equation}
using parts \textbf{(a)} and \textbf{(b)} of Cor \ref{c:PiAndBiotSavartOperators}, which requires the set of constraints 
\begin{equation} \label{Constraints2} \begin{split}
1 < p_1 < \infty \hspace{0.75 cm} \text{ , } \hspace{0.75 cm}
\frac{3}{3-n} < p_3 < \frac{3}{1-\lfloor n \rfloor_1} \hspace{0.75 cm} \text{ , } \hspace{0.75 cm}
\frac{1}{p} \leq \frac{1}{p_1}+ \frac{1}{p_3} -\frac{1}{3} \leq 1
\end{split} \end{equation}
or for $1\leq n \leq 2$ we can obtain the same bound using parts \textbf{(a)} and \textbf{(c)} of Cor \ref{c:PiAndBiotSavartOperators}, which require
\begin{equation} \label{Constraints3} \begin{split}
1 < p_1 < \infty \hspace{0.75 cm} \text{ , } \hspace{0.75 cm} \frac{3}{4-n} < p_3 < \frac{3}{3-n} \hspace{0.75 cm} \text{ , } \hspace{0.75 cm}
\frac{1}{p} \leq \frac{1}{p_1}+ \frac{1}{p_3} -\frac{1}{3} \leq 1
\end{split} \end{equation}
Note that in the overlapping region $1 \leq n < 2$ we can use either bound, but if we use Cor \ref{c:PiAndBiotSavartOperators} \textbf{(a)} and \textbf{(c)} by satisfying the constraints in (\ref{Constraints3}), we are allowed to choose smaller $p_3$ than (\ref{Constraints2}) allow, a fact which we will exploit.  The task then becomes obtaining various choices of $p_1$ and $p_3$ for $I_1$, $I_2$, $1\leq p \leq \infty$, $0 < n \leq 2$.  


For $I_1$ we use the heat estimate to pull the divergence and the $\partial_x^\alpha$ derivative off of the nonlinearity, and use $(\ref{WeightedNonlinearTermBound})$ with $\beta = 0$.  For $0 < n < 1$ we can satisfy the constraints in $(\ref{Constraints2})$ with $q_1 =1$ by taking $p_1 = p_3 = 3/2$, and we obtain
\[ \begin{split}
I_1 & \leq \int \nolimits_0^{t/2} (t-s)^{-\frac{3}{2}(1-\frac{1}{p})-\frac{1+|\alpha|}{2}}(1+t-s)^{\frac{1}{2}-(1-\frac{1}{p})} (1+s)^{-\frac{2}{3}-\lfloor n \rfloor_1 + n}\|(a,\vec{\omega})\|_{X_{n,k}}^2 ds \\ & \leq C\|(a,\vec{\omega})\|_{X_{n,k}}^2 t^{-\frac{3}{2}(1-\frac{1}{p})-\frac{1+|\alpha|}{2}}(1+t)^{\frac{1}{2}-(1-\frac{1}{p})+\frac{1}{3}-\lfloor n \rfloor_1 + n}
\end{split}
\]
whereas for $1 < n < 2$ precisely the same estimate holds by taking $p_1 = p_3 = 3/2$ in (\ref{Constraints3}).  Hence these weighted $L^p$ norms decay sufficiently quickly as $t\to \infty$ for $1\leq p \leq \infty$.  For $t < 1$ this bound becomes
\begin{equation} \label{SmallTimeBehavior} I_1 \leq C\|(a,\vec{\omega})\|_{X_{n,k}}^2 t^{-\frac{3}{2}(1-\frac{1}{p})+\frac{1-|\alpha|}{2}} \end{equation}
hence these norms have the right behavior as $t \to 0$.  For $1\leq n < 3/2$ we can use (\ref{Constraints3}) by taking $p_1 = 2$, $p_3 =6/5$ and for $3/2 < n \leq 2$ we can use $p_1 = 6/5$, $p_3 =2$.  In both cases we have 
\[ \begin{split}
I_1 & \leq \int \nolimits_0^{t/2} (t-s)^{-\frac{3}{2}(1-\frac{1}{p})-\frac{1+|\alpha|}{2}}(1+t-s)^{\frac{1}{2}-(1-\frac{1}{p})} s^{-\frac{1}{4}}(1+s)^{-\frac{5}{6}+n-\frac{7}{12}}\|(a,\vec{\omega})\|_{X_{n,k}}^2 ds \\ & \leq C\|(a,\vec{\omega})\|_{X_{n,k}}^2 t^{-\frac{3}{2}(1-\frac{1}{p})+\frac{1}{4}-\frac{|\alpha|}{2}}(1+t)^{\frac{1}{2}-(1-\frac{1}{p})+n-\frac{4}{3}}
\end{split}
\]
for $0< t < \infty$, hence the weighted $L^p$ norms of this term decay sufficiently fast to remain in $X_{n,k}$ for $1\leq p \leq \infty$.  For $t < 1$ this bound shows that the $\mathring{L}^p(n)$ norms have the right behavior as $t \to 0$ for $1 \leq p < 6/5$.  Then we need only prove that the $L^p(n)$ norms for $6/5 \leq p \leq \infty$ have the right behavior as $t\to 0$ for $n=1$ and $n=2$.  Here we can choose $p_1 = 3/2$, $p_3 = 2$ for $n=1$ using (\ref{Constraints2}) and using (\ref{Constraints3}) for $n=2$ and we again obtain (\ref{SmallTimeBehavior}), so the weighted $L^p$ norms blow up sufficiently slowly for $6/5\leq p \leq \infty$ as $t\to 0$, hence $I_1$ belongs to $X_{n,k}$. 

For $I_2$ we can reuse many of the estimates in the unweighted case, but we have to modify these slightly.  We again use the heat estimate to pull the divergence off the nonlinearity, and we again have to worry about the singularity at $s = t$.  For $0 < n < 1$ we can make the precisely the same choices as in the unweighted case.  Namely that we can obtain the appropriate bounds for the $\mathring{L}^p(n)$ norms using (\ref{Constraints2}) by taking $p_1 = p_3 = 3/2$ for $1\leq p < 3/2$ and we obtain the analogous weighted pointwise bound
\[ \begin{split}
\|\partial_x^{\alpha} \partial_{x_i} (m_j m_l) \|_{\mathring{L}^{1}(n)}& \leq Cs^{-\frac{|\alpha|}{2}}(1+s)^{-\frac{5}{3} +n } \|(a,\vec{\omega})\|_{X_{n,k}}^2
\end{split} \]
We can then make the identical estimate in (\ref{I2Bound}) with this analogous pointwise bound to show that these norms have the correct behavior for $0<t<\infty$.  Similarly we can use (\ref{Constraints2}) by taking $p_1 = p_3 = 2$ for $3/2 \leq p \leq 2$ and taking $p_1 = 8$, $p_3 = 8/3$ for $p = \infty$ and obtain the analogous pointwise bounds, from which it follows in the same way that these norms have the correct behavior for $0<t<\infty$, except for $p=\infty$, $|\alpha|<k$.  We can then match the decay rate for $p=\infty$, $|\alpha|<k$ by keeping all derivatives on the nonlinearity as in (\ref{LargeTimeAltEst}), taking $\beta = \alpha + e_j$ in (\ref{WeightedNonlinearTermBound}) and taking $p_1 = p_3 = 12/5$ in (\ref{Constraints2}).

The case $1 < n < 2$ is also similar, and we can show that the $\mathring{L}^p(n)$ norms have the correct behavior for $1 \leq p < 3/2$ by taking $p_1 = p_3 = 3/2$ in $(\ref{Constraints3})$.  For the $\mathring{L}^p(n)$ norms for $3/2 \leq p \leq 2$ we make a slightly different estimate by taking $p_1 = 3$, $p_3 =3/2$ in $(\ref{Constraints3})$ and we obtain the pointwise bound
\[ \begin{split}
\|\partial_x^{\alpha} \partial_{x_i} (m_j m_l) \|_{\mathring{L}^{\frac{3}{2}}(n)} & \leq Cs^{-\frac{1}{2}-\frac{|\alpha|}{2}}(1+s)^{-\frac{11}{6} +n } \|(a,\vec{\omega})\|_{X_{n,k}}^2
\end{split} \]
and repeating the above analysis.  For $1 < n < 2$ we can set $q_1 = 6$ by choosing $p_1 = 8/3$ and $p_3 = 8$ using $(\ref{Constraints2})$, and show that the $\mathring{L}^{\infty}(n)$ norms have the correct behavior for $0<t<\infty$, except for $p=\infty$, $|\alpha|<k$.  We can then match the decay rate for $p=\infty$, $|\alpha|<k$ by keeping the derivatives on the nonlinearity and using $\beta = \alpha + e_j$ in (\ref{WeightedNonlinearTermBound}) with $p_1 = 2$, $p_3 = 3$ in (\ref{Constraints2}).

It remains to show the $\mathring{L}^p(n)$ norms have the correct behavior for $n=1$ and $n=2$.  We can choose $p_1 = 2$, $p_3 = 6/5$ for $1\leq n < 3/2$ and $p_1 = \frac{6}{5}$, $p_3 = 2$ for $3/2 < n \leq 2$ and we find
\[ \begin{split}
I_2 & \leq \int \nolimits_{t/2}^t (t-s)^{-\frac{3}{2}(1-\frac{1}{p})-\frac{1}{2}}(1+t-s)^{\frac{1}{2}-(1-\frac{1}{p})}s^{-\frac{1}{4}-\frac{|\alpha|}{2}}(1+s)^{-\frac{17}{12}+n} \|(a,\vec{\omega})\|_{X_{n,k}}^2 ds \\ & \leq  Ct^{-\frac{3}{2}(1-\frac{1}{p})+\frac{1}{4}-\frac{|\alpha|}{2}} (1+t)^{\frac{1}{2}-(1-\frac{1}{p})-\frac{17}{12}+n} \|(a,\vec{\omega})\|_{X_{n,k}}^2
\end{split}
\]
for $1\leq p <3/2$, which decays appropriately quickly as $t \to \infty$.  Note also that this bound holds for $t < 1$, and hence the weighted $L^p$ norms tend to zero as $t \to 0$ for $1 \leq p < 6/5$.  For $3/2 \leq p \leq 2$ we can set $q_1 = 3/2$ by choosing $p_1 = p_3 = 2$ using (\ref{Constraints2}) for $1\leq n < 3/2$ and (\ref{Constraints3}) for $3/2 < n \leq 2$ and we find
\[ \begin{split}
I_2 & \leq Ct^{-\frac{3}{2}(\frac{2}{3}-\frac{1}{p})-\frac{|\alpha|}{2}}(1+t)^{\frac{1}{2}-(1-\frac{1}{p})-\frac{5}{3}+ n} \|(a,\vec{\omega})\|_{X_{n,k}}^2
\end{split} \]
for $0 < t < \infty$.  Finally, for $1 \leq n < 3/2$ we choose $p_1 = 8/3$, $p_3 = 8$ using $(\ref{Constraints2})$ and for $15/8 < n \leq 2$ we choose $p_1 = 8$, $p_3 = 8/3$ using $(\ref{Constraints3})$ and we see that the $\mathring{L}^{\infty}(n)$ norm has the right behavior for $t>1$ $|\alpha| = k$ and $t<1$ for all $\alpha$, and we can then match the decay rate for $p=\infty$, $|\alpha|<k$ by keeping the derivatives on the nonlinearity and using $\beta = \alpha + e_j$ in (\ref{WeightedNonlinearTermBound}) with $p_1 = 2$, $p_3 = 3$ in (\ref{Constraints2}) for $1\leq n < 3/2$ and $p_1 = 24/11$, $p_3 = 8/3$ in (\ref{Constraints3}) for $15/8 < n \leq 2$.

The bounds on the Duhamel term for $\vec{\omega}(t)$ can be obtained in a very similar manner.  The only difference is that one need not make the initial step of using Young's inequality.  Namely, we begin by looking at the unweighted norms, and we first split the integral
\[ \begin{split}
\int_0^{t} \Big \| \partial_x^{\alpha}  \mathbb{K}_{\epsilon}(t-s) \ast \big [ \nabla \times N (a(s),\vec{\omega}(s) ) \big ] \Big \|_{\mathbb{L}^p}ds & = \Bigg ( \int_0^{t/2} + \int_{t/2}^{t} \Bigg ) \Big \| \partial_x^{\alpha}  \mathbb{K}_{\epsilon}(t-s) \ast \big [ \nabla \times N (a(s),\vec{\omega}(s) ) \big ] \Big \|_{\mathbb{L}^p}ds  \\ & =: I_1 + I_2
\end{split}
\]
We can then use the heat estimate directly, and for $s\in (0,t/2)$ we pull the divergence and the $\partial_x^{\alpha}$ derivative off the nonlinear term using the heat estimate, whereas $s\in (t/2,t)$ we only pull the divergence off.  By making the exact same estimates as for the Duhamel term for $a(t)$ with the same choices of $p_1$ and $p_3$ we arrive at the analogous bounds.  The weighted norms can be obtained in the same way.  For brevity we omit this.


It remains to obtain continuity for $t > 0$, in which case we would have $F(a,\vec{\omega}) \in Z_{n}^0\cap Z_{n,k}^+$.  Beginning with the Duhamel term for $a(t)$, we note that this is equivalent to showing that
\begin{align} \label{Continuity} & \lim_{h\to 0} \int \nolimits_t^{t+h} \Big \| \partial_t w(t+h-s) \ast \partial_x^{\alpha} K_{\nu}(t+h-s) \ast \big [ \nabla \cdot N (a(s),\vec{\omega}(s) ) \big ] \Big \|_{L^p}ds = 0  \\ & \lim_{h\to 0} \int \nolimits_0^t \Big \|  \big [ \partial_t w(t+h-s) \ast   \partial_x^{\alpha} K_{\nu}(t+h-s)- \partial_t w(t-s) \ast   \partial_x^{\alpha} K_{\nu}(t-s) \big ] \ast \big [ \nabla \cdot N (a(s),\vec{\omega}(s) ) \big ] \Big \|_{L^p}ds = 0 \notag \end{align}
For the first limit we can re-use the methods used to obtain a bound on the $I_2$ term above to show that this limit is zero.  For the second, we can use the estimate
\[ \begin{split}
\Big \|  \big [ \partial_t & w(t+h-s) \ast \partial_x^{\alpha}  K_{\nu}(t+h-s)- \partial_t w(t-s) \ast \partial_x^{\alpha}   K_{\nu}(t-s) \big ] \ast \big [ \nabla \cdot N (a(s),\vec{\omega}(s) ) \big ] \Big \|_{L^p} \\ & \leq \Big \| \partial_t w(t+h-s) \ast \partial_x^{\alpha}   K_{\nu}(\frac{t-s}{2}+h)- \partial_t w(t-s) \ast  \partial_x^{\alpha} K_{\nu}(\frac{t-s}{2}) \Big \|_{L^1} \Big \| K_{\nu}(\frac{t-s}{2}) \ast \big [ \nabla \cdot N (a(s),\vec{\omega}(s) ) \big ] \Big \|_{L^p}
\end{split} \]
and show that this first factor tends to zero uniformly in $s$ as $h\to 0$.  The weighted norms can be bounded similarly, and one can obtain continuity for the Duhamel term corresponding to $\vec{\omega}(t)$ by showing that the limits analogous to $(\ref{Continuity})$ are zero.

\paragraph{Claim Two: $F$ has Lipschitz constant $K = \frac{1}{2}$ on a ball $B(0,R)$ in $X_{n,k}$ }

We need to bound $ \| F(a,\vec{\omega}) - F(\tilde{a},\tilde{\vec{\omega}}) \|_{X_{n,k}}$ for $(a,\vec{\omega}),(\tilde{a},\tilde{\vec{\omega}}) \in B(0,R)$, where $R$ is yet to be chosen.  The analysis is similar to the above, but now we use the bilinear property of the nonlinearity to get the analogous unweighted estimates
\begin{align} \label{UnweightedNonlinearDifferenceBound}  
\big \| \partial_x^{\beta} N & ( a(s),\vec{\omega}(s) ) - \partial_x^{\beta} N(\tilde{a}(s),\tilde{\vec{\omega}}(s)) \big \|_{L^{q_1}} \notag \\ & \leq \max_{ijl} \big \| \partial_x^{\beta} \big [\partial_{x_i}(m_j) (m_l-\tilde{m}_l) \big ] \big \|_{L^{q_1}} + \big \|\partial_x^{\beta} \big [ \partial_{x_i}(m_j-\tilde{m}_j)\tilde{m}_l \big ] \big \|_{L^{q_1}} \\ \leq  C \Big ( \|(\tilde{a}, \tilde{\vec{\omega}})\|_{X_{n,k}} & +  \|(a,\vec{\omega})\|_{X_{n,k}} \Big ) \|(a-\tilde{a},\vec{\omega}-\tilde{\vec{\omega}}) \|_{X_{n,k}} s^{-r_{0,p_1}-r_{0,p_3}-\frac{|\beta|}{2}}(1+s)^{-\min (\ell_{n,p_1,0},\tilde{\ell}_{n,p_1,0}) - \min (\ell_{n,p_3,0},\tilde{\ell}_{n,p_3,0})} \notag \end{align}
corresponding to (\ref{UnweightedNonlinearTermBound}) and (\ref{UnweightedDerivNonlinearTermBound}), which require the set of constraints (\ref{Constraints}), as well as the analogous weighted estimate
\begin{align} \label{WeightedNonlinearDifferenceBound} \big \| & \partial_x^{\beta} \big [ N (a(s), \vec{\omega}(s) ) - N(\tilde{a}(s),\tilde{\vec{\omega}}(s)) \big ] \big \|_{\mathring{L}^{q_1}(n)} \\ \leq  C \big ( \|(\tilde{a},& \tilde{\vec{\omega}})\|_{X_{n,k}} +  \|(a,\vec{\omega})\|_{X_{n,k}} \big ) \|(a-\tilde{a},\vec{\omega}-\tilde{\vec{\omega}}) \|_{X_{n,k}} s^{-r_{0,p_1}-r_{0,p_3}-\frac{|\beta|}{2}}(1+s)^{-\min (\ell_{n,p_1,0},\tilde{\ell}_{n,p_1,0}) - \min (\ell_{n,p_3,n},\tilde{\ell}_{n,p_3,n})} \notag \end{align}
corresponding to (\ref{WeightedNonlinearTermBound}) which requires the set of constraints (\ref{Constraints2}) for $0 < n < 2$ and (\ref{Constraints3}) for $1 \leq n \leq 2$.

The proof then follows exactly the steps used to prove Claim 1 with these analogous estimates.  We begin by looking at the norms of the difference between the Duhamel term corresponding to $a(t)$:
\[ \begin{split}
& \int \nolimits_0^t \Big \| \partial_t w(t-s) \ast \partial_x^{\alpha} K_{\nu}(t-s) \ast \big [ \nabla \cdot \big [ N (a(s),\vec{\omega}(s) ) - N(\tilde{a}(s),\tilde{\vec{\omega}}(s)) \big ] \big ] \Big \|_{L^p}ds \\  & \leq 
\Bigg (\int_0^{t/2} + \int_{t/2}^t \Bigg ) (t-s)^{-\frac{3}{2}(\frac{1}{q_1}-\frac{1}{p})}(1+t-s)^{\frac{1}{2}-(\frac{1}{q_1}-\frac{1}{p})} \Big \| \partial_x^{\alpha} K_{\nu}(\frac{t-s}{2})\ast  \big [ \nabla \cdot \big [ N (a(s),\vec{\omega}(s) ) - N(\tilde{a}(s),\tilde{\vec{\omega}}(s)) \big ]  \big ] \Big \|_{L^{q_1}}ds \\ & =: I_1 + I_2
\end{split}
\]
For $I_1$ we can then use the heat estimate and the bilinearity to obtain 
\[ I_1 \leq C\max_{ijk} \int \nolimits_0^{t/2} (t-s)^{-\frac{3}{2}(\frac{1}{q_1}-\frac{1}{p})-\frac{1}{2}-\frac{|\alpha|}{2}}(1+t-s)^{\frac{1}{2}-(\frac{1}{q_1}-\frac{1}{p})} \Big ( \big \|\partial_{x_i}(m_j) (m_l-\tilde{m}_l) \big \|_{L^{q_1}} + \big \|\partial_{x_i}(m_j-\tilde{m}_j) \tilde{m}_l \big \|_{L^{q_1}} \Big )ds 
\]
We can then repeat the analysis for the Duhamel term above for $a(t)$ line by line for each of these terms, using  $(\ref{UnweightedNonlinearDifferenceBound})$ with $\alpha = 0$ and then making the same choices for $p_1$ and $p_3$ to handle the cases $t\geq 1$ and $t <1$ separately for different values of $p$, and we find
\begin{equation} \label{ContractionFirstDuhamelUnweighted} \begin{split}
\sup_{|\alpha|\leq k} \sup_{1\leq p \leq \infty} \sup_{0 \leq t < \infty} t^{r_{\alpha,p}} (1+t)^{\ell_{n,p,0}+\hat{\ell}_{k,p,\alpha}}I_1 \leq  C\Big ( \|(\tilde{a},\tilde{\vec{\omega}})\|_{X_{n,k}} +  \|(a,\vec{\omega})\|_{X_{n,k}} \Big ) \|(a-\tilde{a},\vec{\omega}-\tilde{\vec{\omega}}) \|_{X_{n,k}}
\end{split}
\end{equation}
Similarly for $I_2$ we can use the heat estimate, $(\ref{UnweightedNonlinearDifferenceBound})$ and the preceding analysis to obtain
\begin{equation} \label{ContractionFirstDuhamelUnweighted2}\begin{split}
\sup_{|\alpha|\leq k} \sup_{1\leq p \leq \infty} \sup_{0 \leq t < \infty} t^{r_{\alpha,p}}(1+t)^{\ell_{n,p,0}+\hat{\ell}_{k,p,\alpha}} I_2 \leq C\Big ( \|(\tilde{a},\tilde{\vec{\omega}})\|_{X_{n,k}} +  \|(a,\vec{\omega})\|_{X_{n,k}} \Big ) \|(a-\tilde{a},\vec{\omega}-\tilde{\vec{\omega}}) \|_{X_{n,k}}
\end{split}
\end{equation}
The bounds on the weighted norms can be obtained by following the steps used in the proof of Claim 1 with the analogous bound (\ref{WeightedNonlinearDifferenceBound}), and the bounds on the Duhamel term for $\vec{\omega}(t)$ can be obtained by repeating this procedure.  By combining  $(\ref{ContractionFirstDuhamelUnweighted})$, $(\ref{ContractionFirstDuhamelUnweighted2})$, the bounds on the weighted norms and the analogue for the Duhamel term for $\vec{\omega}(t)$, we obtain
\[ \Big \| F(a,\vec{\omega}) - F(\tilde{a},\tilde{\vec{\omega}}) \Big \|_{X_{n,k}} \leq C \Big ( \|(\tilde{a},\tilde{\vec{\omega}})\|_{X_{n,k}} + \|(a,\vec{\omega})\|_{X_{n,k}} \Big ) \big \| (a-\tilde{a},\vec{\omega}-\tilde{\vec{\omega}} ) \big \|_{X_{n,k}} \]
so by letting $R = \frac{1}{4C}$ we have our result.

\end{proof}

Having proven the existence of solutions $a(t)$ and $\vec{\omega}(t)$, we now complete the proof of existence of solutions to $(\ref{IntegralForm})$ by proving the existence of a solution $\rho(t)$.  For $n \in \mathbb{R}_{\geq 0 }$ we define the function space
\[ Y_{n,k} = \{\rho : \rho \in \bigcap_{1\leq p < 3/2} C^0\big [ [0,\infty), L^p(n) \big ] \text{ and }  \rho \in \bigcap_{1\leq p \leq \infty} C^0\big [ (0,\infty), W^{k,p}(n) \big ] \} \]
equipped with the norm
\[ \begin{split}
\|\rho \|_{Y_{n,k}} = \sup_{|\alpha|\leq k} \sup_{1\leq p \leq \infty} \sup_{0 \leq \mu \leq n} \sup_{0 < t < \infty}  \Big [ t^{r_{\alpha,p}}(1+t)^{\ell_{n,p,\mu} +\hat{\ell}_{k,p,\alpha}- \frac{1}{2}} \|\partial_{x}^{\alpha}\rho(t)\|_{\mathring{L}^p(\mu)} \Big ] \end{split} \]
where $r_{\alpha,p}, \ell_{n,p,\mu}, \hat{\ell}_{k,p,\alpha}$ are as before.

\begin{Corollary}\label{c:ExistencePt2}
	Fix $n \in [ 0,2 ]$, $k\geq 1$ and let $(\rho_0,a_0,\vec{\omega}_0)$ belong to $W^{1,p}(n)\times L^{p}(n)\times \mathbb{L}_{\sigma}^{p}(n)$ for all $1 \leq p \leq \frac{3}{2}$, where $a_0,\vec{\omega}_0$ have zero total mass and $(\rho_0,a_0,\vec{\omega}_0)$ have sufficiently small norms as in Theorem $\ref{t:Existence}$.  If $(a(t),\vec{\omega}(t))$ is the solution of $(\ref{IntegralForm})$ from Theorem $\ref{t:Existence}$, then the solution $\rho(t)$ defined by $(\ref{IntegralForm})$ belongs to $Y_{n,k}$.
\end{Corollary}

\begin{proof}
	
	As before the decay rates and smoothness properties are chosen to match those of the linear terms hence we need only check the Duhamel term.  We first estimate the unweighted norms
	\[ \begin{split}
	\int \nolimits_0^t & \Big \| w(t-s) \ast \partial_x^{\alpha} K_{\nu}(t-s) \ast  \big [ \nabla \cdot N (a(s),\vec{\omega}(s) ) \big ] \Big \|_{L^p}ds \\  & \leq \max_{ijk} \Bigg ( \int_0^{t/2} + \int_{t/2}^t  \Bigg ) (t-s)^{-\frac{3}{2}(\frac{1}{q_1}-\frac{1}{p}) +1}(1+t-s)^{-(\frac{1}{q_1}-\frac{1}{p})} \Big \| \partial_x^{\alpha} K_{\nu}(\frac{t-s}{2}) \ast \big [ \nabla \cdot N (a(s),\vec{\omega}(s) ) \big ] \Big \|_{L^p}ds =: I_1 + I_2
	\end{split}
	\]
	
	For $I_1$, we pull the divergence and the $\partial_x^{\alpha}$ derivative off of the nonlinearity using the heat estimate, use estimate $(\ref{UnweightedNonlinearTermBound})$, let $p_1 = p_3 = 3/2$ and find
	\[ I_1 \leq C\|(a,\vec{\omega})\|_{X_{n,k}}^2 t^{-\frac{3}{2}(1-\frac{1}{p})+\frac{1}{2}-\frac{|\alpha|}{2}}(1+t)^{-(1-\frac{1}{p})+\max (\frac{1}{3}-\lfloor n \rfloor_1, 0 )}  \]
	which holds for all $t>0$, hence the $L^p$ norms of this term have sufficiently fast decay for $1\leq p \leq \infty$ as $t\to \infty$, tend to zero as $t\to 0$ for $1\leq p < 3/2$, $|\alpha|= 0$ and blow up sufficiently slowly for $3/2 \leq p \leq \infty$.  
	
	For $I_2$, we use the heat estimate to pull the divergence off the nonlinearity, use estimate $(\ref{UnweightedDerivNonlinearTermBound})$ and set $p_1 = p_3 = 3/2$ for $1\leq p \leq 2 $ and find
	\[ I_2 \leq C\|(a,\vec{\omega})\|_{X_{n,k}}^2 t^{\frac{3}{2p}-\frac{|\alpha|}{2}}(1+t)^{-(1-\frac{1}{p})-\frac{2}{3}-\lfloor n \rfloor_1 } \]
	which also holds for all $t$, hence these behave correctly both as $t\to 0$ and as $t\to \infty$ as well.  For $p = \infty$, we can choose $p_1 = 8$, $p_3 = 8/3$ and we obtain
	\[ I_2 \leq C\|(a,\vec{\omega})\|_{X_{n,k}}^2 t^{-\frac{|\alpha|}{2}}(1+t)^{-\frac{7}{6}-\lfloor n \rfloor_1 } \]
	separately for $t > 1$ and $t < 1$ and hence $L^{\infty}$ norm has the correct behavior for $t<1$ and $t > 1$ if $|\alpha| = k$.  We can then match the linear decay rate for $p=\infty$, $|\alpha|<k$ by keeping the derivatives on the nonlinearity and using $\beta = \alpha + e_j$ in (\ref{UnweightedDerivNonlinearTermBound}) with $p_1 = p_3 = 12/5$.
	
	As above, we can bound the weighted norms in terms of the weighted norms of each of the components of the convolution.  For the term in which the weight falls on the heat-wave operator we can repeat the estimates on the unweighted norms of the nonlinearity above.  For the other term, we split the integral into two pieces: 
	\[ \begin{split}
	\int\nolimits_0^t (t-& s)^{-\frac{3}{2}(\frac{1}{q_1}-\frac{1}{p}) +1}(1+t-s)^{-(\frac{1}{q_1}-\frac{1}{p})} \Big \|\partial_x^{\alpha} K_{\nu}(\frac{t-s}{2})\ast  \big [ \nabla \cdot N (a(s),\vec{\omega}(s) ) \big ] \Big \|_{\mathring{L}^{\tilde{q}_1}(n)}ds \\ & \leq \Bigg ( \int_0^{t/2} + \int_{t/2}^t \Bigg ) (t-s)^{-\frac{3}{2}(\frac{1}{q_1}-\frac{1}{p}) +\frac{1}{2}}(1+t-s)^{-(\frac{1}{q_1}-\frac{1}{p})} \Big \|\partial_x^{\alpha} K_{\nu}(\frac{t-s}{2})\ast  \big [ \nabla \cdot N (a(s),\vec{\omega}(s) ) \big ] \Big \|_{\mathring{L}^{\tilde{q}_1}(n)}ds \\ & =: I_1 + I_2
	\end{split}  \]
	We can then make use of 
	(\ref{WeightedNonlinearTermBound}) in each to bound the nonlinear term.  For $I_1$ we as usual pull the divergence off of the nonlinearity, and for $0 < n < 1 $ we use (\ref{Constraints2}) to choose $p_1 = p_3 = 3/2$, whereas for $1 < n < 2$ we use (\ref{Constraints3}) to choose $p_1 = p_3 = 3/2$ and we find 
	\[ I_1 \leq C \|(a,\vec{\omega})\|_{X_{n,k}}^2 t^{-\frac{3}{2}(1-\frac{1}{p})+\frac{1}{2}-\frac{|\alpha|}{2}}(1+t)^{-(1-\frac{1}{p})+\frac{1}{3}-\lfloor n \rfloor_1+n} \]
	which holds for $0 < t < \infty$, $1\leq p \leq \infty$.  Then we use (\ref{Constraints3}) to choose $p_1 = 2$ and $p_3 = 6/5$ for $1\leq n \leq 3/2$ and $p_1 = 6/5$ and $p_3 = 2$ for $3/2 < n \leq 2$ and we obtain
	\[ I_1 \leq C \|(a,\vec{\omega})\|_{X_{n,k}}^2 t^{-\frac{3}{2}(1-\frac{1}{p})+\frac{1}{2}-\frac{|\alpha|}{2}}(1+t)^{-(1-\frac{1}{p})-\frac{7}{12}+n} \]
	for $0 < t < \infty$, $1\leq p \leq \infty$.  Similarly for $I_2$ we use (\ref{Constraints2}) to choose $p_1 = p_3 = 3/2$ for $0 < n < 1$ and we use (\ref{Constraints3}) to choose $p_1 = p_3 = 3/2$ for $1 < n < 2$ and we obtain 
	\[ I_2 \leq C \|(a,\vec{\omega})\|_{X_{n,k}}^2 t^{-\frac{3}{2}(1-\frac{1}{p})+\frac{3}{2}-\frac{|\alpha|}{2}}(1+t)^{-(1-\frac{1}{p})-\frac{2}{3}-\lfloor n \rfloor_1+n} \]
	for $1\leq p \leq 2$ and $0 < t <\infty$. Next we use (\ref{Constraints3}) to choose $p_1 = 2$ and $p_3 = 6/5$ for $1\leq n \leq 3/2$ and $p_1 = 6/5$ and $p_3 = 2$ for $3/2 < n \leq 2$ and we find
	\[ I_2 \leq C \|(a,\vec{\omega})\|_{X_{n,k}}^2 t^{-\frac{3}{2}(1-\frac{1}{p})+\frac{5}{4}-\frac{|\alpha|}{2}}(1+t)^{-(1-\frac{1}{p})-\frac{17}{12}+n} \]
	which holds for $0 < t < \infty$ and $1 \leq p < \infty$. For $p =\infty$ we can set $q_1 = 6$ by choosing $p_1 = 8$ and $p_3 = 8/3$ for $0 < n < 1$ using $(\ref{Constraints2})$, choosing $p_1 = 8/3$ and $p_3 = 8$ for $1 \leq n < 2$ using $(\ref{Constraints2})$ and $p_1 = 8$ and $p_3 = 8/3$ for $15/8 < n \leq 2$ using $(\ref{Constraints3})$ to obtain
	\[ I_2 \leq C \|(a,\vec{\omega})\|_{X_{n,k}}^2 t^{-\frac{|\alpha|}{2}}(1+t)^{-\frac{7}{6}-\lfloor n \rfloor_1 + n} \]
	We can then match the linear decay rate for $p =\infty$, $|\alpha|<k$ by keeping the derivatives on the nonlinearity and using $\beta = \alpha + e_j$ in (\ref{WeightedNonlinearTermBound}) and choosing $p_1 = p_3 = 12/5$ for $0 < n < 1$ using $(\ref{Constraints2})$, choosing $p_1 = 2$ and $p_3 = 3$ for $1 \leq n < 2$ using $(\ref{Constraints2})$ and $p_1 = 24/11$ and $p_3 = 8/3$ for $15/8 < n \leq 2$ using $(\ref{Constraints3})$.  Continuity for $t > 0$ is proven as before.
\end{proof}

\section{Asymptotic approximations to the modified compressible Navier-Stokes}

\label{sect:AsympAnalysis}

With these solutions in hand, we turn to the task of approximating these solutions efficiently and accurately, especially in the regime $t \to \infty$.  If $u(t) = (\rho(t),a(t),\vec{\omega}(t))^T$ is the solution belonging to $Y_{n,k} \times X_{n,k}$ given by Theorem $\ref{t:Existence}$ with initial condition $(\rho_0,a_0,\vec{\omega}_0)^T$, $a_0,\vec{\omega}_0$ with zero total mass, then we can write
\begin{equation} \label{LinearNonLinearDecomp}
u(t) = u_L(t) + u_N(t) \end{equation}
where $u_L(t)$ is the linear evolution defined in (\ref{RhoALinearEvolution}), (\ref{OmegaLinearEvolution}), and $u_N(t) = u(t)-u_L(t)$.  We saw in Prop. \ref{p:HypParaHermiteExp}, \ref{p:VecHermiteExp} that for initial conditions $u_0$ belonging to $L^1(n)$ spaces, we can write
\begin{equation} \label{HermiteExpansion} u_L(t) = u_{H,n}(t) + u_{LR,n}(t) \end{equation}
where the Hermite profiles $u_{H,n}(t)$ are defined as in the linear case by
\begin{equation} \label{HermiteProfiles}
\begin{pmatrix}
\rho_{H,n}(x,t) \\ a_{H,n}(x,t)
\end{pmatrix} =  \sum_{\substack{i\leq 2 \\ |\alpha|\leq \lfloor n \rfloor}} \big \langle H_{\alpha} \vec{e}_i, \begin{pmatrix}
\rho_0 \\ a_0 
\end{pmatrix} \big \rangle \text{ } \partial_{x}^{\alpha} \begin{pmatrix}
\rho_i(x,t) \\ a_i(x,t)
\end{pmatrix} \hspace{0.1 cm} \text{ , } \hspace{0.1 cm} 
\vec{\omega}_{H,n}(x,t) = \sum_{\substack{j\leq 2 \\ |\tilde{\alpha}|\leq \lfloor n \rfloor +1}} \big \langle  \vec{p}_{\tilde{\alpha},j}, \vec{\omega}_0 \big \rangle \mathbb{K}_{\epsilon}(t)\ast \vec{f}_{\tilde{\alpha},j}(x)
\end{equation}
where $\rho_i$, $a_i$ are defined in $(\ref{RhoAHermiteProfiles})$, and where $\vec{p}_{\tilde{\alpha},j}, \vec{f}_{\tilde{\alpha},j}$ are defined in Table \ref{tab:DivFreeHermiteExp}.  For simplicity, we will drop the $n$ dependence in the subscript.  We obtained the temporal behavior of $u_{LR}(t)$ in Prop \ref{p:HypParaHermiteExp} \ref{p:VecHermiteExp}.  In the above existence analysis, we saw that $u_N(t)$ decays faster than $u_L(t)$ in some, but not necessarily all, $L^p$ norms, hence we need to study $u_N(t)$ more closely.  We note that $u_N(t)$ can be written as
\[ u_N(t) = - \int \nolimits_0^t e^{\mathcal{L}(t-s)} \mathcal{Q} \big ( u(s),u(s) \big ) ds \]
so inspired by (\ref{LinearNonLinearDecomp}), (\ref{HermiteExpansion}), we define the Hermite-Picard profiles $u_{HP}(t)$ and nonlinear remainder $u_{NR}(t)$:
\begin{equation} \label{NonlinearDecomposition} \begin{split} u_{HP}(t) & := - \int \nolimits_0^t e^{\mathcal{L}(t-s)} \mathcal{Q} \big ( u_{H}(s),u_{H}(s) \big ) ds \\ u_{NR}(t) & := u_{N}(t) - u_{HP}(t) \end{split} \end{equation}
where $u_I(t) = (\rho_I(t),a_I(t),\vec{\omega}_I(t))^T$, $I=L,HP,NR$.  We have already obtained upper bounds on the temporal behavior of $u_{H}(t)$ in Appendix \ref{app:LinEvol} and $u_{LR}(t)$ in Prop \ref{p:HypParaHermiteExp} and \ref{p:VecHermiteExp}.  In what follows, we will obtain upper bounds for $u_{HP}(t)$ and $u_{NR}(t)$, as well as lower bounds for $u_{H}(t)$.  Our goal is to emphasize the role that the localization of the initial conditions (and consequently, the localization of the solutions) plays in determining the nature of the asymptotics.  In particular, we will establish:
\begin{itemize}
	\item For all $n>0$, $u_H(t)$ captures the leading order behavior for $\rho(t)$, but for $n\leq1$ we need to take $u_{app}(t) = u_L(t)$ to capture the leading order behavior for $a(t)$ and $\vec{\omega}(t)$.
	\item For $n > 1$ we need only evaluate the explicit functions $u_H(t)$ to obtain the leading order behavior.
	\item  For $1 < n < 2$ the next order of behavior is given by $u_{LR}(t)$, and $u_{HP}(t)$ and $u_{NR}(t)$ decay faster still.  Hence we could either use $u_{app}(t) =u_L(t)$ or $u_{app}(t) =u_H(t)$.
	\begin{itemize}
		\item In the first case, the error decays $1/2$ power faster than $u_H(t)$, hence is more accurate, but we need to compute a convolution (this is the result in Theorem \ref{t:Summary}). 
		\item In the second case the error decays $(n-1)/2$ powers faster than $u_H(t)$ but we can explicitly evaluate the approximation. 
	\end{itemize}
	\item  Finally, for $n = 2$ there is no loss in accuracy by taking $u_{app}(t) = u_H(t)$. 
\end{itemize}

\subsection{Temporal behavior of the Hermite and Hermite-Picard profiles}

We can use the substitution $\tilde{x} = \frac{x}{\sqrt{1+\epsilon t}}$ together with the explicit form of the Hermite profiles $\vec{\omega}_H(t)$ in Table \ref{tab:DivFreeHermiteExp} and the explicit form of $B\vec{\omega}_H(t)$ to show that their temporal behavior is given by
\[ \begin{split} \|\partial_x^{\alpha} \mathbb{K}_{\epsilon}(t)\ast \vec{f}_{\tilde{\alpha},j} \|_{\mathring{L}^p(\mu)} & = C_{\alpha} (1+t)^{-\frac{3}{2}(1-\frac{1}{p})+\frac{1-|\tilde{\alpha}|-|\alpha|}{2} +\frac{\mu}{2}} \\ \|\partial_x^{\alpha} B \mathbb{K}_{\epsilon}(t)\ast \vec{f}_{\tilde{\alpha},j} \|_{\mathring{L}^p(\mu)} & = \tilde{C}_{\alpha} (1+t)^{-\frac{3}{2}(1-\frac{1}{p})+\frac{2-|\tilde{\alpha}|-|\alpha|}{2} +\frac{\mu}{2}} \end{split}  \]
The temporal behavior of the Hermite profiles $\rho_H(t), a_H(t)$ are given in the following proposition.  These results follow from explicit calculations of the norms involved, as well as the fact that $\Pi$ commutes with the heat-wave operator, and we leave the proof to the reader.  Note that while these estimates might also hold for higher derivatives, we only require derivatives up to the order shown.

\begin{Proposition}
	\label{p:RhoAProfTempBehav}
	There exist functions $C_{l,\alpha}(t)$, $l=1,2,3$ and constants $m,M \in \mathbb{R}$ such that $0 < m < C_{l,\alpha}(t) < M < \infty$ for all $t > 0$ such that
	\[ \begin{split}
	\| \partial_t^l w(t)\ast K_{\nu}(t) \ast \partial_x^{\alpha} \phi_0 \|_{\mathring{L}^p(\mu)} & = C_{l,\alpha}(t) (1+t)^{-\frac{5}{2}(1-\frac{1}{p})+1-\frac{l+|\alpha|}{2}+\mu}
	\end{split}  \]
	for $|\alpha| \leq 2$, $l=0,1,2$, $\mu \in \mathbb{R}_{\geq 0}$ and $1 \leq p \leq \infty$.  Furthermore we have
	\[ \begin{split}
	\|\Pi \partial_t^lw(t)\ast K_{\nu}(t) \ast \partial_x^{\alpha} \phi_0 \|_{\mathring{\mathbb{L}}^p(\mu)} & \leq C (1+t)^{-\frac{5}{2}(1-\frac{1}{p})+\frac{3-l-|\alpha|}{2}+\mu}
	\end{split}  \]
	for any $\alpha \in \mathbb{N}^3 $, $l=1,2$, $\mu \in \mathbb{R}_{\geq 0}$ and $1 \leq p \leq \infty$, except the case when $(\alpha,l) = (0,1)$ and $1\leq p \leq \frac{3}{2+\mu}$. 
\end{Proposition}

This implies that the linear Hermite profiles have temporal behavior given by
\begin{gather} \label{LinearHermiteProfTempBehav}
\|\partial_x^{\alpha}\rho_H(t) \|_{\mathring{L}^p(\mu)} = \tilde{C}_{1,\alpha}(t)E_n(1+t)^{-\frac{5}{2}(1-\frac{1}{p})+\frac{1-|\alpha|}{2}+\mu} \\ \| \partial_x^{\alpha} a_H(t) \|_{\mathring{L}^p(\mu)}  = \tilde{C}_{2,\alpha}(t)E_n (1+t)^{-\frac{5}{2}(1-\frac{1}{p})-\frac{|\alpha|}{2}+\mu} \hspace{.1 cm} \text{ , } \hspace{.1 cm} \| \partial_x^{\alpha} \Pi a_H(t) \|_{\mathring{L}^p(\mu)} = \tilde{C}_{2,\alpha}(t)E_n (1+t)^{-\frac{5}{2}(1-\frac{1}{p})+\frac{1-|\alpha|}{2}+\mu} \notag \\
\|\partial_x^{\alpha} \vec{\omega}_H(t) \|_{\mathring{\mathbb{L}}^p(\mu)} = \tilde{C}_{3,\alpha}(t)E_n(1+t)^{-\frac{3}{2}(1-\frac{1}{p})-\frac{1+|\alpha|}{2}+\frac{\mu}{2}} \hspace{.1 cm} \text{ , } \hspace{.1 cm} \|\partial_x^{\alpha} B \vec{\omega}_H(t) \|_{\mathring{\mathbb{L}}^p(\mu)} = \tilde{C}_{3,\alpha}(t)E_n(1+t)^{-\frac{3}{2}(1-\frac{1}{p})-\frac{|\alpha|}{2}+\frac{\mu}{2}} \notag 
\end{gather}
where $E_n$ is as in (\ref{InitialNorm}), $|\alpha|\leq 1$ and $\tilde{C}_{l,\alpha}(t)$, $\hat{C}_{l,\alpha}(t)$, $l=1,2,3$ are independent of $(\rho_0,a_0,\vec{\omega}_0)^T$ and are such that there exist constants $m,M \in \mathbb{R}$ such that $0 < m < \tilde{C}_{l,\alpha}(t), \hat{C}_{l,\alpha}(t) < M < \infty$ for all $t > 0$.  We also have the following bounds on the Hermite-Picard profiles:

\begin{Proposition}
	\label{p:NonlinearHermiteProfTempBehav}
	There exists a constant $C$ such that we have
	\[ \begin{split}
	\|\rho_{HP}(t)\|_{\mathring{L}^p(\mu)} \leq C(1+t)^{-\frac{5}{2}(1-\frac{1}{p})+\frac{1}{2}+\mu-\frac{1}{2}} \hspace{.75 cm} & \text{ , } \hspace{.75 cm}  \| a_{HP}(t) \|_{\mathring{L}^p(\mu)} \leq C (1+t)^{-\frac{5}{2}(1-\frac{1}{p})+\mu-\frac{1}{2}} \\
	\|\vec{\omega}_{HP}(t)\|_{\mathring{\mathbb{L}}^p(\mu)} & \leq C(1+t)^{-\frac{3}{2}(1-\frac{1}{p})-\frac{1}{2}+\mu-\frac{1}{2}}
	\end{split}  \]
	for all $t>0$, $|\alpha| \leq 2$, $0\leq \mu \leq 2$ and $1 \leq p \leq \infty$.
\end{Proposition}

\begin{proof}	
	We start with the Hermite-Picard profile $a_{HP}(t)$.  We look at the weighted norms for an arbitrary weight $\mu$.  We first split the convolution:
	\[ \begin{split}
	\int \nolimits_0^t \Big \| \partial_t & w(t-s) \ast K_{\nu}(t-s) \ast \big [ \nabla \cdot N (a_H(s),\vec{\omega}_H(s) ) \big ] \Big \|_{\mathring{L}^p(\mu)}ds \\  & \leq \int \nolimits_0^t \Big \| \partial_t w(t-s) \ast K_{\nu}(\frac{t-s}{2}) \Big \|_{\mathring{L}^{\tilde{q}}(\mu)} \Big \| K_{\nu}(\frac{t-s}{2})\ast  \big [ \nabla \cdot N (a_H(s),\vec{\omega}_H(s) ) \big ] \Big \|_{L^{\tilde{q}_1}}ds \\ & \hspace{1 cm} + \int \nolimits_0^t \Big \| \partial_t w(t-s) \ast K_{\nu}(\frac{t-s}{2}) \Big \|_{L^q} \Big \| K_{\nu}(\frac{t-s}{2})\ast \big [ \nabla \cdot N (a_H(s),\vec{\omega}_H(s) ) \big ] \Big \|_{\mathring{L}^{q_1}(\mu)}ds \end{split}
	\]
	We'll bound the second term, and then as in the existence proof the bounds on the first term follow by repeating the estimates for the second term line by line after using the weighted estimate on the heat-wave operator in Prop \ref{p:HypParaEst} and taking $\mu = 0$ on the nonlinear term.  We first split the second integral into two: 
	\[ \begin{split}
	I_1 + I_2 := \Bigg ( \int_{0}^{t/2} + \int_{t/2}^t \Bigg ) \Big \| \partial_t w(t-s) \ast K_{\nu}(\frac{t-s}{2}) \Big \|_{L^q} \Big \| K_{\nu}(\frac{t-s}{2}) \ast \Big [ \nabla \cdot N(a_H(s),\vec{\omega}_H(s)) \Big ] \Big \|_{\mathring{L}^{q_1}(\mu)} ds \end{split} \]
	For $t<1$ we can choose $q=1$ in both terms, and since our heat estimate and equation (\ref{LinearHermiteProfTempBehav}) can be used to show the resulting integrand is bounded, these remain bounded as $t\to 0$.  Hence we need only consider $t >1$.  For $I_1$ we can use the heat estimate to remove both of the derivatives from the nonlinearity, set $q_1 = 1$, use Cauchy-Schwarz and make use of (\ref{LinearHermiteProfTempBehav}) to bound the norms of $\vec{m}_H$ via
	\[ \begin{split} I_1 & \leq C \int \nolimits_{0}^{\frac{t}{2}} (t-s)^{-\frac{3}{2}(1-\frac{1}{p})}(1+t-s)^{\frac{1}{2}-(1-\frac{1}{p})} \Big \| K_{\nu}(\frac{t-s}{2}) \ast \Big [ \nabla \cdot N(a_H(s),\vec{\omega}_H(s)) \Big ] \Big \|_{\mathring{L}^1(\mu)} ds \\ & \leq C \max_{ij} \int \nolimits_{0}^{\frac{t}{2}} (t-s)^{-\frac{3}{2}(1-\frac{1}{p})-1}(1+t-s)^{\frac{1}{2}-(1-\frac{1}{p})} \| m_{H,i}(s) \|_{L^2} \|m_{H,j}(s) \|_{L^2(\mu)} ds \\ & \leq C E_n^2 (1+t)^{-\frac{5}{2}(1-\frac{1}{p})-\frac{1}{2}+\mu}  \end{split} \]
	For $I_2$ we use the heat estimate but keep all of the derivatives on the nonlinearity and we obtain
	\[ \begin{split} I_2 & \leq C \max_{ij} \int \nolimits_{t/2}^t (t-s)^{-\frac{3}{2}(\frac{1}{q_1}-\frac{1}{p})}(1+t-s)^{\frac{1}{2}-(\frac{1}{q_1}-\frac{1}{p})} \big \| \partial_{x_i}\partial_{x_j} \big ( m_{H,i}(s)m_{H,j}(s) \big ) \big \|_{L^{q_1}(\mu)} ds \\ \leq C & \max_{ijk\ell} \int_{t/2}^t (t-s)^{-\frac{3}{2}(\frac{1}{q_1}-\frac{1}{p})}(1+t-s)^{\frac{1}{2}-(\frac{1}{q_1}-\frac{1}{p})} \Big ( \big \| \partial_{x_i}\partial_{x_j} m_{H,i} \big \|_{L^{p_1}} \big \| m_{H,j} \big \|_{L^{p_2}(\mu)} + \big \| \partial_{x_i} m_{H,j} \big \|_{L^{p_1}} \big \| \partial_{x_k} m_{H,\ell} \big \|_{L^{p_2}(\mu)} \Big ) ds \\ = & J_1 + J_2 \end{split} \]
	For $J_1$ we use Cor \ref{c:PiAndBiotSavartOperators} part \textbf{(a)} to obtain
	\[ \big \| \partial_{x_i}\partial_{x_j} m_{H,i} \big \|_{L^{p_1}} \big \| m_{H,j} \big \|_{L^{p_2}(\mu)} \leq C \max_j \big ( \big \| \partial_{x_j} a_H \big \|_{L^{p_1}} + \big \| \partial_{x_j} \vec{\omega}_H \big \|_{\mathbb{L}^{p_1}} \big ) \big \| m_{H,j} \big \|_{L^{p_2}(\mu)}  \]
	so for $1\leq p \leq 2$ we can set $q_1 =1 $ by choosing $p_1 = p_2 = 2$ and use (\ref{LinearHermiteProfTempBehav}) to obtain
	\[ \begin{split}
	J_1 \leq C E_n^2 t^{-\frac{3}{2}(1-\frac{1}{p})+1}(1+t)^{\frac{1}{2}-(1-\frac{1}{p})-\frac{5}{2} + \mu}
	\end{split} \]
	whereas for $p = \infty$ we can let $q_1 = 3/(2-\delta)$ by setting $p_1=p_2 = 6/(2-\delta)$, where $0 < \delta < 1/5$ is any number and we obtain the following 
	\[ \begin{split}
	J_1 \leq C E_n^2(1+t)^{-\frac{1}{2}(4+\delta)-1+\mu}
	\end{split} \]
	For $J_2$ we can just use (\ref{LinearHermiteProfTempBehav}) directly, and by choosing $p_1 = p_2 = 2$ for $1 \leq p \leq 2$ we obtain
	\[ \begin{split}
	J_2 \leq C \big \| (a,\vec{\omega})\big \|_{X_n}^2 t^{-\frac{3}{2}(1-\frac{1}{p})+1}(1+t)^{-2-(1-\frac{1}{p})+\mu}
	\end{split} \]
	and we can obtain the analogous results for $p = \infty$ by choosing $q_1 = 3/(2-\delta)$ by setting $p_1 = p_2 = 6/(2-\delta)$ for some $0 < \delta < 1/5$.  
	
	The bounds for the Hermite-Picard profiles $\rho_{HP}(t)$ and $\vec{\omega}_{HP}(t)$ can be obtained by similar arguments.
\end{proof}

\subsection{Temporal behavior of the linear and nonlinear remainders}

If one naively uses the estimates in Cor. \ref{c:PiAndBiotSavartOperators} to obtain the decay rate of $\Pi a$, then one finds that the $L^p$ norms of $\Pi a(t)$ grows by a factor of $t^{5/6}$ relative to $a(t)$, whereas one finds that $B\vec{\omega}(t)$ grows by a factor of $t^{1/2}$ relative to $\vec{\omega}(t)$.  We saw in Prop \ref{p:RhoAProfTempBehav} that $\Pi a_{H}(t)$ grows by a factor of $t^{1/2}$ relative to $a_H(t)$, and we now prove that the same holds for remainder $a_{LR}$:

\begin{Proposition}\label{p:PiDecay}
	 Let $n \in [ 0,2 ]$ and let $(\rho_0,a_0)$ belong to $W^{1,p}(n)\times L^{p}(n)$ for all $1 \leq p \leq \frac{3}{2}$.  Then for $a_{LR}(t)$ defined as in Prop. \ref{p:HypParaHermiteExp}, we have
	\[ \| \partial_x^{\alpha} \Pi a_{LR} (t) \|_{\mathring{L}^p(\mu)} \leq CE_n t^{-\frac{5}{2}(1-\frac{1}{p})+\frac{1}{2}+\mu-\frac{n}{2}-\frac{|\alpha|}{2}} \]
	for $t > 1$, $0\leq \mu\leq n$ and any nonzero $\alpha \in \mathbb{N}^3 $, $1 \leq p \leq \infty$.  If $\alpha = 0$ the above estimate holds for $p > \frac{3}{2+\mu}$.  On the other hand for $t < 1$, $0 \leq \mu \leq n$ and for $3/(2-\mu) < p < \infty$ if $n<1$, or $\max(3/2,3/(3-\mu)) < p < \infty$ if $n\geq1$ we have
	\[ \| \partial_x^{\alpha} \Pi a_{LR} (t) \|_{\mathring{L}^p(\mu)} \leq C E_n t^{-r_{\alpha,\tilde{p}}} \]
	  where $p^{-1} = \tilde{p}^{-1}-3^{-1}$.
\end{Proposition}

\begin{proof}
	The estimate for $t <1$ follows from Cor. \ref{c:PiAndBiotSavartOperators} parts \textbf{(c)}, \textbf{(c)}, and from interpolation in the case when $n \geq 1$ and $\mu <1$.  For $t > 1$ the interesting case is when $|\alpha|>0$, and we have
	\[ \| \partial_x^{\alpha} \Pi a_{LR}(t) \|_{\mathring{L}^p(\mu)} \leq \| \Pi \partial_t w(t)\ast \partial_x^{\alpha} K_{\nu}(t)\ast a_{LR}(0) \|_{\mathring{L}^p(\mu)} + \| \Pi \partial_t^2 w(t)\ast \partial_x^{\alpha} K_{\nu}(t)\ast \rho_{LR}(0) \|_{\mathring{L}^p(\mu)} \]
	so if $\alpha = e_i + \beta$ for some $i,\beta$, we can use Young's inequality to obtain
	\[ \begin{split} \| \Pi & \partial_t w(t)\ast \partial_x^{\alpha} K_{\nu}(t)\ast a_{LR}(0) \|_{\mathring{L}^p(\mu)} = \| \pi \ast \partial_t w(t)\ast \partial_{x_i} K_{\nu}(\frac{t}{2}) \ast \partial_x^{\beta} K_{\nu}(\frac{t}{2}) \ast a_{LR}(0) \|_{\mathring{L}^p(\mu)} \\ & \leq \| \Pi \partial_t w(t)\ast \partial_{x_i} K_{\nu}(\frac{t}{2}) \|_{\mathring{L}^p(\mu)} \| \partial_x^{\beta}K_{\nu}(\frac{t}{2}) \ast a_{LR}(0) \|_{L^1} + \| \Pi \partial_t w(t)\ast \partial_{x_i} K_{\nu}(\frac{t}{2}) \|_{L^p} \| \partial_x^{\beta} K_{\nu}(\frac{t}{2}) \ast a_{LR}(0) \|_{\mathring{L}^1(\mu)} \end{split} \]
	where
	\[ \pi(x) = -\frac{1}{4\pi} \frac{x}{|x|^3} \]
	is the integral kernel of the $\Pi$ operator.  The result then follows from our estimates of the $\Pi$ operator acting on the Hermite term in Prop \ref{p:RhoAProfTempBehav}, since the same result applies to the heat-wave operator.  However, for $\alpha = 0$ the heat-wave operator only belongs to $L^p$ for $p > 3/2$.  We leave the remainder of the proof to the reader.
\end{proof}

In the following lemma, we collect the bounds for $u_N(t)$ obtained during the contraction mapping argument in the existence proof and sharpen one of them.  For this purpose we define the rate $b_{n,p}$ to measure the excess decay of $u_N(t)$ above the linear rate as follows, using interpolation for $2 < p < \infty$:
\begin{equation}
\label{ImprovedNonlinearDecayRate}
b_{n,p} = \left \{ \begin{array}{cl} \min(\frac{1}{6} + \frac{\lfloor n \rfloor_1}{2},\frac{3}{10}+\frac{\lfloor n \rfloor_1}{10}) & \text{ for } 1 \leq p \leq 2  \\ \min(\frac{\lfloor n \rfloor_1}{2},\frac{3}{10}+\frac{\lfloor n \rfloor_1}{10})  & \text{ for } p = \infty \\ (b_{n,\infty} - b_{n,2}) \big (1-\frac{2}{p}) + b_{n,2} & \text{ for } 2 < p < \infty\end{array} \right .
\end{equation}

\begin{Lemma}
	\label{lem:NonlinDecay}
	Let $n \in [0,2]$, $k \geq 1$ and let $u_0 = (\rho_0,a_0,\vec{\omega}_0)^T \in \cap_{1\leq p \leq \frac{3}{2}} W^{1,p}(n) \times L^p(n) \times \mathbb{L}_{\sigma}^{p}(n)$.   If $u(t) = (\rho(t),a(t),\vec{\omega}(t))^T$ is the solution in $Y_{n,k}\times X_{n,k}$ given by Theorem \ref{t:Existence} and Corollary \ref{c:ExistencePt2} with initial condition $u_0$, then the nonlinear term $u_N(t)$ in (\ref{LinearNonLinearDecomp}) satisfies
	\[ \begin{split}\|\partial_{x}^{\alpha} \rho_N(t)\|_{\mathring{L}^p(\mu)} & \leq  CE_n^2 t^{-r_{\alpha,p}}(1+t)^{-\ell_{n,p,\mu}+\frac{1}{2}-b_{n,p}} \\ \|\partial_{x}^{\alpha} a_N(t)\|_{\mathring{L}^p(\mu)} & \leq C E_n^2 t^{-r_{\alpha,p}}(1+t)^{-\ell_{n,p,\mu}-b_{n,p}} \\ \|\partial_{x}^{\alpha} \vec{\omega}_{N}(t)\|_{\mathring{\mathbb{L}}^p(\mu)} & \leq CE_n^2 t^{-r_{\alpha,p}}(1+t)^{-\tilde{\ell}_{n,p,\mu}-b_{n,p}}
	\end{split} \]
	for $1\leq p \leq \infty$, $0\leq \mu \leq n$ and $|\alpha|<k$.
\end{Lemma}

\begin{proof}
	The estimates for $t < 1$ are the same as those obtained in the existence proof, hence we need only consider $t>1$.  By inspecting the estimates in the existence proof, we see that all of the bounds obtained already exhibit the extra decay listed in the first argument of the minimum in (\ref{ImprovedNonlinearDecayRate}), with one important exception.  The estimate of the unweighted norm of $I_1$ in (\ref{I1Est}) stops improving relative to the linear rate for $n > 1/3$.  The $|\alpha|=k$ derivative also may decay slower, but we don't estimate this here.
	
	So we need only improve on the bound in (\ref{I1Est}) for $n >1/3$.  We can split $I_1$ into two pieces:
	\[ \begin{split} I_1 = \Bigg ( \int_0^{t^{3/5}} + \int_{t^{3/5}}^{t/2} \Bigg ) (t-s)^{-\frac{3}{2}(\frac{1}{q_1}-\frac{1}{p}) }(1+t-s)^{\frac{1}{2}-(\frac{1}{q_1}-\frac{1}{p})} \Big \| \partial_x^{\alpha} K_{\nu}(\frac{t-s}{2}) \ast \big [ \nabla \cdot N (a(s),\vec{\omega}(s) ) \big ] \Big \|_{L^{q_1}}ds =: J_1 + J_2 \end{split} \]
	Since we are interested in the limit $t\to \infty$ we assume $t/2 > t^{3/5}$ here, but for $1 < t^{2/5} \leq 2$ we can obtain the analogous result.  For $J_1$ we make a modified estimate by taking all of the derivatives off of the nonlinearity and onto the heat-wave propagator by using our heat estimate.  We can then set $q_1 = 1$, use Cauchy-Schwarz and use Cor \ref{c:PiAndBiotSavartOperators} part \textbf{(b)} to obtain 
	\[ \begin{split} J_1 & \leq \int \nolimits_0^{t^{3/5}} (t-s)^{-\frac{3}{2}(\frac{1}{q_1}-\frac{1}{p})-1-\frac{|\alpha|}{2} }(1+t-s)^{\frac{1}{2}-(\frac{1}{q_1}-\frac{1}{p})} \| m_{i}(s)m_{j}(s) \|_{L^{q_1}} ds \\ & \leq C E_n^2 t^{-\frac{5}{2}(1-\frac{1}{p})-\frac{1}{2}-\frac{|\alpha|}{2} } \int \nolimits_0^{t^{3/5}} \big ( \|a(s)\|_{L^{6/5}} +\|\vec{\omega}(s)\|_{\mathbb{L}^{6/5}} \big )^2 ds \leq CE_n^2 t^{-\frac{5}{2}(1-\frac{1}{p})-\frac{1}{2}-\frac{|\alpha|}{2}+\frac{7}{10}-\frac{3\lfloor n \rfloor_1}{5}} \end{split} \]
	For $J_2$, we can use the same estimate as before.  Taking the divergence and $\partial_x^{\alpha}$ off of the nonlinearity by using our heat estimate, setting $q_1 = 1$, using H\"older's inequality and Cor \ref{c:PiAndBiotSavartOperators} parts \textbf{(a)} and \textbf{(c)} we obtain the following for $n > 1/3$:
	\[ \begin{split} J_2 & \leq \int \nolimits_{t^{3/5}}^{t/2} (t-s)^{-\frac{3}{2}(\frac{1}{q_1}-\frac{1}{p})-\frac{1}{2}-\frac{|\alpha|}{2} }(1+t-s)^{\frac{1}{2}-(\frac{1}{q_1}-\frac{1}{p})} \| \partial_{x_i}(m_{i})m_{k} \|_{L^{q_1}} ds \\ & \leq  CE_n^2 t^{-\frac{5}{2}(1-\frac{1}{p})-\frac{|\alpha|}{2} } \int \nolimits_{t^{3/5}}^{t/2} (1+s)^{-\frac{2}{3}-\lfloor n \rfloor_1} ds \leq C E_n^2 t^{-\frac{5}{2}(1-\frac{1}{p})-\frac{|\alpha|}{2}+\frac{1}{5}-\frac{3\lfloor n \rfloor_1}{5}} \end{split} \]
	This same improved bound can be obtained for $\rho_N(t)$ and $\vec{\omega}_N(t)$ as well.
\end{proof}



We now use the estimates just proven, together with a bootstrapping argument, to obtain more refined estimates of the temporal decay of the nonlinear remainder.   For this purpose we define the rate $\mathpzc{b}_{n,p}$ to measure the excess decay of $u_{NR}(t)$ above the linear rate via
\begin{equation}
\label{NonlinearRemainderDecayRate}
\begin{split}
\mathfrak{b}_{n,p} & = \left \{ \begin{array}{cl} \frac{\lfloor n \rfloor_1 - 1}{2} + \min(2n-\frac{1}{3},n,\frac{1}{2}) & \text{ for } 1 \leq p \leq 2  \\ \frac{\lfloor n \rfloor_1 - 1}{2} + \min(n-\frac{1}{2},\frac{1}{2})  & \text{ for } p = \infty \\ (\mathfrak{b}_{n,\infty} - \mathfrak{b}_{n,2}) \big ( 1-\frac{2}{p}) + \mathfrak{b}_{n,2} & \text{ for } 2 < p < \infty\end{array} \right . \\
\mathpzc{b}_{n,p} & = \max (b_{n,p},\mathfrak{b}_{n,p} )
\end{split}
\end{equation}

\begin{Theorem}
	\label{t:Asymp}
	Let $n \in [0,2]$, $k \geq 1$ and let $u_0 = (\rho_0,a_0,\vec{\omega}_0)^T \in \cap_{1\leq p \leq \frac{3}{2}} W^{1,p}(n) \times L^p(n) \times \mathbb{L}_{\sigma}^{p}(n)$.   If $u(t) = (\rho(t),a(t),\vec{\omega}(t))^T$is the solution in $Y_{n,k}\times X_{n,k}$ given by Theorem \ref{t:Existence} and Corollary \ref{c:ExistencePt2} with initial condition $u_0$, then the nonlinear remainder $u_{NR}(t)$ in (\ref{NonlinearDecomposition}) satisfies
	\[ \begin{split}\|\partial_{x}^{\alpha} \rho_{NR}(t)\|_{\mathring{L}^p(\mu)} & \leq  CE_n^2(1+E_n^2) t^{-r_{\alpha,p}}(1+t)^{-\ell_{n,p,\mu}+\frac{1}{2}-\mathpzc{b}_{n,p}} \\ \|\partial_{x}^{\alpha} a_{NR}(t)\|_{\mathring{L}^p(\mu)} & \leq C E_n^2(1+E_n^2) t^{-r_{\alpha,p}}(1+t)^{-\ell_{n,p,\mu}-\mathpzc{b}_{n,p}} \\ \|\partial_{x}^{\alpha} \vec{\omega}_{NR}(t)\|_{\mathring{\mathbb{L}}^p(\mu)} & \leq CE_n^2(1+E_n^2) t^{-r_{\alpha,p}}(1+t)^{-\tilde{\ell}_{n,p,\mu}-\mathpzc{b}_{n,p}}
	\end{split} \]
	for $1\leq p \leq \infty$, $0\leq \mu \leq n$ and $|\alpha|\leq \min(1,k-1)$.
\end{Theorem}

\begin{proof}
	Again the estimates for $t<1$ are identical to those in the existence proof, so we only consider $t>1$.  Here we need only obtain the rate $\mathfrak{b}_{n,p}$, and the result with rate $\mathpzc{b}_{n,p}$ follows from our bounds on $u_N(t)$.  By definition we see that the nonlinear remainder $u_{NR}$ must satisfy the following equation:
	\[ \begin{split}
	u_{NR}(t) = - \int \nolimits_0^t e^{\mathcal{L}(t-s)} \Big [  \mathcal{Q} \big ( u_H,u_{LR}+u_{N} \big ) + \mathcal{Q} \big ( u_{LR}+u_{N},u_{H} \big ) + \mathcal{Q} \big ( u_{LR}+u_N,u_{LR}+u_N \big ) \Big ] ds
	\end{split} \]
	
	We start by looking at the Duhamel term corresponding to $a_{NR}$.  By expanding the nonlinearity, we see that for an arbitrary weight $0 \leq \mu \leq n$ we need to bound the norms of terms of the form
	\[ \begin{split} \int \nolimits_{0}^t \Big \| \partial_{t} w(t-s) & \ast \partial_{x}^{\alpha} K_{\nu}(t-s) \ast \Big [ \partial_{x_i}\partial_{x_j} \big [ m_{I,i}(s) m_{J,j}(s) \big ] \Big ] \Big \|_{\mathring{L}^p(\mu)} ds \\  & \leq \int \nolimits_0^t \Big \| \partial_t w(t-s) \ast K_{\nu}(\frac{t-s}{2}) \Big \|_{\mathring{L}^q(\mu)} \Big \|\partial_{x}^{\alpha}  K_{\nu}(\frac{t-s}{2}) \ast \Big [ \partial_{x_i}\partial_{x_j} \big [m_{I,i}(s) m_{J,j}(s) \big ] \Big ] \Big \|_{L^{q_1}}ds \\ & + \int \nolimits_0^t \Big \| \partial_t w(t-s) \ast K_{\nu}(\frac{t-s}{2}) \Big \|_{L^q} \Big \|\partial_{x}^{\alpha}  K_{\nu}(\frac{t-s}{2}) \ast \Big [ \partial_{x_i}\partial_{x_j} \big [ m_{I,i}(s) m_{J,j}(s) \big ] \Big ] \Big \|_{\mathring{L}^{q_1}(\mu)}ds  \end{split}  \]
	for pairs of indices $(I,J) = (H,LR), (H,N), (LR,LR), (LR,N)$ and $(N,N)$.  We will bound the second term, and the bounds for the first can then be obtained by repeating the same analysis by using the weighted bounds in Prop. \ref{p:HypParaEst} as described previously.  We split the second term into two:
	
	\[ \begin{split}
	\int \nolimits_0^t & \Big \| \partial_t w(t-s) \ast K_{\nu}(\frac{t-s}{2}) \Big \|_{L^q} \Big \| \partial_{x}^{\alpha}  K_{\nu}(\frac{t-s}{2}) \ast \Big [ \partial_{x_i}\partial_{x_j} \big [ m_{I,i}(s) m_{J,j}(s) \big ] \Big ] \Big \|_{\mathring{L}^{q_1}(\mu)}ds \\ & \leq \Bigg ( \int_0^{t/2} + \int_{t/2}^t \Bigg ) \Big \| \partial_t w(t-s) \ast K_{\nu}(\frac{t-s}{2}) \Big \|_{L^q} \Big \|\partial_{x}^{\alpha}  K_{\nu}(\frac{t-s}{2}) \ast \Big [ \partial_{x_i}\partial_{x_j} \big [ m_{I,i}(s) m_{J,j}(s) \big ] \Big ] \Big \|_{\mathring{L}^{q_1}(\mu)}ds = I^{IJ}_1 + I^{IJ}_2
	\end{split} \] 
	Bounds for $I_1^{IJ}$ and $I_2^{IJ}$ can be obtained for $(I,J)= (H,LR), (H,N), (LR,LR),$ and $(LR,N)$ using very similar arguments.  We bound these first, then bound $(N,N)$ later.  For $I_1^{IJ}$ we use the heat estimate to take all of the derivatives off of the nonlinear term, and then use H\"older's inequality as follows:
	\begin{equation} \label{I1AllDerivs} \begin{split}
	I_1^{IJ} \leq \int \nolimits_0^{t/2} (t-s)^{-\frac{3}{2}(\frac{1}{q_1}-\frac{1}{p})-1-\frac{|\alpha|}{2}}(1+t-s)^{\frac{1}{2}-(\frac{1}{q_1}-\frac{1}{p})} \| m_{I,i}(s) \|_{\mathring{L}^{p_1}(\mu)}\| m_{J,j}(s) \|_{L^{p_2}} ds 
	\end{split} \end{equation}
	For $(I,J) = (H,LR), (LR,LR)$ we choose $p_1 = p_2 = 2$ for $1 \leq p \leq \infty$ and use our estimates in (\ref{LinearHermiteProfTempBehav}) and our estimates of the linear remainder in Prop \ref{p:PiDecay}, and we obtain
	\[ \begin{split} I_{1}^{H,LR} \leq CE_n^2t^{-\frac{5}{2}(1-\frac{1}{p})-\frac{1}{2}-\frac{|\alpha|}{2}+\mu} L_{n,0}(t) \hspace{.5 cm} \text{ and } \hspace{.5 cm} I_{1}^{LR,LR} \leq CE_n^2t^{-\frac{5}{2}(1-\frac{1}{p})-\frac{1}{2}-\frac{|\alpha|}{2}+\mu + \max(\frac{1}{2}-n,0)} L_{n,1/2}(t) \end{split}  \] 
	For $(I,J) = (H,N), (LR,N)$ we use Cor. \ref{c:PiAndBiotSavartOperators} \textbf{(b)} and pull the first factors out of the integral to obtain
	\[ I_1^{IJ} \leq   t^{-\frac{3}{2}(\frac{1}{q_1}-\frac{1}{p})-1-\frac{|\alpha|}{2}}(1+t)^{\frac{1}{2}-(\frac{1}{q_1}-\frac{1}{p})}\int \nolimits_0^{t/2} \| m_{I,i}(s) \|_{\mathring{L}^{p_1}(\mu)} \big ( \| a_{J}(s) \|_{L^{p_3}} + \| \vec{\omega}_{J}(s) \|_{\mathbb{L}^{p_3}} \big ) ds  \]
	and we then set $p_1 = p_3 = 3/2$ for $1\leq p \leq \infty$ and use our estimates in (\ref{LinearHermiteProfTempBehav}) and in Prop \ref{p:PiDecay} together with the estimate of the nonlinear term in Lemma \ref{lem:NonlinDecay}.  We find
	\[ \begin{split}
	I_{1}^{H,N} & \leq CE_n^3t^{-\frac{5}{2}(1-\frac{1}{p})-\frac{1}{2}-\frac{|\alpha|}{2}+\mu+\max(\frac{1}{6}-n,0)} L_{n,1/6}(t) \\ I_{1}^{LR,N} & \leq CE_n^3t^{-\frac{5}{2}(1-\frac{1}{p})-\frac{1}{2}-\frac{|\alpha|}{2}+\mu + \max(\frac{5}{6}-n-b_{n,3/2},0)} L_{n,16/33}(t) \end{split}  \] 
	
	For $I_2^{IJ}$ we leave all the derivatives on the nonlinear term and obtain 
	\begin{equation} \label{I2NoDerivs} \begin{split}
	I_2^{IJ} & \leq \int \nolimits_{t/2}^t (t-s)^{-\frac{3}{2}(\frac{1}{q_1}-\frac{1}{p})}(1+t-s)^{\frac{1}{2}-(\frac{1}{q_1}-\frac{1}{p})} \big \|   \partial_{x}^{\alpha} \partial_{x_i} \partial_{x_j} \big [ m_{I,i}(s) m_{J,j}(s) \big ] \big \|_{\mathring{L}^{q_1}(\mu)}ds 
	\end{split} \end{equation}
	Using Liebniz's rule and H\"older's inequality we have
	\begin{equation} \label{Liebniz} \big \| \partial_{x}^{\alpha+e_i+ e_j} \big [m_{I,i}(s) m_{J,j}(s) \big ] \big \|_{\mathring{L}^{q_1}(\mu)} \leq \sum \nolimits_{\gamma_1 + \gamma_2 = \alpha +e_i+ e_j} \| \partial_{x}^{\gamma_1} m_{I,i} \|_{\mathring{L}^{p_1}(\mu)} \| \partial_x^{\gamma_2} m_{J,j} \|_{L^{p_2}} \end{equation}
	For $(I,J)=(H,LR),(LR,LR)$ we choose $p_1=p_2 = 2$ for $1\leq p \leq 2$ and make use of our estimates in Prop \ref{p:PiDecay} and \ref{p:VecHermiteExp}.  Here we obtain 
	\[  I_{2}^{H,LR} \leq CE_n^2t^{-\frac{5}{2}(1-\frac{1}{p})-\frac{1+n}{2}-\frac{|\alpha|}{2}+\mu} \hspace{.5 cm} \text{ and } \hspace{.5 cm} I_{2}^{LR,LR} \leq CE_n^2t^{-\frac{5}{2}(1-\frac{1}{p})-\frac{|\alpha|}{2}+\mu -n}  \]
	whereas for $p= \infty$ we choose $p_1=p_2 = 4$ and use Prop \ref{p:PiDecay} and \ref{p:VecHermiteExp} to obtain 
	\[  I_{2}^{H,LR} \leq CE_n^2t^{-\frac{5}{2}-\frac{n}{2}-\frac{|\alpha|}{2}+\mu} \hspace{.5 cm} \text{ and } \hspace{.5 cm} I_{2}^{LR,LR} \leq CE_n^2t^{-\frac{5}{2}-\frac{|\alpha|}{2}+\mu +\frac{1}{2}-n}  \]
	On the other hand, for $(I,J) = (H,N),(LR,N)$ we use Cor. \ref{c:PiAndBiotSavartOperators} part \textbf{(b)} on $m_{N,j}$ when $\gamma_2 = 0$ and choose $p_1 = p_3 = 3/2$ for $1 \leq p \leq 2$ to obtain
	\[  I_{2}^{H,N} \leq CE_n^3t^{-\frac{5}{2}(1-\frac{1}{p})-\frac{|\alpha|}{2}+\mu-\frac{1}{6}-\frac{\lfloor n\rfloor_1}{2} -b_{n,3/2}} \hspace{.5 cm} \text{ and } \hspace{.5 cm} I_{2}^{LR,N} \leq CE_n^3t^{-\frac{5}{2}(1-\frac{1}{p})-\frac{|\alpha|}{2}+\mu +\frac{1}{3}-\frac{n+\lfloor n \rfloor_1}{2} - b_{n,3/2}}  \]
	whereas we choose $p_1 = 3, p_3 = 2$ for $p = \infty$ to obtain
	\[  I_{2}^{H,N} \leq CE_n^3t^{-\frac{5}{2}-\frac{|\alpha|}{2}+\mu-\frac{\lfloor n\rfloor_1}{2} -b_{n,2}} \hspace{.5 cm} \text{ and } \hspace{.5 cm} I_{2}^{LR,N} \leq CE_n^3t^{-\frac{5}{2}-\frac{|\alpha|}{2}+\mu +\frac{1}{2}-\frac{n+\lfloor n \rfloor_1}{2}-b_{n,2}}  \]
	If $\gamma_2 \neq 0$ then we use Cor. \ref{c:PiAndBiotSavartOperators} part \textbf{(a)} on $m_{N,j}$ and choose $p_1 = p_2 = 2$ for $1 \leq p \leq 2$ to obtain 
	\[  I_{2}^{H,N} \leq CE_n^3t^{-\frac{5}{2}(1-\frac{1}{p})-\frac{|\alpha|}{2}+\mu-\frac{1}{2}-\frac{\lfloor n\rfloor_1}{2} -b_{n,2}} \hspace{.5 cm} \text{ and } \hspace{.5 cm} I_{2}^{LR,N} \leq CE_n^3t^{-\frac{5}{2}(1-\frac{1}{p})-\frac{|\alpha|}{2}+\mu -\frac{n+\lfloor n \rfloor_1}{2} - b_{n,2}}  \]
	whereas we choose $p_1 = \infty,p_2 = 2$ for $p = \infty$ to obtain 
	\[  I_{2}^{H,N} \leq CE_n^3t^{-\frac{5}{2}-\frac{|\alpha|}{2}+\mu-\frac{\lfloor n\rfloor_1}{2} -b_{n,2}} \hspace{.5 cm} \text{ and } \hspace{.5 cm} I_{2}^{LR,N} \leq CE_n^3t^{-\frac{5}{2}-\frac{|\alpha|}{2}+\mu +\frac{1}{2}-\frac{n+\lfloor n \rfloor_1}{2}-b_{n,2}}  \]
	
	We also need to bound the norms of the terms for which $(I,J) = (N,N)$.  For this we will need to bound $\mu = 0$ and $\mu = n$ separately, and the remaining bounds follow from interpolation.  Starting with $\mu = 0$ we first bound $I_{1}^{NN}$ by removing all derivatives from the nonlinearity using the heat estimate and use H\"older's inequality as in (\ref{I1AllDerivs}), but we then use Cor. \ref{c:PiAndBiotSavartOperators} \textbf{(b)} on both terms to obtain
	\[ \begin{split}
	I_1^{NN} \leq \int \nolimits_0^{t/2} (t-s)^{-\frac{3}{2}(\frac{1}{q_1}-\frac{1}{p})-1-\frac{|\alpha|}{2}}(1+t-s)^{\frac{1}{2}-(\frac{1}{q_1}-\frac{1}{p})} \big ( \|a_N\|_{L^{p_3}} + \| \vec{\omega}_N \|_{\mathbb{L}^{p_3}} \big ) \big ( \|a_N\|_{L^{p_4}} + \| \vec{\omega}_N \|_{\mathbb{L}^{p_4}} \big )  ds 
	\end{split} \]
	We can then choose $p_3 = p_4 = 6/5$ for $1\leq p \leq \infty$ to obtain
	\[ I_1^{NN} \leq CE_n^4t^{-\frac{5}{2}(1-\frac{1}{p})-\frac{|\alpha|}{2}+\mu-\frac{1}{2}+ \max ( \frac{7}{6} - \lfloor n\rfloor_1 -2b_{n,6/5},0)}L_{n,\frac{17}{36}}(t) \]
	On the other hand for $I_{2}^{NN}$ we leave all of the derivatives on the nonlinearity and use Liebniz and H\"older as in (\ref{I2NoDerivs}), (\ref{Liebniz}).  Without loss of generality, we assume $|\gamma_1|\geq |\gamma_2|$, and that for some $\tilde{k}$, $\gamma_1 = \tilde{\gamma}_1 + e_k$.  We then use Cor. \ref{c:PiAndBiotSavartOperators} \textbf{(a)} on the first term and Cor. \ref{c:PiAndBiotSavartOperators} \textbf{(b)} on the second to obtain
	\[ \begin{split}
	I_2^{NN} & \leq \int \nolimits_{t/2}^t (t-s)^{-\frac{3}{2}(\frac{1}{q_1}-\frac{1}{p})}(1+t-s)^{\frac{1}{2}-(\frac{1}{q_1}-\frac{1}{p})} \big ( \|\partial_x^{\tilde{\gamma}_1}a_N\|_{L^{p_1}} + \|\partial_x^{\tilde{\gamma}_1} \vec{\omega}_N \|_{\mathbb{L}^{p_1}} \big ) \big ( \|\partial_x^{\gamma_2}a_N\|_{L^{p_3}} + \|\partial_x^{\gamma_2} \vec{\omega}_N \|_{\mathbb{L}^{p_3}} \big )ds 
	\end{split} \]
	For $1 \leq p \leq 2$ we choose $p_1 = p_3 = 3/2$ to obtain
	\[ I_2^{NN} \leq CE_n^4t^{-\frac{5}{2}(1-\frac{1}{p})-\frac{|\alpha|}{2}+\mu+\frac{1}{3}-\lfloor n \rfloor_1 - 2 b_{n,3/2}} \]
	whereas for $p=\infty$ we choose $p_1 = p_3 = 12/5$ to obtain
	\[ I_2^{NN} \leq CE_n^4t^{-\frac{5}{2}-\frac{|\alpha|}{2}+\mu+\frac{1}{2}-\lfloor n \rfloor_1 - 2 b_{n,12/5}} \]
	
	Finally we need to consider the weighted norms when $\mu = n$.  For $I_1^{NN}$ we remove all derivatives from the nonlinearity using the heat estimate, but we need to split the weight between the two terms.  For $0 < n < 1$ we split the weight evenly between the two terms and we can then apply Cor. \ref{c:PiAndBiotSavartOperators} \textbf{(b)} to both terms and pull out the first factors from the integral via
	\[ \begin{split} I_1^{NN} & \leq \int \nolimits_0^{t/2} (t-s)^{-\frac{3}{2}(\frac{1}{q_1}-\frac{1}{p})-1-\frac{|\alpha|}{2}}(1+t-s)^{\frac{1}{2}-(\frac{1}{q_1}-\frac{1}{p})} \| m_{N,i}(s) \|_{\mathring{L}^{p_1}(\frac{n}{2})}\| m_{N,j}(s) \|_{\mathring{L}^{p_2}(\frac{n}{2})} ds \\ & \leq  t^{-\frac{3}{2}(\frac{1}{q_1}-\frac{1}{p})-1-\frac{|\alpha|}{2}}(1+t)^{\frac{1}{2}-(\frac{1}{q_1}-\frac{1}{p})} \int \nolimits_0^{t/2} \big ( \|a_N\|_{L^{p_3}(\frac{n}{2})} + \| \vec{\omega}_N \|_{\mathbb{L}^{p_3}(\frac{n}{2})} \big ) \big ( \|a_N\|_{L^{p_4}(\frac{n}{2})} + \| \vec{\omega}_N \|_{\mathbb{L}^{p_4}(\frac{n}{2})} \big ) ds \end{split} \]
	whereas for $1 \leq n < 2$ we split the weight unevenly between the two terms and apply Cor. \ref{c:PiAndBiotSavartOperators} \textbf{(b)} to the term with less weight and Cor. \ref{c:PiAndBiotSavartOperators} \textbf{(c)} to the term with more weight to obtain
	\[ \begin{split} & I_1^{NN} \leq \int \nolimits_0^{t/2} (t-s)^{-\frac{3}{2}(\frac{1}{q_1}-\frac{1}{p})-1-\frac{|\alpha|}{2}}(1+t-s)^{\frac{1}{2}-(\frac{1}{q_1}-\frac{1}{p})} \| m_{N,i}(s) \|_{\mathring{L}^{p_1}(1+\frac{n-1}{2})}\| m_{N,j}(s) \|_{\mathring{L}^{p_2}(\frac{n-1}{2})} ds \\ & \leq  t^{-\frac{3}{2}(\frac{1}{q_1}-\frac{1}{p})-1-\frac{|\alpha|}{2}}(1+t)^{\frac{1}{2}-(\frac{1}{q_1}-\frac{1}{p})} \int \nolimits_0^{t/2} \big ( \|a_N\|_{L^{p_3}(1+\frac{n-1}{2})} + \| \vec{\omega}_N \|_{\mathbb{L}^{p_3}(1+\frac{n-1}{2})} \big ) \big ( \|a_N\|_{L^{p_4}(\frac{n-1}{2})} + \| \vec{\omega}_N \|_{\mathbb{L}^{p_4}(\frac{n-1}{2})} \big ) ds \end{split} \]
	In both cases the choice of $p_3 = p_4 = 6/5$ satisfies the constraints imposed by the use of Cor. \ref{c:PiAndBiotSavartOperators} \textbf{(b)}, \textbf{(c)}, so for $1 \leq p \leq \infty$ we obtain
	\[ I_1^{NN} \leq CE_n^4 t^{-\frac{5}{2}(1-\frac{1}{p})-\frac{|\alpha|}{2}+n-\frac{1}{2}+ \max( \frac{7}{6} - \lfloor n \rfloor_1 - 2 b_{n,6/5},0)}L_{n,\frac{17}{36}}(t) \]
	For $1 < n \leq 2$ we can obtain a different bound, and note that in the overlapping region $1< n < 2$ we can use the better of the two estimates.  We split the weight unevenly in a different way and apply Cor. \ref{c:PiAndBiotSavartOperators} \textbf{(b)}, \textbf{(c)} to the terms with respectively less and more weight to obtain
	\[ \begin{split} & I_1^{NN} \leq \int \nolimits_0^{t/2} (t-s)^{-\frac{3}{2}(\frac{1}{q_1}-\frac{1}{p})-1-\frac{|\alpha|}{2}}(1+t-s)^{\frac{1}{2}-(\frac{1}{q_1}-\frac{1}{p})} \| m_{N,i}(s) \|_{\mathring{L}^{p_1}(1+\frac{n-1}{3})}\| m_{N,j}(s) \|_{\mathring{L}^{p_2}(\frac{2(n-1)}{3})} ds \\ & \leq  t^{-\frac{3}{2}(\frac{1}{q_1}-\frac{1}{p})-1-\frac{|\alpha|}{2}}(1+t)^{\frac{1}{2}-(\frac{1}{q_1}-\frac{1}{p})} \int \nolimits_0^{t/2} \big ( \|a_N\|_{L^{p_3}(1+\frac{n-1}{3})} + \| \vec{\omega}_N \|_{\mathbb{L}^{p_3}(1+\frac{n-1}{3})} \big ) \big ( \|a_N\|_{L^{p_4}(\frac{2(n-1)}{3})} + \| \vec{\omega}_N \|_{\mathbb{L}^{p_4}(\frac{2(n-1)}{3})} \big ) ds \end{split} \]
	In this case the choice of $p_3 = 15/13$, $p_4 = 5/4$ satisfies the constraints imposed by the use of Cor. \ref{c:PiAndBiotSavartOperators} \textbf{(b)}, \textbf{(c)}, so for $1 \leq p \leq \infty$ we have
	\[ I_1^{NN} \leq CE_n^4 t^{-\frac{5}{2}(1-\frac{1}{p})-\frac{|\alpha|}{2}+n-\frac{1}{2}} \]
	For $I_2^{NN}$ we leave all derivatives on the nonlinearity as in (\ref{I2NoDerivs}), use Liebniz and H\"older as in (\ref{Liebniz}) and put the weight on the term having fewer derivatives to obtain
	\[ I_2^{NN} \leq \int \nolimits_{t/2}^t (t-s)^{-\frac{3}{2}(\frac{1}{q_1}-\frac{1}{p})}(1+t-s)^{\frac{1}{2}-(\frac{1}{q_1}-\frac{1}{p})} \| \partial_{x}^{\gamma_1} m_{N,i} \|_{L^{p_1}} \| \partial_x^{\gamma_2} m_{N,j} \|_{\mathring{L}^{p_2}(n)} ds \]
	where without loss of generality we assume $|\gamma_1| \geq |\gamma_2|$.  We can use Cor. \ref{c:PiAndBiotSavartOperators} \textbf{(a)} on the first term, and either Cor. \ref{c:PiAndBiotSavartOperators} \textbf{(b)} or \textbf{(c)} on the second term, depending on $n$.  In either case one obtains
	\[ I_2^{NN} \leq \int \nolimits_{t/2}^t (t-s)^{-\frac{3}{2}(\frac{1}{q_1}-\frac{1}{p})}(1+t-s)^{\frac{1}{2}-(\frac{1}{q_1}-\frac{1}{p})} \big ( \|\partial_{x}^{\gamma_1-e_{\tilde{k}}} a_N\|_{L^{p_1}} + \| \partial_{x}^{\gamma_1-e_{\tilde{k}}} \vec{\omega}_N \|_{\mathbb{L}^{p_1}} \big ) \big ( \|\partial_{x}^{\gamma_2}a_N\|_{L^{p_3}(n)} + \| \partial_{x}^{\gamma_2} \vec{\omega}_N \|_{\mathbb{L}^{p_3}(n)} \big ) ds \]
	for some index $\tilde{k}$.  For $0 < n < 1$ we use Cor. \ref{c:PiAndBiotSavartOperators} \textbf{(b)} and choose $p_1 = p_3 = 3/2$ for $1\leq p \leq 2$ to obtain
	\[ I_2^{NN} \leq CE_n^4 t^{-\frac{5}{2}(1-\frac{1}{p})-\frac{|\alpha|}{2}+n+\frac{1}{3}-\lfloor n \rfloor_1 - 2 b_{n,3/2}} \]
	whereas we can obtain the exact same bound for $1 < n < 2$ using Cor. \ref{c:PiAndBiotSavartOperators} \textbf{(c)} with $p_1=p_3 =3/2$.  We can also obtain the same bound for $1\leq n < 3/2$ using Cor. \ref{c:PiAndBiotSavartOperators} \textbf{(c)} with $p_1=2$, $p_3 =6/5$ for $1\leq p \leq 2$, and also for $3/2 < n \leq 2$ using Cor. \ref{c:PiAndBiotSavartOperators} \textbf{(c)} with $p_1=6/5$, $p_3 =2$.  Finally, for $p = \infty$ and $0< n < 7/4$ we can use Cor. \ref{c:PiAndBiotSavartOperators} \textbf{(b)} by choosing $p_1=p_3=12/5$ and we obtain
	\[ I_2^{NN} \leq CE_n^4 t^{-\frac{5}{2}-\frac{|\alpha|}{2}+n+\frac{1}{2}-\lfloor n \rfloor_1 - 2 b_{n,12/5}} \]
	and for $7/4 < n \leq 2$ we can obtain the same bound by using Cor. \ref{c:PiAndBiotSavartOperators} \textbf{(c)}.  For $3/2 < n \leq 2$ we can use Cor. \ref{c:PiAndBiotSavartOperators} \textbf{(c)} with $p_1 =3$, $p_3 = 2$ and we obtain
	\[ I_2^{NN} \leq CE_n^4 t^{-\frac{5}{2}-\frac{|\alpha|}{2}+n+\frac{1}{2}-\lfloor n \rfloor_1 - b_{n,3}-b_{n,2}} \]
	
	The excess decay rate $\mathfrak{b}_{n,p}$ can therefore be found by collecting these results and finding the slowest decay, and the rate the bounds for the terms $\rho_{NR}$ and $\vec{\omega}_{NR}$ can be obtained similarly.
\end{proof}

\appendix  

\section{Proof of the estimates on $\Pi$ and $B$ in Proposition \ref{p:BiotSavart}}

\label{app:BiotSavart}

We begin the proof with the following Lemmas:

\begin{Lemma}
	\label{l:WeightedHardyLittlewoodSobolev1}
	For $p_2,p_3$ and $n$ chosen as in part $\textbf{(b)}$ above, and given $f,g$ such that
	\[ f(x) = \int \nolimits_{\mathbb{R}^3} \frac{g(y)}{|x-y|^2}dy \]
	we have
	\[ \|f \|_{L^{p_2}(n)} \leq C \| g \|_{L^{p_3}(n)} \]
\end{Lemma}

\begin{proof}  
	The proof is based on a dyadic decomposition
	\[ \mathbb{R}^3 = \cup_{j=0}^\infty A_j \]
	where $A_0 = \{x \in \mathbb{R}^3: |x| \leq 1 \}$ and $A_j = \{ x \in \mathbb{R}^3: 2^{j-1} < |x| < 2^{j} \}$ for $ j \in \mathbb{N}$.  Let $f_i = f\chi_{A_i}$ and $g_j = g \chi_{A_j}$.  Clearly $f_i = \sum_{j \in \mathbb{N}} \Delta_{ij}$, where
	\[ \Delta_{ij}(x) = \chi_{A_i}(x)\int \nolimits_{A_j} \frac{g_j(y)}{|x-y|^2}dy  \]
	
	For the case $|i-j|\leq 1$ note that if
	\[ h_j(x) = \int \nolimits_{A_j} \frac{g_j(y)}{|x-y|^2}dy \]
	then by the Hardy-Littlewood-Sobolev inequality (\cite{stein} Theorem V.1), we have
	\[ \|\Delta_{ij}\|_{L^{p_2}} \leq \|h_j(x)\|_{L^{p_2}} \leq C\|g_j\|_{L^{p_3}}\leq \tilde{C}2^{-\alpha|i-j|} \|g_j\|_{L^{p_3}} \]
	for an $\alpha \in (0,1)$ of our choosing.  Next, we consider the case $i \geq j+2$.  By the triangle inequality
	\[ \|\Delta_{ij}\|_{L^{p_2}} \leq \Big ( \int \chi_{A_i}(x) \Big ( \int \chi_{A_j}(y) \frac{|g_j(y)|}{|x-y|^2} dy \Big )^{p_2} dx \Big )^{1/p_2}  \]
	Since $i \geq j+2$ we have $ |x-y| \geq 2^{i-2} $ and hence
	\begin{eqnarray} \nonumber
	\|\Delta_{ij}\|_{L^{p_2}} & \leq  &\Big ( \int \chi_{A_i}(x) \Big ( \int \chi_{A_j}(y) \frac{|g_j(y)|}{|x-y|^2} dy \Big )^{p_2} dx \Big )^{1/p_2}  \le  \frac{16}{2^{2i}} \Big ( \int \chi_{A_i}(x) \Big ( \int \chi_{A_j}(y) |g_j(y)| dy \Big )^{p_2} dx \Big )^{1/p_2} \\ \nonumber 
	& =  &\frac{16}{2^{2i}} \Big ( \int \chi_{A_i}(x) dx \Big )^{1/p_2} \int \chi_{A_j}(y) |g_j(y)| dy  \le\frac{16}{2^{2i}} \Big ( \int \chi_{A_i}(x) dx \Big )^{1/p_2} \Big ( \int \chi_{A_j}(y) dy \Big )^{1-1/p_3} \| g_j \|_{L^{p_3}} \\ & = & \frac{C}{2^{2i}} 2^{3i/p_2} 2^{3j(1-\frac{1}{p_3})} \|g_j\|_{L^{p_3}} = C 2^{-3(1-\frac{1}{p_3})(i-j)}\|g_j\|_{L^{p_3}}\ ,
	\end{eqnarray}	
	where in the last step we used (\ref{HardyLittlewoodSobolevIndexCondition}).  By a very similar argument, if $j \geq i+2$ we have 
	\begin{equation}
	\|\Delta_{ij}\|_{L^{p_2}} \le C 2^{-3(\frac{1}{p_3}-\frac{1}{3})(j-i)}\|g_j\|_{L^{p_3}}\ . \notag 
	\end{equation}
	Recalling the limits on the support of $f_i$ and its decomposition in terms of $\Delta_{i,j}$ we have
	the inequality:
	\begin{eqnarray}\label{eq:WeightedDecomposition} \nonumber
	\|f_i\|_{L^{p_2}(n)} & \leq  & C2^{ni} \|f_i\|_{L^{p_2}}  \le
	C2^{ni} \sum_{j\in \mathbb{N}} 2^{-|i-j| - 3(\frac{2}{3} -\frac{1}{p_3})(i-j) } \|g_j\|_{L^{p_3}} \\ & \leq & C \sum_{j\in \mathbb{N}} 2^{-|i-j| - 3(\frac{2}{3} -\frac{1}{p_3} -\frac{n}{3} )(i-j)} \|g_j\|_{L^{p_3}(n)}
	\leq \sum_{j\in \mathbb{N}} C 2^{-\alpha|i-j|} \|g_j\|_{L^{p_3}(n)}\ ,
	\end{eqnarray}
	for some $\alpha>0$, since $-1 < 3(\frac{2-n}{3}-\frac{1}{p_3}) < 1$.
	Considering now $f$ itself, we have
	\[ \|f\|_{L^{p_2}(n)}^{p_2} = \sum_i \|f_i\|_{L^{p_2}(n)}^{p_2} \]
	and since
	\begin{equation}\label{eq:BiotSavart}
	\sum_i \|f_i\|_{L^{p_2}(n)}^{p_2}  \leq \sum_{i} \Big ( \sum_j C2^{-\alpha|i-j|} \|g_j \|_{L^{p_3}(n)} \Big )^{p_2}  = \sum_{i} \Big ( \sum_j C2^{-\alpha|i-j|(1-\frac{1}{p_3})} 2^{-\alpha|i-j|\frac{1}{p_3}} \|g_j \|_{L^{p_3}(n)} \Big )^{p_2} \end{equation}
	we can then apply H\"older's inequality and interchange the order of summation to obtain
	\begin{equation} \label{eq:BiotSavart2} \|f\|_{L^{p_2}(n)}^{p_2} \leq C \sum_{i} \Big [ \sum_{j} 2^{-\alpha|i-j|} \|g_j \|_{L^{p_3}(n)}^{p_3} \Big ]^{\frac{p_2}{p_3}}
	\leq C  \Big [ \sum_{i} \sum_{j} 2^{-\alpha|i-j|} \|g_j \|_{L^{p_3}(n)}^{p_3} \Big ]^{\frac{p_2}{p_3}}  \le 
	C  \|g \|_{L^{p_3}(n)}^{p_2}
	\end{equation}
	where in the last step we compute the geometric sum and use convexity since
	$\frac{p_2}{p_3} = \frac{3}{3-p_3} > 1$.
\end{proof}

\begin{Lemma}
	\label{l:WeightedHardyLittlewoodSobolev2}
	For $1 < p_3 < p_2 < \infty$ and $n \in [0,2)$ chosen such that
	\[ \frac{2-n}{3} < \frac{1}{p_3}  < \frac{3-n}{3} \]
	and given $f,g$ such that
	\[ f(x) = \int \nolimits_{\mathbb{R}^3} \frac{g(y)}{|x-y|}dy \]
	we have
	\[ \|f \|_{L^{p_2}(n-1)} \leq C \| g \|_{L^{p_3}(n)} \]
\end{Lemma}

\begin{proof}
	Defining $f_i$, $g_j$, $\Delta_{ij}$ and $h_j$ analogously to the above much of the proof follows in almost identical fashion.  The key difference arises from the fact that $p_3$ lies in a different range in this case.  In the step analogous to \eqref{eq:WeightedDecomposition}, we have
	\begin{eqnarray}
	\|f_i\|_{L^{p_2}(n-1)}  &\leq & C2^{(n-1)i} \|f_i\|_{L^{p_2}}   \leq C2^{(n-1)i} \sum_{j\in \mathbb{N}} 2^{-\frac{1}{2}|i-j| - 3(\frac{1}{2} -\frac{1}{p_3})(i-j) } \|g_j\|_{L^{p_3}(1)} \notag \\ & \leq & C \sum_{j\in \mathbb{N}} 2^{-\frac{1}{2}|i-j| - 3(\frac{5}{6} -\frac{1}{p_3} -\frac{n}{3})(i-j)} \|g_j\|_{L^{p_3}(n)} \leq \sum_{j\in \mathbb{N}} C 2^{-\alpha|i-j|} \|g_j\|_{L^{p_3}(n)}\ , \notag 
	\end{eqnarray}
	for some $\alpha>0$, since $-\frac{1}{2} < 3(\frac{5}{6}-\frac{1}{p_3} - \frac{n}{3} ) < \frac{1}{2}$.
	The estimate in the lemma now follows by a summation similar to that in \eqref{eq:BiotSavart} and \eqref{eq:BiotSavart2}.
\end{proof}
\begin{proof}[Proof of Proposition \ref{p:BiotSavart}]
	The operators $\partial_{x_i} \Pi$ and $\partial_{x_i} B$ are singular integral operators formed by kernels of Calderon-Zygmund type, so part \textbf{(a)} follows from Theorem II.3 in \cite{stein}.  Examining the form of the $\Pi$ and $B$ operators, we see that part \textbf{(b)} follows directly the result of the Lemma $\ref{l:WeightedHardyLittlewoodSobolev1}$.
	
	For part \textbf{(c)}, we use our modified versions of $B.2$, $B.3$ to complete the analogous proof in \cite{gallaywayne2}.  Write
	\[ \big ( \Pi a \big )_{i} = - \frac{1}{4\pi} \int \nolimits_{\mathbb{R}^3} \left ( \frac{x_i-y_i}{|x-y|^3} - \frac{x_i}{|x|^3} \right ) a(y)dy  \]
	using the moment zero condition.  Using the identity
	\[ |x|^3(x_i-y_i) - |x-y|^3x_i = (x_i-y_i)|x|^2(|x|-|x-y|) + |x-y|(2x_i(x\cdot y)-y_i|x|^2-x_i|y|^2) \]
	it follows that
	\[ \begin{split}
	\big | |x|^3(x_i-y_i) - |x-y|^3x_i\big | & \leq C |x-y||x||y|(|x|+|y|) \leq C (|x-y||x|^2|y|+|x-y|^2|x||y|)
	\end{split} \]
	and hence $|(\Pi a)_i| \leq C(u_1 + u_2)$ where
	\[ u_1(x) = \frac{1}{|x|} \int \frac{|y||a(y)|}{|x-y|^2}dy \hspace{.75 cm} \text{ , } \hspace{.75 cm} u_2(x) = \frac{1}{|x|^2} \int \frac{|y||a(y)|}{|x-y|}dy \]
	Therefore, using Lemma $\ref{l:WeightedHardyLittlewoodSobolev1}$, $\ref{l:WeightedHardyLittlewoodSobolev2}$ with $f_1 = |x|u_1$, $f_2 = |x|^2u_2$ and $g_1 = g_2 = |y|$ $ |a(y)|$ we have
	\begin{eqnarray}
	\|\Pi a\|_{L^{p_2}(n)} & \leq  &C \big \|\chi_{|\cdot|\leq 1} \Pi a \big \|_{L^{p_2}} + C \big \|\chi_{|\cdot|>1} |\cdot|^{n}\Pi a \big \|_{L^{p_2}} \notag \\ \nonumber & \leq & C \big \| \Pi a \big \|_{L^{p_2}} + C \Big \|\chi_{|\cdot|>1} |\cdot|^{n}u_1 \Big \|_{L^{p_2}} + C \Big \|\chi_{|\cdot|>1} |\cdot|^{n}u_2 \Big \|_{L^{p_2}} \\ \nonumber
	& \leq  &C \big \| \Pi a \big \|_{L^{p_2}(n-1)} + C \big \|f_1 \big \|_{L^{p_2}(n-1)} + C \big \|f_2 \big \|_{L^{p_2}(n-2)} \\ \nonumber & \leq &C \big \| a \big \|_{L^{p_3}(n-1)} + C \big \|g_1 \big \|_{L^{p_3}(n-1)} + C \big \|g_2 \big \|_{L^{p_3}(n-1)} \leq C \big \| a \big \|_{L^{p_3}(n)} 
	\end{eqnarray}
	The proof for $B\vec{\omega}$ is analogous.
	
\end{proof}

\section{Proof of heat estimate in Proposition \ref{p:HeatEst}}

\label{app:HeatEst}

\begin{proof}
	We prove that in general dimension $d\geq 1$
	\[ \| \partial_x^{\alpha} K_{\nu}(t) \ast f \|_{\mathring{L}^p(\mu)} \leq C (\nu t)^{-\frac{|\alpha|}{2}-\frac{d}{2}(\frac{1}{q}-\frac{1}{p})- \frac{n-\mu}{2}} \| f\|_{L^q(n)} \]
	and the result then holds by estimating the $\mathring{L}^p(\mu)$ norms separately for $\nu t < 1$ and $\nu t \geq 1$ using different values for $q$.  Write
	\[ \begin{split}
	\|\partial_{x}^{\alpha} K_{\nu}(x,t)\ast f \|_{\mathring{L}^p_x(\mu)} = \| \Big ( \int_{|y| \geq \sqrt{\nu t}} + \int_{|y| < \sqrt{\nu t}} \Big ) \partial_{x}^{\alpha} K_{\nu}(x-y,t) f(y) dy \|_{\mathring{L}^p_x(\mu)}  \leq  S_1 + S_2 + S_3 
	\end{split} \]
	where 
	\[ \begin{split} 
	S_1 & = \| \int \nolimits_{\mathbb{R}^d} \partial_{x}^{\alpha} K_{\nu}(x-y,t) f(y) \chi_{|y| \geq \sqrt{\nu t}} dy \|_{\mathring{L}^p_x(\mu)} \\
	S_2 & = \sum_{|\beta| \leq \tilde{n}} \frac{1}{\beta!}\| \int \nolimits_{\mathbb{R}^d} \partial_{x}^{\alpha+\beta} K_{\nu}(x,t) y^{\beta} f(y) \chi_{|y| < \sqrt{\nu t}} dy \|_{\mathring{L}^p_x(\mu)} \\
	S_3 & = \sum_{|\beta| = \tilde{n} +1 } \frac{ \tilde{n}+1}{\beta!} \Big \| \int \nolimits_{|y| <\sqrt{\nu t}} y^{\beta} f(y) \int \nolimits_0^1 (1-s)^{\tilde{n}}\partial_x^{\alpha+\beta}K_{\nu}(x-sy,t) dsdy \Big \|_{\mathring{L}^p_x(\mu)}
	\end{split} \]
	and where we used Taylor's theorem
	\[ \begin{split} \partial_{x}^{\alpha}K_{\nu}(x-y,t) = \sum_{|\beta| \leq \tilde{n}} (-1)^{|\beta|} \frac{\partial_{x}^{\alpha+\beta} K_{\nu}(x,t)}{\beta !} y^{\beta} + \sum_{|\beta| = \tilde{n} +1 }  (-1)^{\tilde{n}+1} \frac{\tilde{n}+1}{\beta!} y^{\beta} \int \nolimits_0^1 (1-s)^{\tilde{n}}\partial_x^{\alpha+\beta}K_{\nu}(x-sy,t) ds \end{split} \]
	
	For $S_1$, we change variables and use $|\tilde{x}| \leq |\tilde{x}-\tilde{y}|+|\tilde{y}|$:
	\[ \begin{split}
	\|\int_{\mathbb{R}^d} \partial_{x}^{\alpha} K_{\nu}(x-y,t) f(y)\chi_{|y|\geq \sqrt{\nu t}} dy \|_{\mathring{L}^p_x(\mu)} & \leq 
	(\nu t)^{\frac{\mu}{2}-\frac{|\alpha|}{2}+\frac{d}{2p}} \Big [  \|\int_{\mathbb{R}^d} |\tilde{x}-\tilde{y}|^{\mu} \partial_{\tilde{x}}^{\alpha}K_{\nu}(\tilde{x}-\tilde{y}) f(\sqrt{\nu t}\tilde{y})\chi_{|\tilde{y}|\geq 1} d\tilde{y} \|_{\mathring{L}^p_{\tilde{x}}}  \\ & \hspace{3 cm} +
	\|\int_{\mathbb{R}^d} |\tilde{y}|^{\mu} \partial_{\tilde{x}}^{\alpha}K_{\nu}(\tilde{x}-\tilde{y}) f(\sqrt{\nu t}\tilde{y})\chi_{|\tilde{y}|\geq 1} d\tilde{y} \|_{\mathring{L}^p_{\tilde{x}}} \Big ] \end{split}  \] 
	We can then use Young's inequality and change back to our original variables:
	\[ \begin{split}
	\|\int_{\mathbb{R}^d} \partial_{x}^{\alpha} K_{\nu}(x-y,t) f(y)\chi_{|y|\geq \sqrt{\nu t}} dy \|_{\mathring{L}^p_x(\mu)} & \leq (\nu t)^{\frac{\mu}{2}-\frac{|\alpha|}{2}+\frac{d}{2p}} \Big [ \|\partial_{\tilde{x}}^{\alpha} K_{\nu}(\tilde{x}) \|_{\mathring{L}^{\frac{pq}{pq+q-p}}_{\tilde{x}}(\mu)} \|f(\sqrt{\nu t}\tilde{y})\chi_{|\tilde{y}|\geq 1} \|_{L^q_{\tilde{y}}} \\ & \hspace{3 cm} + \|\partial_{\tilde{x}}^{\alpha}K_{\nu}(\tilde{x}) \|_{L^{\frac{pq}{pq+q-p}}_{\tilde{x}}} \||\tilde{y}|^{\mu}f(\sqrt{\nu t}\tilde{y})\chi_{|\tilde{y}|\geq 1} \|_{L^q_{\tilde{y}}} \Big ] \\ & \leq C (\nu t)^{\frac{\mu-n}{2}-\frac{|\alpha|}{2}+\frac{d}{2}\left ( \frac{1}{p}-\frac{1}{q} \right )} \|f(y) \|_{\mathring{L}^q_{y}(n)}
	\end{split} \]
	
	For $S_2$, we can factor out the $y$ dependent terms from the $L^p_{x}$ norm
	\[ \begin{split}
	\| \int_{\mathbb{R}^d} \partial_{x}^{\alpha+\beta} K_{\nu}(x,t) y^{\beta} f(y) \chi_{|y| < \sqrt{\nu t}} dy \|_{\mathring{L}^p_x(\mu)} \leq (\nu t)^{-\frac{|\alpha|+|\beta|}{2}+\frac{\mu}{2}+\frac{d}{2} (\frac{1}{p}-1)} \| \partial_{\tilde{x}}^{\alpha+\beta} K_{\nu}(\tilde{x} ) \|_{L^p_{\tilde{x}}(\mu)} \Big |\int_{|y|< \sqrt{\nu t}} y^{\beta} f(y) dy \Big | 
	\end{split} \]
	Since $|\beta| \leq \tilde{n} < n - d(1-\frac{1}{q})$ we use the zero moment property to obtain 
	\[ \begin{split}
	\Big |\int \nolimits_{|y|< \sqrt{\nu t}} y^{\beta} f(y) dy \Big | & = 
	\Big |\int \nolimits_{|y|\geq \sqrt{\nu t}} \text{ } \frac{|y|^{n-|\beta|}}{|y|^{n-|\beta|}}|y|^{|\beta|} |f(y) | dy \Big | \\ & \leq \| |y|^{-n+|\beta|} \chi_{|y|\geq \sqrt{\nu t}} \|_{L^{\frac{q}{q-1}}_y} \| f \|_{\mathring{L}^q(n)} \leq C (\nu t)^{-\frac{n}{2} + \frac{|\beta|}{2} - \frac{d}{2} (\frac{1}{q} -1)}\|f\|_{\mathring{L}^q(n)}
	\end{split} \]
	
	For $S_3$, write
	\[ \begin{split}
	\| \int_{|y|<\sqrt{\nu t}} \Phi_{\beta} y^{\beta} f(y)dy \|_{\mathring{L}^p_x(\mu)} & = \Big \| |x|^{\mu} \int_{|y|< \sqrt{\nu t}} y^{\beta} f(y) \Big [\int_0^1 (1-s)^{\tilde{n}} \partial_{x}^{\alpha+\beta}K_{\nu}(x-sy,t) ds \Big ] dy \Big \|_{L^p_{\tilde{x}}} \\ & = (\nu t)^{\frac{\mu}{2}-\frac{|\alpha|}{2}+\frac{d}{2p}} \Big \| |\tilde{x}|^{\mu} \int_{|\tilde{y}|< 1} \tilde{y}^{\beta} f(\sqrt{\nu t}\tilde{y}) \Big [\int_0^1 (1-s)^{\tilde{n}} \partial_{\tilde{x}}^{\alpha+\beta}K_{1}(\tilde{x}-s\tilde{y}) ds \Big ] d\tilde{y} \Big \|_{L^p_{\tilde{x}}} \\ & \leq (\nu t)^{-\frac{|\alpha|}{2}+\frac{\mu}{2}+\frac{d}{2p}} \Big \| |\tilde{x}|^{\mu} \int_{|\tilde{y}|< 1} |\tilde{y}|^{\tilde{n}+1} |f(\sqrt{\nu t}\tilde{y})| \Big [\int_0^1  \Big | \partial_{\tilde{x}}^{\alpha+\beta}K_{1}(\tilde{x}-s\tilde{y}) \Big | ds \Big ] d\tilde{y} \Big \|_{L^p_{\tilde{x}}}
	\end{split} \]
	Now using the fact that $s \leq 1,|\tilde{y}| \leq 1$ we have
	\[ \Big | \partial_{\tilde{x}}^{\alpha+\beta} K_{1}(\tilde{x}-s\tilde{y}) \Big | = \Big | \sum_{j=0}^{\tilde{n}+1+|\alpha|} c_j (\tilde{x}_j-s\tilde{y}_j)^{j} \exp [ -\frac{|\tilde{x}-s\tilde{y}|^2}{4} ] \Big | \leq C (1+|\tilde{x}|)^{\tilde{n}+1+|\alpha|} \exp [ -\frac{|\tilde{x}-s\tilde{y}|^2}{4} ] \]
	and
	\[ \exp \big [ -\frac{|\tilde{x}-s\tilde{y}|^2}{4} \big ] = \exp \big [ -\frac{|\tilde{x}|^2}{8} \big ]\exp \big [ -\frac{|\tilde{x}|^2}{8} + \frac{s\tilde{x}\cdot \tilde{y} }{2} - \frac{s^2|\tilde{y}|^2}{4} \big ] \leq C \exp \big [ -\frac{|\tilde{x}|^2}{8} \big ] \]
	If we let $\delta > 0$ be such that $n - d(1-\frac{1}{q}) + \delta < \tilde{n} +1$, then we have
	\[ \begin{split}
	\| \int_{|y|<\sqrt{\nu t}} \Phi_{\beta} y^{\beta} f(y)dy \|_{\mathring{L}^p_x(\mu)} \leq (\nu t)^{-\frac{|\alpha|}{2}+\frac{\mu}{2}+\frac{d}{2p}} \Big \| (1+|\tilde{x}|)^{\mu+\tilde{n}+1+|\alpha|} \int_{|\tilde{y}|< 1} |\tilde{y}|^{\tilde{n}+1} |f(\sqrt{\nu t}\tilde{y})| \exp [ -\frac{|\tilde{x}|^2}{8}] d\tilde{y} \Big \|_{L^p_x} \end{split} \]
	and since the integral in $\tilde{y}$ no longer depends on $\tilde{x}$ we have
	\[ \begin{split}
	\| \int_{|y|<\sqrt{\nu t}} \Phi_{\beta} y^{\beta} f(y)dy \|_{\mathring{L}^p_x(\mu)} & = (\nu t)^{-\frac{|\alpha|}{2}+\frac{\mu}{2}+\frac{d}{2p}} \Big \| (1+|\tilde{x}|)^{\mu+\tilde{n}+1+|\alpha|} \exp [ -\frac{|\tilde{x}|^2}{8}] \Big \|_{L^p_x} \int_{|\tilde{y}|< 1} |\tilde{y}|^{\tilde{n}+1} |f(\sqrt{\nu t}\tilde{y})| d\tilde{y} \\ & \leq C(\nu t)^{-\frac{|\alpha|}{2}+\frac{\mu}{2}+\frac{d}{2p}} \int \nolimits_{|\tilde{y}|< 1} |\tilde{y}|^{n - d(1-\frac{1}{q}) +\delta} |f(\sqrt{\nu t}\tilde{y})| d\tilde{y} \\ & \leq C(\nu t)^{-\frac{|\alpha|}{2}-\frac{n-\mu}{2}-\frac{d}{2}(\frac{1}{q}-\frac{1}{p})} \|f\|_{\mathring{L}^q(n)} 
	\end{split} \]
	
\end{proof}

\section{Proof of the heat-wave estimates}

\subsection{Proof of Proposition \ref{p:HypParaEst}}

\label{app:HeatWaveEst}

We first obtain point-wise estimates.  Recalling the form of the Kirchhoff formula, we need a bound on the spherical integral of the Gaussian, so we begin with the following estimate:

\begin{Lemma}
	\label{l:WaveIntExp} 
	In any dimension $d\geq 1$,
	\[ \int \nolimits_{|z| = 1} \re^{\frac{-|x+ctz|^2}{at}} dS(z) \leq C(d) \left ( 1 + t \right)^{-\frac{d-1}{2}} \re^{-\frac{(|x|-ct)^2}{4at}}. \]
\end{Lemma}
\begin{proof}
	We recall the proof given by \cite{hoffzumbrun2}.  First note that the integral above is rotationally invariant so that we may, without loss of generality, set $x = |x|e_1$.  It then suffices to integrate over the set $\{z\,:\, |z| = 1, \,\, z_1 \leq 0 \}$, since the other part is smaller, and we'll relabel $z$ with $-z$ for convenience.  Note that for such $x$ and $z$,
	\[ \begin{split}
	3|x - t z|^2   &\geq (|x| - ct)^2 + 2|x - ctz|^2 \geq  (|x| - ct)^2 + 2 (|x|^2 - 2 |x| z_1 ct + c^2t^2 |z|^2)\notag\\
	&\geq (|x| - ct)^2 + 2c^2t^2(1-z_1^2 ) \notag \\
	&\geq (|x| - ct)^2 + c^2t^2(1-z_1 ) \notag
	\end{split}
	\]
	This can then be used to obtain the estimate
	\[ \begin{split}
	\int_{|z| = 1, z_1\geq0} \re^{\frac{-||x| e_1 - ctz|^2}{at}} dS(z) &\leq \re^{\frac{-(|x| - ct)^2}{3at}}\int_{|z| = 1, z_1\geq0} \re^{\frac{c^2t(z_1 - 1)}{3a}} dS(z) =C(d) \Big ( \frac{ct}{a} \Big )^{-\frac{d-1}{2}} \re^{\frac{-(|x| - ct)^2}{3at}}
	\end{split}	\]
	by a simple calculation using the parameterization $z_1 = \sqrt{1-(z_2^2+...+z_d^2)}$ of the hemispherical integral.  We can improve this somewhat when $0<ct<1$.  Note that for $|z|=1$, we have
	\[ |x+ct z|^2 = |x|^2 + c^2t^2 -2 ct |x| \ge \frac{|x|^2}{4} - 2 c^2t^2\ ,
	\]
	so if $0<ct<a$,
	\[ \int \nolimits_{|z|=1} e^{-\frac{|x+ctz|^2}{at} } dS(z) \le \int \nolimits_{|z|=1} e^{-\frac{|x|^2}{4 at}} e^{\frac{2c^2t}{a}} dS(z)
	\le C(d) e^{-\frac{|x|^2}{4 at}}\ . \]
	
	On the other hand, for $a \le ct \le 1$, $(ct/a)^{-\frac{d-1}{2}} \le 1$, so we have
	\[ \int \nolimits_{|z|=1} e^{-\frac{|x+ctz|^2}{at} } dS(z) \le C(d) e^{-\frac{(|x|-ct)^2}{3at}}\ \]
	
\end{proof}

\begin{proof}[ Proof of Proposition \ref{p:HypParaEst}]
	Working in general odd dimension $d\geq 3$, we first derive pointwise bounds for the Green's functions $w*K_{\nu t}, \partial_t w*K_{\nu t},$ and $\partial_{t}^2w*K_{\nu t}$.  Using $(\ref{Kirchhoff})$ and the above lemmas, we find
	\begin{align}
	|w*K_{\nu t}(x) | & \leq \sum_{0\leq |\alpha| \leq \frac{d-3}{2} } \Big | b_{\alpha,0} (ct)^{|\alpha|+1} \int \nolimits_{|z| = 1} D^{\alpha} K_{\nu t}(x + ctz) z^{\alpha} dS(z) \Big | \notag\\
	& \leq C \sum_{0\leq |\alpha| \leq \frac{d-3}{2} } (ct)^{\frac{|\alpha|}{2}+1-\frac{d}{2}} \int \nolimits_{|z| = 1} e^{-\frac{|x+ctz|^2}{5\nu t}} dS(z) \notag \\ 
	&\leq  C(1+t)^{\frac{d-3}{4}} t^{1-\frac{d}{2}}(1+t)^{-\frac{d-1}{2}} \re^{-\frac{(|x| - ct)^2}{20 \nu t}}     \notag
	\end{align}
	where we have used the fact that 
	\[ D^{\alpha}_{\xi} e^{-\frac{|\xi|^2}{4}} \leq C \re^{-\frac{|\xi|^2}{5}}  \]
	for some constant $C$.  Using the analogous bounds we then find 
	\begin{align}
	|\partial_tw*K_{\nu t}(x)| \leq  t^{-\frac{d}{2}}(1+t)^{-\frac{d-1}{4}} \re^{-\frac{(|x| - ct)^2}{20 \nu t}} \hspace{0.5 cm} \text{ and } \hspace{0.5 cm} 
	|\partial_t^2w*K_{\nu t}(x)| \leq  C t^{-\frac{d+1}{2}}(1+t)^{-\frac{d-1}{4}} \re^{-\frac{(|x| - ct)^2}{20 \nu t}}     \notag
	\end{align}
	The desired $\mathring{L}^q(n)$ bounds then follow from an estimate of the $\mathring{L}^{q}(n)$ norm of the translating exponential:
	\begin{align}
	\| \re^{-\frac{(|\cdot| - ct)^2}{20 \nu t}}\|^q_{\mathring{L}^q(n)} &= \int \nolimits_{\mathbb{R}^d} (|x|^n)^q\re^{-p\frac{(|x| - ct)^2}{20 \nu t}} dx = C\int \nolimits_0^\infty (r^n)^q \re^{-q\frac{(r - ct)^2}{20 \nu t}} r^{d-1} dr\notag\\
	&=\int \nolimits_{0}^\infty (t^{1/2}\tl r)^{nq+d-1} \re^{-\frac{q(\tl r-ct^{1/2})^2}{20 \nu}}  t^{1/2}d\tl r,\qquad\qquad r = \tl r t^{1/2}\notag\\
	&= t^{(nq+d)/2 }\int_{0}^\infty \tl r^{nq+d-1} \re^{-\frac{q(\tl r-ct^{1/2})^2}{20 \nu}}  d\tl r = t^{(nq+d)/2 }\int_{-t^{1/2}}^\infty (\rho + ct^{1/2})^{nq+d-1}\re^{-\frac{q\rho^2}{20\nu}}  d\rho, \qquad \tl r = \rho + ct^{1/2}\notag\\
	&\lesssim t^{(nq+d)/2} \int \nolimits_{\R} (\rho^{nq+d-1} + t^{(nq+d-1)/2})\re^{-\frac{q\rho^2}{20\nu}}  d\rho   \lesssim C t^{(nq+d)/2}(1 + t^{(nq+d-1)/2}) \notag
	\end{align}
	and hence 
	\[ \| \re^{-\frac{(|\cdot| - ct)^2}{20 \nu t}}\|_{\mathring{L}^q(n)}\lesssim t^{\frac{n}{2} + \frac{d}{2q}}(1+t)^{\frac{n}{2}+\frac{d-1}{2q}} \]		
\end{proof}

\subsection{Proof of Proposition \ref{p:HypParaEst2}}

\label{app:HeatWaveEst2}

\begin{proof}
	The proof follows by putting one of the derivatives in (\ref{Kirchhoff}) on $\rho_0$.  Specifically, we have
	\begin{align}
	\Big | \partial_t^2w*K_{\nu t} \ast \rho_0 \Big | & \leq \sum_{1 \leq |\tilde{\alpha}|\leq \frac{d+1}{2}} c_{\tilde{\alpha}} t^{|\tilde{\alpha}|-1} \Big | \int \nolimits_{|z|=1} D^{\tilde{\alpha}}_x \big [ K_{\nu t} \ast \rho_0 (x+ctz) \big ] z^{\tilde{\alpha}} dS(z) \Big | \notag \\ & \leq \sum_{0 \leq |\alpha|\leq \frac{d-1}{2}} \sum_{j= 1}^{d} c_{\alpha + e_j} t^{|\alpha|}\int \nolimits_{|z|=1} D^\alpha_x D_{x_j} \big [ K_{\nu t} \ast \rho_0 (x+ ctz) \big ] z_j z^\alpha dS(z) \notag \\ & \leq \sum_{0 \leq |\alpha|\leq \frac{d-1}{2}} \sum_{j= 1}^{d} c_{\alpha + e_j} t^{\frac{|\alpha|}{2}-\frac{d}{2}}\int \nolimits_{\mathbb{R}^3} \Big [ \int \nolimits_{|z|=1} \exp \big [ -\frac{|x-y+ctz|^2}{5\nu t} \big ]dS(z) \Big ] \big | D_{x_j} \rho_0(y) \big | dy \notag \end{align}
	and we can then use Lemma $\ref{l:WaveIntExp}$
	\begin{align} \Big | \partial_t^2w*K_{\nu t} \ast \rho_0 \Big | & \leq C \sum_{j= 1}^{d} t^{-\frac{d}{2}}(1+t)^{\frac{d-1}{4}-\frac{d-1}{2}}\int \nolimits_{\mathbb{R}^3} \exp \big [ -\frac{(|x-y|-ct)^2}{20\nu t} \big ] \big | D_{x_j} \rho_0(y) \big | dy \notag \\ & = C \sum_{j= 1}^{d} t^{-\frac{d}{2}}(1+t)^{\frac{d-1}{4}-\frac{d-1}{2}} \exp \big [ -\frac{(|\cdot|-ct)^2}{20\nu t} \big ] \ast \big | D_{x_j} \rho_0 \big | (x) \notag
	\end{align}
	hence for small times we can make the estimate
	\[ \begin{split}
	\|\partial_t^2w*K_{\nu t} \ast \rho_0\|_{\mathring{L}^p(\mu)} & \leq C \sum_{j= 1}^{d} t^{-\frac{d}{2}}(1+t)^{\frac{d-1}{4}-\frac{d-1}{2}} \big \| \exp \big [ -\frac{(|\cdot|-ct)^2}{20\nu t} \big ] \ast \big | D_{x_j} \rho_0 \big | \big \|_{\mathring{L}^p(\mu)} \notag \\ & \leq  C \sum_{j=1}^{d} t^{-\frac{d}{2}(\frac{1}{q}-\frac{1}{p})}(1+t)^{\frac{d-1}{4}-\frac{d-1}{2}(\frac{1}{q}-\frac{1}{p})} \Big [ \| D_{x_j} \rho_0 \|_{\mathring{L}^q(\mu)} + t^{\frac{\mu}{2}}(1+t)^{\frac{\mu}{2}} \| D_{x_j} \rho_0 \|_{L^q} \Big ] \notag 
	\end{split} \]
	whereas for large times we use the Young's inequality together with the estimate in Prop. \ref{p:HypParaEst}.
	
\end{proof}

\section{Explicit calculations of the Hermite profiles}

\label{app:AsympProfs}

\subsection{Explicit functional form for the hyperbolic-parabolic Hermite profiles}

The functions $\rho_1$, $a_1$, $\rho_2$, $a_2$, $\Pi a_1$ and $\Pi a_2$ are given by the following explicit formulas:

\begin{gather}
\rho_1(x,t) = \frac{(|x|-ct)e^{-\frac{(|x|-ct)^2}{4(1+\nu t)}}+(|x|+ct)e^{-\frac{(|x|+ct)^2}{4(1+\nu t)}}}{2|x|(4\pi (1+\nu t))^{3/2}} \notag \\ a_1(x,t) = \frac{c}{2|x|(4\pi(1+\nu t))^{3/2}} \Big [ \big [ \frac{(|x|+ct)^2}{2(1+\nu t)}-1 \big ] e^{-\frac{(|x|+ct)^2}{4(1+\nu t)}} -\big [\frac{(|x|-ct)^2}{2(1+\nu t)}-1 \big ]e^{-\frac{(|x|-ct)^2}{4(1+\nu t)}} \Big ] \notag
\end{gather}

\begin{gather} \rho_{2}(x,t) = \frac{1}{(4\pi)^{3/2}(1+\nu t)^{1/2}} \frac{e^{-\frac{(|x|+ct)^2}{4(1+\nu t)}}-e^{-\frac{(|x|-ct)^2}{4(1+\nu t)}}}{c|x|} \notag \\ a_{2}(x,t) = \frac{1}{(4\pi(1+\nu t))^{3/2}} \frac{(|x|-ct)e^{-\frac{(|x|-ct)^2}{4(1+\nu t)}}+(|x|+ct)e^{-\frac{(|x|+ct)^2}{4(1+\nu t)}}}{2|x|} \notag \end{gather}

\[ \begin{split} \Pi a_1 = \frac{cx}{(4\pi)^{3/2}|x|^3(1+\nu t)^{1/2} } \Big [ e^{-\frac{(|x|-ct)^2}{4(1+\nu t)}}\Big ( \frac{|x|(|x|-ct)}{2(1+\nu t)}+1 \Big ) -e^{-\frac{(|x|+ct)^2}{4(1+\nu t)}}\Big ( \frac{|x|(|x|+ct)}{2(1+\nu t)}+1 \Big ) \Big ] \end{split} \]

\[ \begin{split} \Pi a_{2} = \frac{1}{(4\pi)^{3/2}}  \frac{x}{|x|^3} \Big ( - |x| \frac{e^{-\frac{(|x|+ct)^2}{4(1+\nu t)}}+e^{-\frac{(|x|-ct)^2}{4(1+\nu t)}}}{(1+\nu t)^{1/2}} +\mathsf{Erf}(\frac{|x|-ct}{2(1+\nu t)^{1/2}}) + \mathsf{Erf}(\frac{|x|+ct}{2(1+\nu t)^{1/2}}) \Big ) \end{split} \]
where
\[ \mathsf{Erf}(\zeta) = 2 \int \nolimits_0^{\zeta} e^{-z^2}dz \]

Given a spherically symmetric initial condition $(u_0,0)^T$, the solution to the wave equation is given by
\begin{equation} \label{RadialWaveSoln} u(x,t) = \frac{(|x|-ct)u_0(||x|-ct|)+(|x|+ct)u_0(|x|+ct)}{2|x|} \end{equation}

Taking $u_0$ to be $K_{\nu}(t) \ast \phi_0$, we obtain the equation for $\rho_1$.  We compute $a_1$ by plugging $u_0 = K_{\nu}(s) \ast \phi_0$ into $(\ref{RadialWaveSoln})$, taking the derivative of $u(x,t)$ with respect to $t$, multiplying by $-1$ and then setting $s = t$.

To compute $\Pi a_1$, note that
\[ \Pi a_1 = \nabla (\Delta^{-1} a_1) \]
and that since $a_1$ is spherically symmetric it suffices to compute $\nabla u$, where
\[ \frac{1}{r^2} \frac{\partial }{\partial r} \big [ r^2 \frac{\partial u}{\partial r} \big ] = a_1 \]
The result follows by computing an indefinite radial integral, ensuring the integal is zero at the origin, and making use of
\begin{equation} \label{Gradient} \nabla u = \frac{x}{r} \frac{\partial u}{\partial r}  \end{equation}

To calculate the explicit forms of $\rho_2$ and $a_2$ we use the fact that the solution of the wave equation with a spherically symmetric initial condition of the form $(0,u_0(r))^T$ is given by
\begin{equation} \label{RadialWaveSolnIC2} u(x,t) = -\int \nolimits_0^t \frac{(|x|-cs)u_0(||x|-cs|)+(|x|+cs)u_0(|x|+cs)}{2|x|} ds \end{equation}
hence we have the result above for $\rho_2$, and $a_2$ is found by using $(\ref{RadialWaveSoln})$.  $\Pi a_2$ is computed using the same method used for $\Pi a_1$.

\subsection{Explicit functional form for the divergence-free vector field Hermite profiles}

We compute $B\vec{g}_i$ where
\[ \vec{g}_i :=   \frac{1}{(4\pi)^{3/2}(1+\epsilon t)^{3/2}} \nabla \times \big ( e^{-\frac{|x|^2}{4(1+\epsilon t)}} \text{ } \vec{e}_{i}  \big ) \]
and note that in view of the definitions in Table \ref{tab:DivFreeHermiteExp} the terms $B \vec{f}_{\tilde{\alpha},j}$ can be computed by taking appropriate derivatives.  One can check that the function
\[  \frac{1}{(4\pi)^{3/2}(1+\epsilon t)^{3/2}} \Big [ e^{-\frac{|x|^2}{4(1+\epsilon t)}}\vec{e}_i - \partial_{x_i} \nabla (\Delta^{-1} e^{-\frac{|x|^2}{4(1+\epsilon t)}} )  \Big ] \]
has curl equal to $\vec{g}_i$ since the second term is a gradient, hence has zero curl.  Furthermore the divergence of the above expression is zero, since the divergence and gradient cancel the inverse Laplacian in the second term.  As before we can compute the inverse Laplacian of the Gaussian term by exploiting the spherical symmetry and we get
\[ \frac{\partial u}{\partial r} = -\frac{2(1+\epsilon t)}{r}e^{-\frac{r^2}{4(1+\epsilon t)}}+ \frac{2(1+\epsilon t)}{r^2} \int \nolimits_0^r e^{-\frac{z^2}{4(1+\epsilon t)}}dz \]
so using (\ref{Gradient}) we have
\[ B \vec{g}_{i} = \frac{1}{(4\pi)^{3/2}} \Big [ \frac{e^{-\frac{|x|^2}{4(1+\epsilon t)}}}{(1+\epsilon t)^{3/2}}\vec{e}_i - \partial_{x_i}\Big [ \frac{x}{|x|^3} \big [ -\frac{2|x|e^{-\frac{|x|^2}{4(1+\epsilon t)}}}{(1+\epsilon t)^{1/2}}+ 2\mathsf{Erf}(\frac{|x|}{2(1+\epsilon t)^{1/2}}) \big ] \Big ] \Big ]  \] 
	
\section{Analysis of the linear evolution}
	
	\label{app:LinEvol}
	
	We show that $\rho_L(t), a_L(t)$ and $\vec{\omega}_L(t)$ defined in defined by (\ref{RhoALinearEvolution}) and (\ref{OmegaLinearEvolution}) for $t > 0$ and $(\rho_L(t), a_L(t), \vec{\omega}_L(t))^T = (\rho_0,a_0,\vec{\omega}_0)^T$ for $t = 0$ map $[0,\infty)$ continuously into $L^p(n)$ for initial conditions in $L^p(n)$, and that these define differentiable functions of space and time for $t > 0$.  We also determine bounds on the temporal evolution of the norms of these terms. 
	
\subsection{Smoothness properties}

	\begin{Proposition}\label{p:LinSmoothness}
		\textbf{(a)} Let $n \in \mathbb{R}_{\geq 0}$, $p\geq 1$ and $(\rho_0,a_0,\vec{\omega})^T \in W^{1,p}(n)\times L^p(n) \times \mathbb{L}_{\sigma}^{p}(n)$.  Then 
		\[ (\rho_L(t), a_L(t), \vec{\omega}_L(t))^T \in C^0\big [[0,\infty),L^p(n)\times L^p(n)\times \mathbb{L}_{\sigma}^p(n) \big ] \]
		
		\textbf{(b)} Let $n \in \mathbb{R}_{\geq 0}$ and $(\rho_0,a_0,\vec{\omega})^T \in W^{1,1}(n)\times L^1(n) \times \mathbb{L}_{\sigma}^{1}(n)$.  Then
		\[ (\partial_x^{\alpha}\rho_L(t), \partial_x^{\alpha} a_L(t), \partial_x^{\alpha} \vec{\omega}_L(t))^T \in C^0\big [(0,\infty),L^p(n)\times L^p(n)\times \mathbb{L}_{\sigma}^p(n) \big ] \]
		for every $1\leq p \leq \infty$ and $\alpha \in \mathbb{N}^3$. 
		
	\end{Proposition}
	
	\begin{proof}
		
		We prove continuity at $t=0$ for part \textbf{(a)}, then prove part \textbf{(b)}, and the continuity for $t>0$ follows from the fact that solutions are differentiable in time, and that these time derivatives can be written in terms of the spatial derivatives by virtue of the differential equation that the solutions satisfy.  Starting with $\vec{\omega}_L$ we show continuity at $t=0$ by first noting that it suffices to consider $\vec{\omega}_0$ which is smooth and has compact support by a density argument, together with the linearity of the heat operator, Young's inequality and the heat estimates in Proposition $\ref{p:HeatEst}$.  Standard arguments show that for such $\vec{\omega}_0$ we have $\mathbb{K}_{\epsilon}(t)\ast \vec{\omega}_0\to \vec{\omega}_0$ uniformly as $t\to 0$, and the result follows.  For $t > 0$ one obtains $\partial_x^{\alpha} \mathbb{K}_{\epsilon}(t)\ast \vec{\omega}_0 \in \mathbb{L}^p(n)$ via Young's inequality and the differentiability as a map into $\mathbb{L}^p_{\sigma}$ follows from the fact that
		\[ \lim_{h \to 0} \big \| \frac{\partial_x^{\alpha}K_{\epsilon}(t+h)-\partial_x^{\alpha} K_{\epsilon}(t)}{h} - \partial_t \partial_x^{\alpha} K_{\epsilon}(t) \big \|_{L^1(\mu)} =0 \]
		for all $\mu$, together with Young's inequality.
		
		For $\rho_L(t)$ we start with $\partial_t w(t) \ast K_{\nu}(t) \ast \rho_0$.  Again we can assume that $\rho_0$ is smooth and has compact support using Proposition $\ref{p:HypParaEst}$.  For such $\rho_0$ the uniform convergence of $\partial_t w(t) \ast K_{\nu}(t)\ast \rho_0$ to $\rho_0$ as $t\to 0$ is immediate from the formula 
		\[ \partial_{t}w(t) \ast K_{\nu}(t)\ast \rho_0 = \frac{1}{4\pi} \int \nolimits_{|z|=1} K_{\nu}(t)\ast \rho_0(x+ctz)dS(z) \]
		and from the result for $K_{\nu}(t)\ast \rho_0$.  The continuity in $L^p(n)$ then follows.  For $t > 0$ the differentiability follows by the same reasoning as above.  The proofs for the smoothness properties of the other terms are similar.
	\end{proof}
	
\subsection{Linear evolution decay rates}
	
	Let $r_{\alpha,p}, \ell_{n,p,\mu}$ and $\tilde{\ell}_{n,p,\mu}$ be as defined in $(\ref{LinDecayRates1})$ and $(\ref{LinDecayRates2})$.  
	
	\begin{Proposition}\label{p:LinDecayRates}
		Let $n \in \mathbb{R}_{\geq 0}$ be given.  Suppose $(\rho_0,a_0,\vec{\omega}_0)^T \in \bigcap \limits_{1\leq \tilde{p} \leq 3/2} W^{1,\tilde{p}}(n) \times L^{\tilde{p}}(n) \times \mathbb{L}_{\sigma}^{\tilde{p}}(n)$.  If $n>0$, suppose also that $a_0$ and $\vec{\omega}_0$ have zero total mass.  Then
		\begin{equation} \label{LinearDecayRates} \begin{split}
		\|\partial_x^{\alpha} \rho_L(t)\|_{\mathring{L}^p(\mu)} & \leq C t^{-r_{\alpha,p}} (1+t)^{-\ell_{n,p,\mu}+\frac{1}{2}} \sup_{1\leq \tilde{p} \leq 3/2 } \big ( \|\rho_0 \|_{W^{1,\tilde{p}}(n)} + \|a_0 \|_{L^{\tilde{p}}(n)}  \big )  \\
		\|\partial_x^{\alpha} a_L(t)\|_{\mathring{L}^p(\mu)} & \leq C t^{-r_{\alpha,p}} (1+t)^{-\ell_{n,p,\mu}} \sup_{1\leq \tilde{p} \leq 3/2 } \big ( \|\rho_0 \|_{W^{1,\tilde{p}}(n)} + \|a_0 \|_{L^{\tilde{p}}(n)}  \big ) \\
		\|\partial_x^{\alpha} \vec{\omega}_L(t)\|_{\mathring{\mathbb{L}}^p(\mu)} & \leq C t^{-r_{\alpha,p}} (1+t)^{-\tilde{\ell}_{n,p,\mu} } \sup_{1\leq \tilde{p} \leq 3/2  } \big ( \|\vec{\omega}_0 \|_{\mathbb{L}^{\tilde{p}}(n)} \big ) 
		\end{split} \end{equation}
		holds for all $t \in (0,\infty)$, $1 \leq p \leq \infty$, $0 \leq \mu \leq n$ and $\alpha \in \mathbb{N}^3$
	\end{Proposition}
	
	\begin{proof}
		
		In the following computations we ignore constant proportionality factors for simplicity.  The proof follows from Young's inequality, together with the the fact that we can split the weight via $(1+|x|)^\mu \leq (1+|y|)^\mu + (1+|x-y|)^\mu$ and estimate in different $L^p$ norms.  For the first term in (\ref{RhoALinearEvolution}), this is as follows.  For large times $t > 1$ we have 
		\[ \begin{split}
		\|\partial_t w\ast \partial_x^{\alpha} K_{\nu} \ast \rho_0 \|_{\mathring{L}^p(\mu)} & \leq \|\partial_t w \ast \partial_x^{\alpha} K_{\nu}(t) \|_{\mathring{L}^p(\mu)} \|\rho_0\|_{L^1} + \|\partial_t w \ast \partial_x^{\alpha} K_{\nu}(t) \|_{L^p} \|\rho_0\|_{\mathring{L}^1(\mu)} \\ & \leq t^{\frac{\mu}{2}-\frac{3}{2}(1-\frac{1}{p})-\frac{|\alpha|}{2}}(1+t)^{\frac{\mu}{2}+\frac{1}{2}-(1-\frac{1}{p})}\|\rho_0\|_{L^1} + t^{-\frac{3}{2}(1-\frac{1}{p})-\frac{|\alpha|}{2}}(1+t)^{\frac{1}{2}-(1-\frac{1}{p})}\|\rho_0\|_{\mathring{L}^1(\mu)}
		\end{split} \]
		whereas for small times $t < 1$ we have
		\[ \begin{split}
		\|\partial_t w\ast \partial_x^{\alpha} K_{\nu} \ast \rho_0 \|_{\mathring{L}^p(\mu)} & \leq \|\partial_t w \ast \partial_x^{\alpha} K_{\nu}(t) \|_{\mathring{L}^{\tilde{p}}(\mu)} \|\rho_0\|_{L^{3/2}} + \|\partial_t w \ast \partial_x^{\alpha} K_{\nu}(t) \|_{L^{\tilde{p}}} \|\rho_0\|_{\mathring{L}^{3/2}(\mu)} \\ & \leq t^{\frac{\mu}{2}-\frac{3}{2}(\frac{2}{3}-\frac{1}{p})-\frac{|\alpha|}{2}}(1+t)^{\frac{\mu}{2}+\frac{1}{2}-(\frac{2}{3}-\frac{1}{p})}\|\rho_0\|_{L^{3/2}} + t^{-\frac{3}{2}(\frac{2}{3}-\frac{1}{p})-\frac{|\alpha|}{2}}(1+t)^{\frac{1}{2}-(\frac{2}{3}-\frac{1}{p})}\|\rho_0\|_{\mathring{L}^{3/2}(\mu)}
		\end{split} \]
		for $p\geq 3/2$ and 
		\[ \begin{split}
		\|\partial_t w\ast \partial_x^{\alpha} K_{\nu} \ast \rho_0 \|_{\mathring{L}^p(\mu)} & \leq \|\partial_t w \ast \partial_x^{\alpha} K_{\nu}(t) \|_{\mathring{L}^1(\mu)} \|\rho_0\|_{L^p} + \|\partial_t w \ast \partial_x^{\alpha} K_{\nu}(t) \|_{L^{1}} \|\rho_0\|_{\mathring{L}^{p}(\mu)} \\ & \leq t^{\frac{\mu}{2}-\frac{|\alpha|}{2}}(1+t)^{\frac{\mu}{2}+\frac{1}{2}}\|\rho_0\|_{L^{p}} + t^{-\frac{|\alpha|}{2}}(1+t)^{\frac{1}{2}}\|\rho_0\|_{\mathring{L}^{p}(\mu)}
		\end{split} \]
		for $1\leq p \leq 3/2$ hence these norms blow up at the rate $t^{-\frac{3}{2}(\frac{2}{3}-\frac{1}{p})-\frac{|\alpha|}{2}}$ as $t\to 0$ for $p \geq 3/2$, blow up at the rate $t^{-\frac{|\alpha|}{2}}$ as $t\to 0$ for $1\leq p \leq 3/2$ and decay at the rate $t^{k-\frac{5}{2}(1-\frac{1}{p})+\frac{1}{2}-\frac{|\alpha|}{2}}$ as $t \to \infty$ for all $1 \leq p \leq \infty$.  For the next term in $\rho_L$ we find
		\[ \begin{split}
		\| w\ast \partial_x^{\alpha} & K_{\nu} \ast a_0 \|_{\mathring{L}^p(\mu)}  \leq \| w \ast \partial_x^{\alpha} K_{\nu}(\frac{t}{2}) \|_{\mathring{L}^{p}(\mu)} \|K_{\nu}(\frac{t}{2})\ast a_0\|_{L^{1}} + \| w \ast \partial_x^{\alpha} K_{\nu}(\frac{t}{2}) \|_{L^{p}} \|K_{\nu}(\frac{t}{2})\ast a_0\|_{\mathring{L}^{1}(\mu)} \\ & \leq t^{1+\frac{\mu}{2}-\frac{3}{2}(1-\frac{1}{p})-\frac{|\alpha|}{2}}(1+t)^{\frac{\mu}{2}-(1-\frac{1}{p})-\frac{\lfloor n \rfloor_1 }{2}}\|a_0\|_{\mathring{L}^{1}(\lfloor n \rfloor_1)} + t^{1-\frac{3}{2}(1-\frac{1}{p})-\frac{|\alpha|}{2}}(1+t)^{-(1-\frac{1}{p})}\|K_{\nu}(\frac{t}{2})\ast a_0\|_{\mathring{L}^{1}(\mu)}
		\end{split} \]
		for large times.  For the case $\mu=0$ note that the second term on the right hand side does not appear since we can use Young's inequality directly, and if $0 < \mu \leq n$ then we can use 
		\[ \|K_{\nu}(\frac{t}{2}) \ast a_0\|_{\mathring{L}^{1}(\mu)} 
		\leq t^{-\frac{\lfloor n \rfloor_1-\lfloor \mu \rfloor_1}{2}} \|a_0\|_{L^{1}(n)} \]
		For small times we have
		\[ \begin{split}
		\| w\ast \partial_x^{\alpha} K_{\nu} \ast a_0 \|_{\mathring{L}^p(\mu)} & \leq \| w \ast \partial_x^{\alpha} K_{\nu}(t) \|_{\mathring{L}^{\tilde{p}}(\mu)} \|a_0\|_{L^{3/2}} + \| w \ast \partial_x^{\alpha} K_{\nu}(t) \|_{L^{\tilde{p}}} \|a_0\|_{\mathring{L}^{3/2}(\mu)} \\ & \leq t^{1+\frac{\mu-|\alpha|}{2}-\frac{3}{2}(\frac{2}{3}-\frac{1}{p})}(1+t)^{\frac{\mu}{2}-(\frac{2}{3}-\frac{1}{p})}\|a_0\|_{L^{3/2}} + t^{1-\frac{3}{2}(\frac{2}{3}-\frac{1}{p})-\frac{|\alpha|}{2}}(1+t)^{-(\frac{2}{3}-\frac{1}{p})}\|a_0\|_{\mathring{L}^{3/2}(\mu)}
		\end{split} \]
		for $p\geq 3/2$ and 
		\[ \begin{split}
		\|w\ast \partial_x^{\alpha} K_{\nu} \ast a_0 \|_{\mathring{L}^p(\mu)} & \leq \| w \ast \partial_x^{\alpha} K_{\nu}(t) \|_{\mathring{L}^1(\mu)} \|a_0\|_{L^p} + \| w \ast \partial_x^{\alpha} K_{\nu}(t) \|_{L^{1}} \|a_0\|_{\mathring{L}^{p}(\mu)} \\ & \leq t^{1+\frac{\mu-|\alpha|}{2}}(1+t)^{\frac{\mu}{2}}\|\rho_0\|_{L^{p}} + t^{1-\frac{|\alpha|}{2}}\|\rho_0\|_{\mathring{L}^{p}(\mu)}
		\end{split} \]
		for $1\leq p \leq 3/2$.  The time estimates of the other terms in (\ref{RhoALinearEvolution}), (\ref{OmegaLinearEvolution}) are obtained similarly.  Note the weighted estimates in (\ref{LinearDecayRates}) aren't sharp for $\vec{\omega}_L$, but instead match the decay rate of solutions of (\ref{IntegralForm}).
	\end{proof}

\end{document}